\documentclass{elsarticle} 

\makeatletter
	\def\ps@pprintTitle{%
 	\let\@oddhead\@empty
	\let\@evenhead\@empty
	\def\@oddfoot{\centerline{\thepage}}%
	\let\@evenfoot\@oddfoot}
\makeatother

\usepackage{etoolbox}
\patchcmd{\MaketitleBox}{\footnotesize\itshape\elsaddress\par\vskip36pt}{\footnotesize\itshape\elsaddress\par\parbox[b][36pt]{\linewidth}{\vfill\hfill\textnormal{\today}\hfill\null\vfill}}{}{}%
\patchcmd{\pprintMaketitle}{\footnotesize\itshape\elsaddress\par\vskip36pt}{\footnotesize\itshape\elsaddress\par\parbox[b][36pt]{\linewidth}{\vfill\hfill\textnormal{\today}\hfill\null\vfill}}{}{}%
\usepackage{tikz}
\usepackage{tikz-imagelabels}
\usepackage{pgfplots} 
\usepackage{pgfgantt}
\usepackage{pdflscape}
\pgfplotsset{compat=newest} 
\pgfplotsset{plot coordinates/math parser=false} 
\newlength\fwidth
\newlength\fheight
\usepackage[
            colorlinks=true,%
            breaklinks=true,%
            linkcolor=blue,
            urlcolor=blue,%
            citecolor=blue,%
            pdftitle={Symplectic model reduction of Hamiltonian systems using data-driven quadratic manifolds}, 
            pdfkeywords={Symplectic model reduction, Hamiltonian systems, Data-driven modeling, Quadratic manifolds, Scientific machine learning},	
            pdfauthor={Authors},		
            bookmarksopen=false,
            pdfpagemode=UseNone]{hyperref}
            
\usepackage[margin = 1.2in]{geometry}
\usepackage[T1]{fontenc}
\usepackage[english]{babel}
\usepackage[utf8]{inputenc} 
\usepackage{mathtools}
\usepackage{amsmath}
\usepackage{amsfonts}
\usepackage{amsthm}
\usepackage{amssymb}
\usepackage{array}
\usepackage{algorithm} 
\usepackage{algorithmicx}
\usepackage{algpseudocode}
\usepackage{subcaption}
\usepackage{siunitx}

\usepackage{color}
\usepackage{url}
\usepackage{cleveref}
\usepackage{graphicx}
\usepackage{breakcites}
\usepackage{soul} 
\usepackage{enumitem}
\usepackage[title,titletoc,toc]{appendix}

\newcommand{\Ra}[1]{\color{black}{#1}}  
\newcommand{\Rb}[1]{\color{black}{#1}}        
\newcommand{\Rc}[1]{\color{black}{#1}}         

\usepackage[textsize=scriptsize]{todonotes}

\newtheorem{theorem}{Theorem}

\newtheorem{remark}{Remark}

\newtheorem{proposition}{Proposition}

\newtheorem{definition}{Definition}
\newtheorem{lemma}{Lemma}

\renewenvironment{proof}%
{\noindent {\textbf{Proof}:} }%
{\hfill $\Box$ \\[1ex] }

\DeclareMathOperator{\colspan}{colspan}

\newcommand{\bit}{\begin{itemize}}
\newcommand{\eit}{\end{itemize}}
\newcommand{\ben}{\begin{enumerate}}
\newcommand{\een}{\end{enumerate}}

\newcommand {\real} {\mathbb{R}}

\newcommand{\dH}{ H}
\newcommand{\In}{\mathbf I_n}
\newcommand{\Ir}{\mathbf I_r}
\newcommand{\Jn}{\mathbf J_{2n}}
\newcommand{\Jr}{\mathbf J_{2r}}
\newcommand{\Vt}{\bV^{+}}

\newcommand{\inv}[1]{#1^{-1}}
\newcommand{\invb}[1]{\inv{\left( #1 \right)}}

\DeclareMathOperator*{\argmin}{arg\,min}%
\DeclareMathOperator*{\skewPart}{skew}

\newcommand{\kronFun}{\bB}

\newcommand{\paramSet}{\mathcal{P}}

\newcommand{\ns}{n_{\mathrm{s}}}
\newcommand{\bYe}{\bY_{\mathrm e}}

\newcommand{\Sp}{\mathrm{Sp}}

\newcommand{\state}{\by}
\newcommand{\stateq}{\bq}
\newcommand{\statep}{\bp}

\newcommand{\bgammaQ}{\bGamma_{\mathrm{Q}}}
\newcommand{\Dq}{\bD_{\stateqRed}}
\newcommand{\MCL}{\mathrm{MCL}}
\newcommand{\LSL}{\mathrm{LSL}}
\newcommand{\qMCL}{\mathrm{QMCL}}
\newcommand{\bgammaMCL}{\bGamma_{\MCL}}
\newcommand{\bgammaLSL}{\bGamma_{\LSL}}
\newcommand{\bgammaQMCL}{\bGamma_{\qMCL}}
\newcommand{\bqgammaQMCL}{\bGamma_{\qMCL,q}}
\newcommand{\bpgammaQMCL}{\bGamma_{\qMCL,p}}
\newcommand{\bqgammaMCL}{\bGamma_{\MCL,q}}
\newcommand{\bpgammaMCL}{\bGamma_{\MCL,p}}
\newcommand{\vRed}{\widetilde{\bv}}
\newcommand{\stateV}{\bc_{\bV_2}}
\newcommand{\stateVbar}{\bc_{\Vbar_2}}
\newcommand{\VV}{\bV_2}
\newcommand{\VbarVbar}{\Vbar_2}
\newcommand{\bHVV}{\widetilde{\bH}_{\bV\bV}}
\newcommand{\bHVVbar}{\widetilde{\bH}_{\bV\Vbar}}
\newcommand{\bHVyref}{\widetilde{\bH}_{\bV\yref}}
\newcommand{\bHVbarV}{\widetilde{\bH}_{\Vbar\bV}}
\newcommand{\bHVbarVbar}{\widetilde{\bH}_{\Vbar\Vbar}}
\newcommand{\bHVbaryref}{\widetilde{\bH}_{\Vbar\yref}}
\newcommand{\tbHq}{\widetilde{\bH}_{q}}
\newcommand{\tbHp}{\widetilde{\bH}_{p}}
\newcommand{\tbHqbar}{\widetilde{\overline{\bH}}_{q}}
\newcommand{\tbHpbar}{\widetilde{\overline{\bH}}_{p}}
\newcommand{\tbHrefq}{\widetilde{\bH}_{\mathrm{ref},q}}
\newcommand{\tbHrefp}{\widetilde{\bH}_{\mathrm{ref},p}}
\newcommand{\eval}[2][\state]{\left.#2\right|_{#1}}

\newcommand{\rT}[1]{#1^\top}
\newcommand{\rTb}[1]{\rT{\left( #1 \right)}}

\newcommand{\reduce}[1]{\widetilde{#1}}
\newcommand{\stateRed}{\reduce{\state}}
\newcommand{\stateqRed}{\reduce{\stateq}}
\newcommand{\statepRed}{\reduce{\statep}}
\newcommand{\mpinv}[1]{#1^\dagger}

\newcommand{\dHRed}{ \reduce{H}_{{\rm d}}}
\newcommand{\Vq}{\bV_{\stateq}}
\newcommand{\Vqbar}{\overline{\bV}_{\stateq}}

\newcommand{\bzero}{\ensuremath{\mathbf{0}}} 

\newcommand{\bA}{\ensuremath{\mathbf{A}}}
\newcommand{\bB}{\ensuremath{\mathbf{B}}}

\newcommand{\bD}{\ensuremath{\mathsf{D}}}

\newcommand{\bH}{\ensuremath{\mathbf{H}}}
\newcommand{\bI}{\ensuremath{\mathbf{I}}}

\newcommand{\bP}{\ensuremath{\mathbf{P}}}
\newcommand{\bQ}{\ensuremath{\mathbf{Q}}}

\newcommand{\bV}{\ensuremath{\mathbf{V}}}

\newcommand{\bY}{\ensuremath{\mathbf{Y}}}

\newcommand{\bc}{\ensuremath{\mathbf{c}}}

\renewcommand{\bf}{\ensuremath{\mathbf{f}}}

\newcommand{\bp}{\ensuremath{\mathbf{p}}}
\newcommand{\bq}{\ensuremath{\mathbf{q}}}
\newcommand{\br}{\ensuremath{\mathbf{r}}}

\newcommand{\bv}{\ensuremath{\mathbf{v}}}

\newcommand{\bx}{\ensuremath{\mathbf{x}}}
\newcommand{\by}{\ensuremath{\mathbf{y}}}

\newcommand {\bGamma} {\boldsymbol{\Gamma}}

\newcommand {\bPhi} {\mbox{\boldmath $\Phi$}}

\newcommand{\bmu}{\ensuremath{\boldsymbol{\mu}}}

\newcommand{\cC}{\ensuremath{\mathcal{C}}}

\newcommand{\cH}{\ensuremath{\mathcal{H}}}

\newcommand{\cO}{\ensuremath{\mathcal{O}}}
\newcommand{\bigO}{\cO}
\newcommand{\cP}{\ensuremath{\mathcal{P}}}

\newcommand{\Vbar}{\overline{\mathbf V}}
\newcommand{\bqgamma}{\ensuremath{\mathbf{\Gamma}}}
\newcommand {\yref} {\by_{\mathrm{ref}}}
\newcommand {\qref} {\bq_{\mathrm{ref}}}
\newcommand {\pref} {\bp_{\mathrm{ref}}}
\newcommand{\BQ}{\mathrm{BQ}}
\newcommand{\bgammaBQ}{\bGamma_{\BQ}}
\newcommand{\bqgammaBQ}{\bGamma_{\BQ,q}}
\newcommand{\bpgammaBQ}{\bGamma_{\BQ,p}}

\graphicspath{{./figures/}}

\begin{document}
	
\begin{frontmatter}
    
    \title{Symplectic model reduction of Hamiltonian systems using data-driven quadratic manifolds}
		
		
 		\author[ucsd]{Harsh Sharma\corref{cor1}}
		 		\cortext[cor1]{Corresponding author}
		\ead{hasharma@ucsd.edu}
            \author[affil3]{Hongliang Mu}
            \author[affil4]{Patrick Buchfink}
 		\author[affil2]{Rudy Geelen}
 		\author[affil3]{Silke Glas}

 		\author[ucsd]{Boris Kramer}

			\address[ucsd]{Department of Mechanical and Aerospace Engineering, University of California San Diego, CA, United States}
		
 		\address[affil3] {Department of Applied Mathematics, University of Twente, Enschede, The Netherlands}
				
 		\address[affil4] {Institute of Applied Analysis and Numerical Simulation, University of Stuttgart, Stuttgart, Germany}
 		\address[affil2] {Oden Institute for Computational Engineering and Sciences, University of Texas at Austin, TX, United States}

\date{}
		
\begin{abstract}
This work presents two novel approaches for the symplectic model reduction of high-dimensional Hamiltonian systems using data-driven quadratic manifolds. Classical symplectic model reduction approaches employ linear symplectic subspaces for representing the high-dimensional system states in a reduced-dimensional coordinate system. While these approximations respect the symplectic nature of Hamiltonian systems, {\Rb{linear basis approximations can suffer from slowly decaying Kolmogorov $N$-width, especially in wave-type problems, which then requires a large basis size.}} We propose two different model reduction methods based on recently developed quadratic manifolds, each presenting its own advantages and limitations. The addition of quadratic terms {\Rb{to}} the state approximation, which sits at the heart of the proposed methodologies, enables us to better represent intrinsic low-dimensionality in the problem at hand. Both approaches are effective for issuing predictions in settings well outside the range of their training data while providing more accurate solutions than the linear symplectic reduced-order models.
\end{abstract}	


\end{frontmatter}

\section{Introduction} \label{sec:intro}
Computational modeling, simulation, and control of dynamical systems characterized by Hamiltonian mechanics are essential for many science and engineering applications such as robotics, quantum mechanics, solid state physics, and climate modeling, see~\cite{leimkuhler_reich_2005,marsden2013introduction}. Hamiltonian partial differential equations (PDEs) demonstrate complex dynamic behavior while possessing underlying mathematical structures in the form of symmetries, Casimirs, symplecticity, and energy conservation. The structure-preserving spatial discretization of Hamiltonian PDEs with finite element methods, spectral methods, or finite difference methods leads to finite-dimensional Hamiltonian systems with large state-space dimensions. 

For a high-dimensional Hamiltonian system, it is desirable to construct a reduced-order model (ROM) that preserves the underlying geometric structure and in particular, the stability of the full-order model (FOM). To this end, the field of structure-preserving model reduction has developed efficient tools for reducing the computational complexity of high-dimensional FOMs while preserving the underlying geometric structure. Symplectic model reduction of Hamiltonian systems was introduced in \cite{peng2016symplectic}, where the Galerkin projection-based ROM was modified so that the ROM retains the underlying symplectic structure. In another research direction, the work in \cite{afkham2017structure} presented a greedy approach for symplectic model reduction of parametric Hamiltonian systems. A variety of symplectic model reduction approaches have been developed along these two research directions recently, e.g., basis generation for symplectic model reduction of Hamiltonian systems~\cite{bendokat2022,buchfink2020a,buchfink2019symplectic,Buchfink2022}, modified projection techniques~\cite{gong2017structure,hesthaven2022rank,pagliantini2021dynamical}, and data-driven operator inference~\cite{sharma2022hamiltonian,sharma2022preserving,gruber2023canonical}. However, all of the aforementioned methods for symplectic model reduction project onto linear subspaces{\Rb{ which often requires simulating high-dimensional ROMs to obtain accurate approximate solutions.}} 

Classical projection-based model reduction methods encounter substantial difficulties whenever the solution manifold has so-called Kolmogorov $N$-widths that decay slowly with increasing $N$ \cite{kolmogoroff1936uber}. The Kolmogorov $N$-width quantifies the degree of accuracy by which a set can be approximated using linear subspaces of dimension $N$. If the $N$-widths decay slowly with increasing $N$, the associated ROMs can provide accurate approximate solutions only for large values of $N$, which in turn defeats the original purpose of achieving computational speedups. This poses substantial limitations to the spectrum of possible applications. Many physics-based systems do not exhibit a global low-rank structure and are characterized by slowly decaying Kolmogorov $N$-widths, such as those featuring coherent structures that propagate in space-time \cite{greif2019decay, ohlberger2016reduced}. Strategies aimed at overcoming this bottleneck typically fall into one of three categories: (1) dictionary approaches (local reduced basis methods, divide-and-conquer, \cite{https://doi.org/10.1002/nme.4371, geelen2022localized, peherstorfer2014localized, daniel2020model, DANIEL2022111120}), (2) time-dependent basis or other transformations \cite{ISSAN2023111689, papapicco2022neural, reiss2018shifted, nonino2019overcoming, ohlberger2013nonlinear, taddei2015reduced, koch2007dynamical, hesthaven2022rank, sapsis2009dynamically}, (3) and the use of nonlinear approximation techniques (deep learning and other data-driven manifold techniques) \cite{fresca2021comprehensive, FRESCA2022114181, GAO2020132614, gonzalez2018deep, SAN2019271, pawar2019deep, ahmed2021nonlinear,barnett2022neural}, which are sometimes referred as model reduction on manifolds. In the present work, we focus on the latter class of methods. Following the pioneering work from \cite{LEE2020108973}, projection-based reduced-order modeling approaches based on convolutional autoencoders have seen a surge of activity in recent years \cite{fresca2021comprehensive, FRESCA2022114181, KIM2022110841, KADEETHUM2022104098, ahmed2021nonlinear}. With various levels of success, these methods can produce ROMs of lower dimensionality---compared to classical approaches based on linear subspaces---for the same level of accuracy. In general, however, nonlinear model reduction methods are not tailored to preserve the underlying geometric structure and are thus ill-suited for Hamiltonian systems. 

An exception in that regard is the work \cite{buchfink2021symplecticmanifold}, which uses a nonlinear approximation while preserving the structure of Hamiltonian systems.
There, the symplectic manifold Galerkin (SMG) projection is introduced which generalizes the structure-preserving model reduction with linear subspaces from \cite{peng2016symplectic} to symplectic (nonlinear) approximation mappings.
The SMG projection requires a symplectic approximation mapping. While the SMG projection works for all symplectic approximation mappings, the numerical example in~\cite{buchfink2021symplecticmanifold} however only considers a weakly symplectic mapping. Moreover, the mapping is based on a deep convolutional autoencoder which leads to a ROM that lacks interpretability. Designing autoencoders requires careful considerations towards reproducible research including the choice of network size and depth, activation function, optimizer, and weight initialization. Each physical problem usually requires its own hyperparameter tuning. Furthermore, different latent space sizes might lead to different optimal hyperparameter configurations, rendering the training process tedious and often unstable. Solutions to this problem have been proposed in the form of automated parameter choices, but more principled techniques have yet to be established.

In a similar direction, data-driven quadratic manifolds have emerged to build ROMs for transport-dominated problems in~\cite{BARNETT2022111348, GEELEN2023115717}. Combining quadratic manifolds with Galerkin projection leads to a highly scalable approach for deriving efficient, yet interpretable, ROMs. However, the application of these approaches to Hamiltonian systems, in their standard formulation, leads to a violation of the underlying geometric structure. This makes the associated ROMs prone to displaying nonphysical behavior (see Section~\ref{sec:motivation}). The main goal of this work is to build on the recent successes of data-driven quadratic manifolds to construct structure-preserving state approximations, characterized by quadratic nonlinearities, for symplectic model reduction. We focus particularly on high-dimensional dynamical-system models that arise from the semi-discretization of Hamiltonian PDEs. The key contributions of this paper are:
     \begin{enumerate}
        \item We propose the quadratic manifold cotangent lift (QMCL) and the resulting SMG-QMCL-ROM in Section \ref{sec:symp_quad_jacobian} that use a quadratic approximation to build a symplectic approximation mapping while using the SMG as projection. The resulting ROM is a Hamiltonian system of reduced dimension. Thus, the SMG-QMCL is strictly preserving the structure of Hamiltonian systems throughout the reduction.
        \item We propose the blockwise-quadratic (BQ) approximation mapping and the resulting Galerkin-BQ-ROM in Section~\ref{sec:symp_quad_pod} that increases the accuracy of linear-subspace Hamiltonian ROMs by introducing quadratic terms to the state approximation. This approach augments the linear-subspace Hamiltonian ROM with interpretable and physics-preserving higher-order terms to derive an approximately Hamiltonian ROM with improved computational efficiency compared to the SMG-QMCL-ROM.
        \item Numerical experiments on the challenging cases of a thin moving pulse and {\Rb{a two-dimensional nonlinear wave equation}} illustrate that the proposed ROMs of the same dimensionality yield more accurate results than the symplectic linear subspace ROMs. In a parameter extrapolation study {\Rb{for the nonlinear example}}, the proposed methods yield accurate predictions of the high-dimensional state even for parameter values that are {\Rc{outside the range of the training parameters}}. In a time extrapolation study, the proposed ROMs provide long-time stable and accurate predictions outside the training time interval.
    \end{enumerate}
    
This paper is structured as follows. In Section~\ref{sec:background} we recall the main concepts and definitions of Hamiltonian systems and symplectic model reduction on linear subspaces. We also describe the data-driven construction of quadratic manifolds for nonlinear dimension reduction. Section~\ref{sec:method} presents the proposed SMG-QMCL-ROM and Galerkin-BQ-ROM for deriving ROMs of high-dimensional, Hamiltonian systems. In Section~\ref{sec:results} we present various numerical experiments for the {\Rb{parametric one-dimensional linear wave equation and two-dimensional nonlinear wave equation.}} Finally, Section~\ref{sec:conclusions} summarizes the contributions and suggests a number of interesting avenues for future research. 

\section{Background} \label{sec:background}

Section~\ref{sec:hamiltonian} provides an introduction to Hamiltonian PDE models, followed by details about their structure-preserving space discretization. In Section~\ref{sec:sympl_mor_lin_subspace} we review the basic principles of symplectic model reduction for Hamiltonian systems using linear symplectic subspaces. The data-driven learning of quadratic manifolds directly from high-dimensional data is recapitulated in Section~\ref{sec:quadratic}.
\subsection{Hamiltonian systems} \label{sec:hamiltonian}
In the following we consider parametric Hamiltonian systems{\Rb{.}} Let $\paramSet \subset \real^{n_{\bmu}}$ denote a parameter set for which we denote the $n_{\bmu}$-dimensional elements $\bmu \in \paramSet$ as the parameter vector. We consider infinite-dimensional Hamiltonian systems described by evolutionary PDEs of the form 
\begin{equation}
\frac{\partial y(x,t;\bmu)}{\partial t} = {\Ra{\mathcal{S}}} \frac{\delta \mathcal{H}}{\delta y}(y(x,t;\bmu);\bmu),
\label{eq:HPDE}
\end{equation}
where $x$ is the spatial coordinate, $t$ is time, $\mathcal{S}$ is a skew-symmetric operator, and $\delta \mathcal{H}/\delta y$ is the variational derivative\footnote{%
The variational derivative of $\mathcal{H}$ is defined through
$\frac{\text{d}}{\text{d}\varepsilon}\mathcal{H}[y+ \varepsilon v; \bmu]|_{\varepsilon=0}=\left\langle \frac{\delta \mathcal{H}}{\delta y}(y; \bmu),v\right\rangle$ where $v$ is an arbitrary function.
} of the Hamiltonian energy functional $\mathcal{H}$ with respect to the state variable $y$, and $\bmu \in \paramSet$ is a fixed but arbitrary parameter vector.

Hamiltonian PDEs possess important geometric properties in the form of symplecticity and conservation laws. These geometric properties are intimately related to the ability of the space-discretized FOMs to reproduce the long-time behavior of the solutions of Hamiltonian PDEs~\cite{bridges2006numerical,leimkuhler_reich_2005}. Therefore, the underlying symplectic structure and the conservative nature of Hamiltonian PDEs should be respected in their spatial discretization. Space-discretized Hamiltonian FOMs are finite-dimensional systems that are derived from the Hamiltonian PDE via structure-preserving semi-discretization. Contrary to the traditional approach in which the governing PDEs are discretized directly, Hamiltonian FOMs are obtained through a discretization of the space-time continuous Hamiltonian functional $\mathcal{H}$ and the skew-symmetric operator $\mathcal{S}$ in space~\cite{bridges2006numerical,celledoni2012preserving}. The resulting Hamiltonian FOM is given by 
\begin{equation}\label{eq:noncan}
\dot{\by}(t;\bmu)
={\Ra{\mathbf{S}}} \nabla_{\by}H(\by(t;\bmu);\bmu), \qquad \by(0;\bmu)=\by_0(\bmu),
\end{equation}
where $\by(t;\bmu) \in \real^{2n}$ is the high-dimensional state vector, $\by_0:\mathcal{P}\mapsto \real^{2n}$ is the parametric initial condition, ${\Ra{\mathbf{S}=-\mathbf{S}^\top}}$ is the skew-symmetric matrix approximation to ${\Ra{\mathcal{S}}}$ in~\eqref{eq:HPDE}, and $H(\cdot; \bmu) \in \cC^1(\real^{2n})$ is the space-discretized parametric Hamiltonian function. In the above equation, and for the remainder of this paper, the dot notation denotes differentiation with respect to time.

We focus on high-dimensional Hamiltonian systems with canonical symplectic structure ${\Ra{\mathbf S}}=\Jn=\small \begin{pmatrix}
       \bzero & \In \\
       -\In & \bzero
\end{pmatrix}$, where $\In$ is the $n$-dimensional identity matrix. For canonical Hamiltonian systems, the state vector can be partitioned as $\by(t;\bmu)=(\bq(t;\bmu)^\top,\bp(t;\bmu)^\top)^\top \in \real^{2n}$ where $\bq(t;\bmu) \in \real^n$ is the generalized position vector, and $\bp(t;\bmu) \in \real^n$ is the generalized momentum vector. The governing equation for the canonical Hamiltonian FOM equals
      \begin{equation}
    \dot{\by}(t;\bmu)=\begin{pmatrix}
       \dot{\bq}(t;\bmu) \\
       \dot{\bp}(t;\bmu)
     \end{pmatrix}=\Jn \nabla_{ \by }H(\bq(t;\bmu),\bp(t;\bmu);\bmu), \qquad \by(0;\bmu)=\by_0(\bmu).
     \label{eq:HFOM}
\end{equation}
The flow map\footnote{For a dynamical system $\dot{\by}(t)=\mathbf{f}(\by(t))$  described by a differentiable vector field $\bf$ with initial value $\by(0) = \by_0$, the flow map $\varphi_t(\by_0) := \by(t)$ describes the solution as a function of the initial value.} $\varphi_t(\by(0;\bmu);\bmu)=\by(t;\bmu)$ for Hamiltonian FOMs \eqref{eq:HFOM} preserves a skew-symmetric, bilinear form known as the canonical symplectic form $\omega = \sum_{i=1}^n \text{d}p_i \wedge \text{d}q_i$ and conserves the Hamiltonian function $\dH$, i.e., $H(\bq(0;\bmu),\bp(0;\bmu);\bmu)=H(\bq(t;\bmu),\bp(t;\bmu);\bmu)$ for all $t$.

The field of geometric numerical integration~\cite{hairer2006geometric, leimkuhler_reich_2005,sharma2020review} has developed a wide variety of time integrators for Hamiltonian systems that are designed to respect some of the geometric features, such as energy, momentum, or the symplectic form. Structure-preserving integrators for finite-dimensional Hamiltonian systems can be divided into two categories: (i) symplectic methods~\cite{sanz1992symplectic,yoshida1990construction,mclachlan1993symplectic} that preserve the canonical symplectic form $\omega$ and (ii) energy-preserving integrators~\cite{gonzalez1996stability,celledoni2012preserving, sharma2022performance} that conserve the space-discretized Hamiltonian $\dH$ exactly. 
\subsection{Symplectic model reduction using linear subspaces} \label{sec:sympl_mor_lin_subspace}
The main goal in symplectic model reduction of high-dimensional Hamiltonian systems is the preservation of symplectic structure during the projection step. The methodologies described in Section \ref{sec:method} may be viewed as nonlinear extensions of the symplectic model reduction method from \cite{peng2016symplectic}, which covered linear subspaces and is briefly reviewed here.

To approximate the high-dimensional system state $\by(t;\bmu) \in \real^{2n}$ we seek an \emph{approximation mapping} $\bqgamma: \real^{2r}\mapsto \real^{2n}$, with $r \ll n$, as follows:
 \begin{equation}
\state(t; \bmu) \approx \bqgamma ( \stateRed ( t ; \bmu ) ) \in \real^{2n},
 \end{equation}
 where $\stateRed(t; \bmu)\in\real^{2r}$ denotes the reduced state vector of dimension $2r$ at time $t$ and parameter vector $\bmu$. We continue by defining the linear symplectic lift (LSL) as 
\begin{equation}
\bqgamma_{\text{LSL}} ( \stateRed ( t ; \bmu )) := \by_{\text{ref}} + \bV\stateRed(t; \bmu),
\label{eq:linear_approximation}
\end{equation}
where $\mathbf{y}_\text{ref}\in \real^{2n}$ is a reference state that is used for centering the training data and the basis matrix $\bV \in \real^{2n \times 2r}$ is a symplectic matrix, that is a matrix that satisfies 
\begin{equation}
\bV^\top\Jn\bV=\Jr.
\end{equation}
The set of all $2n \times 2r$ symplectic matrices is referred to as the symplectic Stiefel manifold and denoted with $\Sp(2r,\real^{2n})$. The symplectic inverse $\Vt \in \real^{2r \times 2n}$ of a symplectic matrix $\bV$ is defined by 
\begin{equation}\label{eq:sympl_inv}
\Vt:=\Jr^\top\bV^\top\Jn
\qquad\text{such that}\qquad
\Vt \bV = \bI_{2r} \in \real^{2r \times 2r},
\end{equation}
which allows one to define the \emph{symplectic Galerkin projection} with 
$$\stateRed(t;\bmu)=\Vt\left(\by(t;\bmu)-\by_{\text{ref}}\right).$$ Assuming $\by(t;\bmu)-\by_{\text{ref}} \in \colspan(\bV)$ for all $t$, the time evolution of the reduced state $\dot{\stateRed}(t;\bmu)$ is derived from the symplectic projection of the FOM \eqref{eq:HFOM} on the symplectic subspace via $\Vt$. This results in the ROM equations
\begin{equation} \label{eq:Hr_linear_subs}
\begin{aligned}
    \dot{\stateRed}(t;\bmu)&=\Vt\Jn \nabla_{ \by }\dH(\bqgamma_{\text{LSL}} ( \stateRed ( t ; \bmu )) ) \\
    &=\Jr\bV^\top\nabla_{\by}\dH(\bqgamma_{\text{LSL}} ( \stateRed ( t ; \bmu )) ) \\
    &=\Jr\nabla_{\stateRed}\dH(\bqgamma_{\text{LSL}} ( \stateRed ( t ; \bmu ))),
\end{aligned}
\end{equation}
where we used $\Vt \bV = \bI_{2r}$ on the left-hand side, $\bV^+\Jn=\Jr\bV^\top$ in the second step, and the chain rule
$\nabla_{\stateRed} \dH(\bqgamma_{\text{LSL}} ( \stateRed ( t ; \bmu )))
=\bV^\top \nabla_{\by} \dH(\bqgamma_{\text{LSL}} ( \stateRed ( t ; \bmu )))$ for $\stateRed(t;\bmu) \in \real^{2r}$ in the last step. Even if $\by(t;\bmu)-\by_{\text{ref}} \notin \colspan(\bV)$ for some $t$, the rightmost expression in~\eqref{eq:Hr_linear_subs} is still well-defined. 

By defining the reduced Hamiltonian $\widetilde{H}: \real^{2r}\mapsto \real$ as $\widetilde H(\stateRed(t;\bmu);\bmu):=\dH (\bqgamma_{\text{LSL}} ( \stateRed ( t ; \bmu )))$, we can rewrite~\eqref{eq:Hr_linear_subs} as a canonical Hamiltonian system. Thus, the symplectic Galerkin projection of a $2n$-dimensional Hamiltonian system \eqref{eq:HFOM} is given by a $2r$-dimensional \textit{{\Rc{LSL-ROM}}} with
\begin{equation}
\dot{\stateRed}(t;\bmu)=\Jr\nabla_{\stateRed}\widetilde{H}(\stateRed(t;\bmu);\bmu).
\label{eq:Hr_linear}
\end{equation}

The basis matrix can be computed using one of many snapshot-based basis generation techniques. Such techniques collect solutions of the Hamiltonian FOM (the snapshots) for different time instances and parameter vectors. The \emph{proper symplectic decomposition}~\cite{peng2016symplectic} is a snapshot-based basis generation method to find a symplectic basis matrix $\bV \in \real^{2n\times 2r}$. The basis is required to minimize the projection error of the symplectic projection in the mean over all snapshots which is characterized by the optimization problem
\begin{equation}
\min_{\substack{\bV\in \real^{2n\times 2r} \\ \text{s.t.} \ \bV^\top\Jn\bV=\Jr}} \left\|
    \bY_{\bmu} - \bV\Vt \bY_{\bmu}
\right\|_F,
\label{eq:PSD}
\end{equation}
where the centered snapshot data matrix $\bY_{\bmu}$ is defined as
\begin{equation}\label{eq:cotangent_snapshot}
    \bY_{\bmu} := \left(
        \bY_{\bmu_1}, \dots, \bY_{\bmu_M} \right)
    \in \mathbb{R}^{2n \times \ns},
\end{equation}
which stacks $\ns := MK$ centered snapshots in its columns for $M$ parameter instances $\bmu_1,\dots,\bmu_M$, each with $K$ time steps with
\begin{equation}
    \bY_{\bmu_j} :=
    \begin{pmatrix}
        | &  & | \\
        (\state(t_1; \bmu_j)-\by_{\text{ref}}) &  \dots & (\state(t_{K}; \bmu_j)-\by_{\text{ref}}) \\
        | &  & | \\
    \end{pmatrix} \in \mathbb{R}^{2n \times K}.
\end{equation}
For canonical Hamiltonian systems, the centered snapshot data matrix $\bY_{\bmu}$ can be partitioned as $\bY_{\bmu}=(\bQ_{\bmu}^\top,\bP_{\bmu}^\top)^\top$ with the (centered) generalized position and momentum data matrices $\bQ_{\bmu} \in \mathbb{R}^{n \times \ns}$ and $\bP_{\bmu} \in \mathbb{R}^{n \times \ns}$, respectively. We denote the columns of $\bY_{\bmu}$, $\bQ_{\bmu}$,  and $\bP_{\bmu}$ as $\by_j$, $\bq_j$, and $\bp_j$, respectively, with $j=1, \dots, \ns$ in the following. 

\begin{algorithm}[tbp]
\caption{Cotangent lift algorithm \cite{peng2016symplectic}}
\begin{algorithmic}[1]
\Require  Centered snapshot data matrix $\bY_{\bmu}=(\bQ_{\bmu}^\top,\bP_{\bmu}^\top)^\top \in \real^{2n \times \ns}$ arranged as in~\eqref{eq:cotangent_snapshot} and reduced dimension $r$
\Ensure Linear symplectic basis matrix $\bV \in \real^{2n \times 2r}$
\State Assemble the extended snapshot matrix $\bYe:=(\bQ_{\bmu},\bP_{\bmu})\in \real^{n \times 2\ns}$

\State Compute the SVD of $\bYe$
\State $\mathbf{\Phi} \in \mathbb{R}^{n \times r}$ $\leftarrow$ the $r$ leading left-singular vectors of $\bYe$
\State $\bV= \begin{pmatrix}
\mathbf{\Phi} & \bzero \\
\bzero & \mathbf{\Phi}
\end{pmatrix} \in \real^{2n \times 2r}$ with $\Vt=\bV^\top$ \Comment{Assemble linear symplectic basis matrix}
\end{algorithmic}
\label{alg:cotangent_lift}
\end{algorithm}

Since no general solution is known for \eqref{eq:PSD}, we seek to find solutions in subsets of $\Sp(2r,\real^{2n})$ by imposing further assumptions on $\bV$. In this paper the \emph{cotangent lift {\Ra{algorithm}}}\footnote{{\Ra{The cotangent lift algorithm~\cite{peng2016symplectic} is an SVD-based algorithm for constructing a symplectic basis matrix $\bV$. Unlike the cotangent lift of a smooth map defined for general nonlinear diffeomorphisms~\cite{abraham2008foundations}, the cotangent lift algorithm in~\cite{peng2016symplectic} only considers linear transformations.}}}, summarized in Algorithm \ref{alg:cotangent_lift}, is considered for computing a linear symplectic basis
\begin{align}\label{eq:basis_cotan_lift}
    \bV = \begin{pmatrix}
        \bPhi & \bzero\\
        \bzero & \bPhi
    \end{pmatrix} \in \real^{2n\times 2r}
\end{align}
of block-diagonal structure. {\Rb{For problems where $\bQ_{\bmu}$ and $\bP_{\bmu}$ have very disparate scales, one can scale the extended snapshot data matrix in Algorithm~\ref{alg:cotangent_lift}, i.e., $\bYe:=(\bQ_{\bmu},\gamma\bP_{\bmu})\in \real^{n \times 2\ns}$, before computing $\bPhi$ via SVD of the extended snapshot matrix, see~\cite{peng2016symplectic} for more details.}} Other symplectic basis generation techniques consider greedy basis generation \cite{afkham2017structure,buchfink2020a}, non-orthogonal, symplectic basis generation via the SVD-like decomposition \cite{buchfink2019symplectic}, iterative basis generation via optimization on manifolds \cite{bendokat2022}, or optimal symplectic basis generation under certain assumptions on the Hamiltonian system \cite{Buchfink2022}.

\subsection{Data-driven learning of quadratic state approximations}
\label{sec:quadratic}
Many datasets in problems involving high-dimensional state spaces are amenable to dimension reduction using linear techniques such as the POD. However, driving the projection error to an acceptable value can require a large reduced dimension $r$, rendering the associated ROMs inefficient. We therefore follow \cite{GEELEN2023115717} in introducing nonlinear state approximations with quadratic dependence on the reduced state vector as
\begin{equation}
\state(t; \bmu) \approx \bgammaQ ( \stateRed ( t ; \bmu )) := \by_{\text{ref}}+ \bV \stateRed(t; \bmu) + \overline{\bV} (\stateRed(t; \bmu) \otimes \stateRed(t; \bmu)),
\label{eq:quadratic_approximation}
\end{equation}
where $\bV\in \real^{2n \times 2r}$ is a basis matrix that has $\mathbf{v}_j$ as its $j$th column; $\overline{\bV} \in \real^{2n \times r(2r+1)}$ is a matrix whose columns are populated by the vectors $\overline{\mathbf{v}}_j$; and $\otimes$ denotes the Kronecker product without the redundant terms\footnote{For a column vector $\bx=[x_1,x_2, \dots, x_m]^\top \in \real^m$, the column-wise Kronecker product without the redundant terms is defined by
$ \bx \otimes \bx := \begin{bmatrix}
x_{1}^{2} & x_{1}x_{2} & \dots & x_1x_m  & x_{2}^{2} & x_2x_3 & \dots & x_2x_m & \dots & x_m^2\end{bmatrix}^{\top} \in \real^{m(m+1)/2}.$}. Learning the quadratic state approximation \eqref{eq:quadratic_approximation} from the centered snapshot set $\{ \state_j \}_{j=1}^{\ns} \subset \real^{2n}$ amounts to a representation learning problem in which the matrices $\bV$, $\overline{\bV}$, and the reduced-order state representations $\{ \stateRed_j \}_{j=1}^{\ns} \subset \real^{2r}$ are determined numerically through the solution of the optimization problem
\begin{equation}
    \argmin_{\bV, \overline{\bV}, \{ \stateRed_j \}_{j=1}^{\ns} } \left( J(\bV,\overline{\bV},\{ \stateRed_j \}_{j=1}^{\ns}) + \gamma \left\| \overline{\bV} \right\|_F^2 \right),
    \label{eq:regularized}
\end{equation}
where
\begin{equation}
    J(\bV,\overline{\bV},\{ \stateRed_j \}_{j=1}^{\ns}) = \sum_{j=1}^{\ns} \left\| \mathbf{y}_j  - \bV \stateRed_j - \overline{\bV} ( \stateRed_j  \otimes \stateRed_j ) \right\|_2^2.
    \label{eq:optim_problem2}
\end{equation}
A scalar regularization factor $\gamma \geq 0$ is used to avoid the overfitting of data.

We learn state approximations of the form \eqref{eq:quadratic_approximation} in a two-step fashion. First, the columns of the basis matrix $\bV$ are chosen as the left-singular vectors corresponding to the $2r$ largest singular values of the centered snapshot data matrix $\bY_{\bmu}$. The representation of the data in the reduced coordinate system can then be computed via the projection $\stateRed_j = \bV^\top\mathbf{y}_j$ for every centered snapshot $\mathbf{y}_j$ with $j=1,\dots,\ns$. This leaves only the matrix operator $\overline{\bV}$ to be inferred from the data. Given the basis matrix $\mathbf{V}$ and snapshot set $\{ \stateRed_j \}_{j=1}^{\ns}$, the optimization problem \eqref{eq:regularized} simplifies to
\begin{equation}
    \argmin_{\overline{\bV} } \left( J(\bV,\overline{\bV},\{ \stateRed_j \}_{j=1}^{\ns}) + \gamma \left\| \overline{\bV} \right\|_F^2 \right),
    \label{eq:regularized2}
\end{equation}
which is a linear least-squares problem with the explicit solution
\begin{equation}\label{eq:def_Vbar}
    \overline{\mathbf{V}} = (\mathbf{I} - \mathbf{V}\mathbf{V}^\top) \bY_{\text{\bmu}} \mathbf{W}^\top ( \mathbf{W} \mathbf{W}^\top + \gamma \mathbf{I})^{-1} \in \mathbb{R}^{2n \times r(2r+1)},
\end{equation}
where
\begin{equation}
    \mathbf{W} :=
    \begin{pmatrix}
    | &  & | \\
    \stateRed_1 \otimes \stateRed_1 & \dots & \stateRed_{\ns} \otimes \stateRed_{\ns}  \\
    | &  & |
    \end{pmatrix} \in \mathbb{R}^{r(2r+1) \times \ns}.
    \label{eq:quadratic_data}
\end{equation}
By construction of $\Vbar$ in \eqref{eq:def_Vbar}, each column of $\overline{\mathbf{V}}$ is in the column space of the orthogonal complement of $\mathbf{V}$ in $\real^{2n}$ so that the orthogonality condition $\bV^\top\overline{\mathbf{V}}=\bzero$ holds, see~\cite{GEELEN2023115717} for more details. The inferred basis matrix then approximates $\state(t; \bmu) \approx \bgammaQ ( \stateRed ( t ; \bmu ))$ in a least-squares sense. Equation~\eqref{eq:regularized} amounts to a regularized linear least-squares problem that can be solved for each row of $\overline{\bV}$ for a total of $2n$ such problems, each inferring $r(2r+1)$ unique entries.
\section{Symplectic model reduction of Hamiltonian
systems using quadratic manifolds} 
\label{sec:method}

We now present two different approaches that use, at their core, data-driven quadratic manifolds for building projection-based ROMs of Hamiltonian systems. Section~\ref{sec:motivation} motivates the need for preserving the symplectic structure while deriving ROMs when  using quadratic manifold constructions. Section~\ref{subsec:problem_formulation} formalizes the problem considered in this paper. In Section~\ref{sec:symp_quad_jacobian}, we employ quadratic manifolds to develop a novel nonlinear mapping that is guaranteed to be symplectic and then use this mapping to derive a Hamiltonian ROM of~\eqref{eq:HFOM} via the SMG projection. In Section~\ref{sec:symp_quad_pod}, we use data-driven quadratic manifolds to derive approximately Hamiltonian ROMs that are computationally more efficient than the Hamiltonian ROMs derived in Section~\ref{sec:symp_quad_jacobian}. In this alternate approach, we augment the linear symplectic approximation mapping (based on the cotangent lift) with terms of the quadratic Kronecker
product to derive a quadratic approximation mapping and then combine it with a Galerkin projection to derive an approximately Hamiltonian reduced model of~\eqref{eq:HFOM}.

\subsection{A motivating example} \label{sec:motivation}
We motivate the need for symplectic model reduction using quadratic manifolds by demonstrating how the conventional approach for deriving ROMs using quadratic manifolds~\cite{GEELEN2023115717} violates the
underlying Hamiltonian structure. To demonstrate this we consider ROMs for a nonparametric linear wave equation with periodic boundary conditions. The space-time continuous Hamiltonian for the linear wave equation, see definition~\eqref{eq:functional}, is discretized using $n=2048$ grid points. This leads to snapshot vectors with $2n=4096$ entries. The scalar parameter value in this numerical experiment is fixed to $\mu=0.5$. Since we aim to reproduce the dynamics of only a single trajectory, we omit the dependence on the scalar parameter $\mu$ in our derivation of the governing ROM equations in this subsection. The corresponding FOM equations for the linear wave equation are
\begin{equation}
         \dot{\by}(t)=\begin{pmatrix}
             \dot{\bq}(t) \\
             \dot{\bp}(t)  
         \end{pmatrix}=\bA
         \begin{pmatrix}
             \bq(t) \\ 
             \bp(t)
         \end{pmatrix}, \quad \text{with } \bA := \begin{pmatrix}
             \bzero & \In \\
             0.25\mathbf{D}_{\text{fd}} & \bzero
         \end{pmatrix},
\end{equation}
where $\mathbf{D}_{\text{fd}}$ denotes the finite difference approximation of the spatial derivative $\partial_{xx}$. The Hamiltonian FOM is integrated using a symplectic integrator based on the implicit midpoint rule with a fixed time step of $\Delta t=10^{-3}$. For this numerical example, we center the training data about the initial condition \eqref{eq:parameteric_ic}, that is $\by_{\text{ref}}=\by(0)$. Based on $\ns=K=4000$ FOM data snapshots in the range $t \in [0,4]$, we compute a POD basis matrix $\bV$ by stacking the leading $2r$ left-singular vectors of the centered snapshot data matrix. We then obtain $\Vbar$ by numerically solving the optimization problem \eqref{eq:regularized2} via~\eqref{eq:def_Vbar}. Employing the quadratic manifold approximation of the full state $\state(t)$, together with a Galerkin projection step, yields the following ROM equations
     \begin{equation}
         \dot{\stateRed}(t)
=\bV^\top \bA\by_{\text{ref}}+\bV^\top \bA\bV \stateRed(t) + \bV^\top \bA\Vbar \left( \stateRed(t) \otimes \stateRed(t) \right).
    \label{eq:motivation-rom}
    \end{equation}
We then solve \eqref{eq:motivation-rom} using the same fixed time step as the FOM, after which we carry out a reconstruction of the reduced-state data in the original state space using~\eqref{eq:quadratic_approximation}. 

\begin{figure}[tbp]
\small
\captionsetup[subfigure]{oneside,margin={1.8cm,0 cm}}
\begin{subfigure}{.42\textwidth}
       \setlength\fheight{6 cm}
        \setlength\fwidth{\textwidth}
%
%
\definecolor{mycolor1}{rgb}{0.92941,0.69412,0.12549}%
\definecolor{black}{rgb}{0,0,0}%
\begin{tikzpicture}

\begin{axis}[%
width=0.951\fheight,
height=0.596\fheight,
at={(0\fheight,0\fheight)},
scale only axis,
unbounded coords=jump,
xmin=0,
xmax=40,
xlabel style={font=\color{white!15!black}},
xlabel={Reduced dimension $2r$},
ymode=log,
ymin=0.001,
ymax=10,
yminorticks=true,
ylabel style={font=\color{white!15!black}},
ylabel={Relative state error},
axis background/.style={fill=white},
 xmajorgrids,
 ymajorgrids,
]
\addplot [color=black!80, dotted, line width=2.0pt, mark size=4pt, mark=pentagon, mark options={solid, black!80}, forget plot]
  table[row sep=crcr]{%
4	1.2369966222406\\
nan	nan\\
12	0.0722530887782904\\
16	0.0184687046594212\\
20	0.0103096225330478\\
24	0.228696779315956\\
28	0.649673942658989\\
32	0.883231481532633\\
36	3.34651332230691\\
};
\end{axis}
\end{tikzpicture}%
\caption{State error (training data)}
\label{fig: separate_state}
    \end{subfigure}
    \hspace{0.9cm}
    \begin{subfigure}{.42\textwidth}
           \setlength\fheight{6 cm}
           \setlength\fwidth{\textwidth}
\raisebox{-15mm}{\input{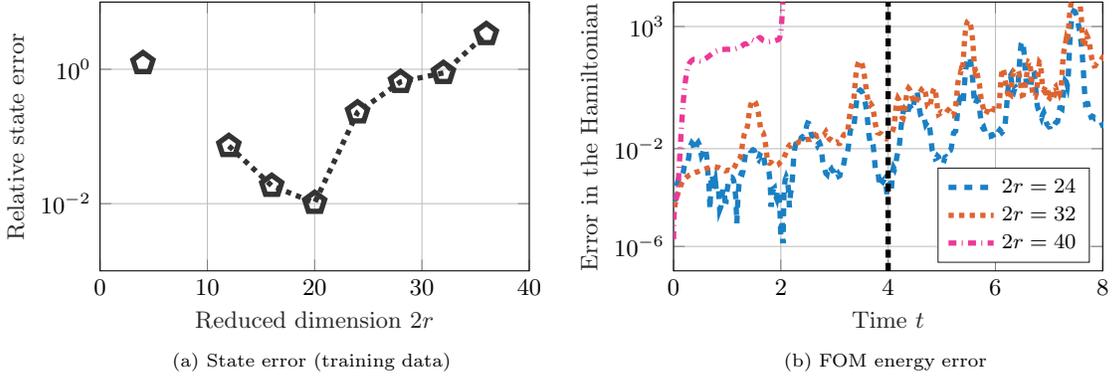}}

\caption{FOM energy error}
\label{fig:separate_Hd}
    \end{subfigure}
   \caption{{\Rb{One-dimensional linear wave equation.}} Even though plot (a) shows low relative state error~\eqref{eq:err_rel_sim_train} for some of the quadratic ROMs in the training data regime, the corresponding {\Rc{energy error}}~\eqref{eq:err_hamiltonian} behavior in plot (b) reveals that the conventional model reduction using quadratic manifolds violates the underlying Hamiltonian structure. The black dashed line indicates end of training time interval.}
 \label{fig:motivation}
\end{figure}
     
Figure~\ref{fig: separate_state} shows the state approximation error in the training data regime for ROMs of different reduced dimensions. We observe that the relative state error (see~\eqref{eq:err_rel_sim_train} {\Rc{for a precise definition}}) in the training data regime decreases from $2r=12$ to $2r=20$ and then increases from $2r=20$ to $2r=36$. For $2r=8$ and $2r=40$, {\Ra{the ROM solutions become unstable}} in the training data regime, which is the reason there are no markers for these points in Figure~\ref{fig: separate_state}. For $2r=32$ and $2r=36$, we observe similar instabilities at approximately $t=8$ (not shown here). Despite some of the quadratic ROMs~\eqref{eq:motivation-rom} exhibiting a relative state error below $10^{-1}$ in the training data regime, they do not conserve the space-discretized Hamiltonian $\dH$. The {\Rc{Hamiltonian error}} plots (see~\eqref{eq:err_hamiltonian} for a precise definition) in Figure~\ref{fig:separate_Hd} show that, as expected, {\Ra{the FOM energy is not preserved by}} standard quadratic ROMs of different sizes. This behavior can be attributed to the reduced quadratic operators not preserving the Hamiltonian structure. In the remainder of this section, we propose two novel approaches that employ data-driven quadratic manifolds to derive structure-preserving ROMs of Hamiltonian systems{\Rc{. The}} numerical experiments in Section~\ref{sec:results} show that the proposed methods lead to accurate and stable ROMs with predictive capabilities.

\subsection{Problem formulation} \label{subsec:problem_formulation}
We consider a canonical high-dimensional parametric Hamiltonian system
\begin{equation}
        \dot{\by}(t;\bmu)=\Jn\nabla_{\by}\dH(\by(t;\bmu);\bmu), \qquad \by(0;\bmu)=\by_0(\bmu).
       \label{eq:FOM_general}
\end{equation}
Specifically, we focus on parametric Hamiltonian systems with FOM Hamiltonians of the form
\begin{equation}\label{eq:quad_sep_Ham}
    \dH(\bq,\bp; \bmu)
= \frac{1}{2} \rT\bp \bH_p(\bmu) \bp
+ \frac{1}{2} \rT\bq \bH_q(\bmu) \bq,
\end{equation}
where $\bH_q(\bmu), \bH_p(\bmu) \in \real^{n \times n}$ are obtained by structure-preserving spatial discretization of parametric Hamiltonian PDEs. The corresponding FOM equations are
\begin{equation}
    \dot{\by}(t; \bmu)
=\begin{pmatrix}
    \dot{\bq}(t; \bmu) \\ \dot{\bp}(t; \bmu)
\end{pmatrix}
=\Jn\nabla_{\by}\dH(\by(t; \bmu); \bmu)
=\begin{pmatrix}
    \bzero & \In \\ -\In & \bzero 
\end{pmatrix}
\begin{pmatrix}
    \bH_q(\bmu) & \bzero\\
    \bzero & \bH_p(\bmu) 
\end{pmatrix}
\begin{pmatrix}
    \bq(t; \bmu) \\
    \bp(t; \bmu)
\end{pmatrix}.
    \label{eq:lin_sep_Ham_sys}
\end{equation}
 {\Ra{The block-diagonal structure in the FOM equations is related to the separable nature of~\eqref{eq:quad_sep_Ham} where the FOM Hamiltonian is additively separable with respect to generalized position vector $\bq$ and  generalized momentum vector $\bp$. Separable Hamiltonian FOMs of the form~\eqref{eq:lin_sep_Ham_sys} appear as models in many science and engineering applications such as the linear wave equation in physics, the linear elasticity equation in solid mechanics, and Maxwell's equations in electromagnetics.}} 

We consider an intrusive model reduction setting in this work where the FOM operators $\bH_q(\bmu)$ and $\bH_p(\bmu)$ are assumed {\Rb{to be}} known.
In the scope of this work, a projection-based model reduction technique consists of (a) an approximation mapping $\bqgamma: \real^{2r} \mapsto \real^{2n}$ and (b) a projection of the residual
\begin{equation}
\br(t; \bmu) :=  \frac{\text{d}}{\text{d}t}\left(\bqgamma(\stateRed(t;\bmu))\right) -  \Jn \nabla_{ \by }\dH(\bqgamma(\stateRed(t;\bmu));\bmu ) \in \real^{2n},
\label{eq:mor_residual}
\end{equation}
where $\stateRed(t; \bmu) \in \real^{2r}$ denotes the reduced state vector. 

In this paper, we propose two novel approaches that employ data-driven quadratic manifolds, as introduced in Section \ref{sec:quadratic}, to derive accurate and stable ROMs for Hamiltonian systems. The two distinct model reduction techniques proposed in Sections \ref{sec:symp_quad_jacobian} and \ref{sec:symp_quad_pod} are named by the approximation map and the projection that are used to define it, SMG-QMCL (with the SMG projection and the QMCL approximation map) and Galerkin-BQ (with a Galerkin projection and the BQ approximation map).
The reduced-order matrix operators in both the SMG-QMCL-ROM and the Galerkin-BQ-ROM may be precomputed in the offline part of the approach. 

Since the parameter vectors $\bmu \in \cP$ are fixed (but arbitrary) for each FOM or ROM, we simplify the notation by omitting the explicit dependence on the parameter: the FOM state vector $\by(t;\bmu)$ at time $t$ and parameter $\bmu$ is therefore often denoted as $\by(t)$. Analogously, the generalized position and momentum vectors $\bq(t;\bmu), \bp(t;\bmu)$ are denoted as $\bq(t)$ and $\bp(t)$, respectively.
\begin{remark}
{\Rb{For noncanonical Hamiltonian systems with a constant and non-degenrate skew-symmetric matrix $\mathbf{S}$ in~\eqref{eq:noncan}, there exists a congruent transformation which can transform the noncanonical Hamiltonian system to the canonical form, see~\cite{peng2016symplectic,gruber2023canonical} for more details \textcolor{black}{about the model reduction of noncanonical Hamiltonian systems.} Thus, the proposed approaches can be applied to noncanonical systems with constant and non-dengerate $\mathbf{S}$ after transforming the noncanonical FOM to a canonical form.}}
\end{remark}
\subsection{Quadratic manifold cotangent lift (QMCL) state approximation and the SMG-QMCL-ROM}
\label{sec:symp_quad_jacobian}

We introduce a model reduction technique that strictly preserves the structure of the Hamiltonian system.
To this end, an approximation mapping is constructed---based on data-driven quadratic manifolds---which is guaranteed to be symplectic. We begin by considering a quadratic approximation mapping for the generalized position vector $\stateqRed \in \real^r$ in the form of
    \begin{align}\label{eq:bqgammaQMCL}
        \bqgammaQMCL(\stateqRed)
:= \qref + \Vq \stateqRed + \Vqbar (\stateqRed \otimes \stateqRed),
    \end{align}
with a constant vector $\qref \in \real^{n}$ and matrices $\Vq \in \real^{n \times r}$, $\Vqbar \in \real^{n \times r(r+1)/2}$ with 
\begin{align}\label{eq:bqgammaQMCL_emb_matrices}
        \rT\Vq \Vq =& \Ir
&\text{and}&
&\rT\Vq \Vqbar = \bzero \in \real^{r \times r(r+1)/2}.
\end{align}
Based thereon, we introduce the \emph{quadratic manifold cotangent lift} (QMCL) as the approximation mapping of the state $\stateRed = (\stateqRed^\top, \statepRed^\top)^\top \in \real^{2r}$ with
    \begin{align}\label{eq:quad_mnf_co_tan}
        \bgammaQMCL(\stateqRed, \statepRed)
:= \begin{pmatrix}
    \bqgammaQMCL(\stateqRed)\\
    \bpgammaQMCL(\stateqRed, \statepRed)
\end{pmatrix}
:= \begin{pmatrix}
    \qref + \Vq \stateqRed + \Vqbar (\stateqRed \otimes \stateqRed)\\
    \pref + \Big( \Vq + \Vqbar \kronFun(\stateqRed) \Big) \vRed(\stateqRed, \statepRed)
\end{pmatrix},
    \end{align}
with $\vRed(\stateqRed, \statepRed)
:= \invb{
   \Ir + \rTb{\kronFun(\stateqRed)} \rT\Vqbar \Vqbar \kronFun(\stateqRed)
}
\statepRed \in \real^r$ where $\kronFun(\stateqRed) := \eval[\stateqRed]{\bD (\cdot \otimes \cdot)} \in \real^{r(r+1)/2 \times r}$ denotes the derivative of the Kronecker product and a constant vector $\pref \in \real^{n}$. We show later in this subsection (in a more general framework) that the approximation mapping $\bgammaQMCL$ is (i) continuously differentiable and (ii) a symplectic map, i.e.\ the Jacobian $\eval[\stateRed]{\bD \bgammaQMCL} \in \real^{2n \times 2r}$ is a symplectic matrix
\begin{align}\label{eq:sympl_map}
        \rT{\eval[\stateRed]{\bD \bgammaQMCL}} \Jn \eval[\stateRed]{\bD \bgammaQMCL} = \Jr
        \qquad \text{for all } \stateRed \in \real^{2r}.
\end{align}
These two properties allow us to formulate the symplectic inverse \eqref{eq:sympl_inv} of the Jacobian $\eval[\stateRed]{\bD \bgammaQMCL}$ which we denote as $\eval[{\stateRed(t)}]{\bD \bgammaQMCL}^+$. Thus, we can directly apply the techniques from \cite{buchfink2021symplecticmanifold} to use the SMG projection.

\begin{definition}[SMG projection~\cite{buchfink2021symplecticmanifold}]
The \emph{SMG projection} requires that the symplectic projection of the residual \eqref{eq:mor_residual} with the Jacobian of the approximation mapping vanishes, that is
\begin{equation}\label{eq:smg_proj_mor_residual}
\eval[{\stateRed(t)}]{\bD \bgammaQMCL}^+ \br(t) = \bzero \in \real^{2r} \quad  {\rm{for\ all\ }} t.
\end{equation}
\end{definition}
\begin{proposition}
    The ROM obtained via the SMG projection~\eqref{eq:smg_proj_mor_residual} of a canonical Hamiltonian system with $\bgammaQMCL$~\eqref{eq:quad_mnf_co_tan} is again a canonical Hamiltonian system.
    \end{proposition}
    \begin{proof}
    We apply the proof from~\cite[Section~2.3]{buchfink2021symplecticmanifold} to the QMCL approximation mapping proposed in the present work: Substituting the residual~\eqref{eq:mor_residual} into the definition of the SMG projection~\eqref{eq:smg_proj_mor_residual} leads to
    \begin{align*}
        \dot{\stateRed}(t)
        =&\eval[{\stateRed(t)}]{\bD \bgammaQMCL}^+ \Jn \nabla_{ \by } \dH(\bgammaQMCL(\stateRed(t)))\\
        =&\Jr \eval[{\stateRed(t)}]{\bD \bgammaQMCL}^\top \nabla_{\by} \dH(\bgammaQMCL(\stateRed(t)))\\
        =&\Jr\nabla_{\stateRed}\dH(\bgammaQMCL(\stateRed(t))),
    \end{align*}
    where we use that for all $\stateRed \in \real^{2r}$ it holds that
    $\eval[{\stateRed}]{\bD \bgammaQMCL}^+
    \eval[{\stateRed}]{\bD \bgammaQMCL}
    = \bI_{2r}$ on the left-hand side,
    $\eval[{\stateRed}]{\bD \bgammaQMCL}^+ \Jn
    = \Jr \eval[{\stateRed}]{\bD \bgammaQMCL}^\top$
    in the second step, and the chain rule for the gradient of the Hamiltonian
    $\nabla_{\stateRed} \dH(\bgammaQMCL(\stateRed))
    = \eval[{\stateRed}]{\bD \bgammaQMCL}^\top \nabla_{\by} \dH(\bgammaQMCL(\stateRed))$ in the last step.
    With the definition of the reduced Hamiltonian $\dHRed := \dH \circ \bgammaQMCL \in \cC^1(\real^{2r})$, we see that this is indeed a canonical Hamiltonian system
    \begin{align}\label{eq:jac_based_nonlin_Ham}
        \dot{\stateRed}(t) = \Jr \nabla_{\stateRed}\dHRed(\stateRed(t)) \in \real^{2r},
    \end{align}
    of reduced dimension $2r$.
\end{proof}
   
\begin{remark}[SMG-QMCL-ROM for the linear Hamiltonian system \eqref{eq:lin_sep_Ham_sys}]
        For the special case of a linear Hamiltonian system of the form~\eqref{eq:lin_sep_Ham_sys}, the ROM obtained using the SMG-QMCL approach can be written as
        \begin{align}\label{eq:jac_based_lin_Ham}
            \dot{\stateRed}(t)
=& \Jr \rT{ \eval[\stateRed(t)]{\bD \bgammaQMCL} } \bH \bgammaQMCL(\stateRed(t)),&
            &\text{with}&
            \bH
=&\begin{pmatrix}
    \bH_q & \bzero\\
    \bzero & \bH_p 
\end{pmatrix}
\in \real^{2n \times 2n}.
        \end{align}
        To write the specific operations out more explicitly, we introduce the abbreviations
        \begin{align}
           \VV
&:=\begin{pmatrix}
   \Vq                 & \bzero\\
   \bzero & \Vq
\end{pmatrix}
\in \real^{2n \times 2r},&
       \VbarVbar
&:= \begin{pmatrix}
   \Vqbar                 & \bzero\\
   \bzero & \Vqbar
\end{pmatrix}
\in \real^{2n \times r(r+1)},&
        \yref
&:= \begin{pmatrix}
    \qref\\
    \pref
\end{pmatrix}
\in \real^{2n},
       \end{align}
       and for $\stateRed = (\stateqRed^\top, \statepRed^\top)^\top \in \real^{2r}$ the respective nonlinear, but low-dimensional, coefficient functions
       \begin{align*}
           \stateV(\stateqRed, \statepRed)
&:= \begin{pmatrix}
   \stateqRed\\
   \vRed(\stateqRed, \statepRed)\\ 
\end{pmatrix}
\in \real^{2r},&
   \stateVbar(\stateqRed, \statepRed)
&:= \begin{pmatrix}
   (\stateqRed \otimes \stateqRed)\\
   \kronFun(\stateqRed)
   \vRed(\stateqRed, \statepRed)
\end{pmatrix}
\in \real^{r(r+1)},
       \end{align*}
       such that the mapping $\bgammaQMCL$ and its derivative are
       \begin{align*}
           \bgammaQMCL(\stateRed)
=& 
\yref
+ \VV \stateV(\stateRed)
+ \VbarVbar \stateVbar(\stateRed)
,&
           \eval[\stateRed]{\bD\bgammaQMCL}
=& \VV \eval[\stateRed]{\bD\stateV}
+ \VbarVbar \eval[\stateRed]{\bD\stateVbar}.
       \end{align*}
       The \emph{SMG-QMCL-ROM} for the linear Hamiltonian system \eqref{eq:quad_sep_Ham} can be written as
       \begin{align}\label{eq:jac_based_lin_Ham_quad_gamma}
           \begin{split}
               \dot{\stateRed}(t)
=&\; \Jr
\rT{\eval[\stateRed(t)]{\bD\stateV}}
\left(
\bHVyref
+
\bHVV\;
\stateV(\stateRed(t))
+
\bHVVbar\;
\stateVbar(\stateRed(t))
\right)
\\
& + \Jr
\rT{\eval[\stateRed(t)]{\bD\stateVbar}}
\left(
\bHVbaryref
+ 
\rT{\bHVVbar}\;
\stateV(\stateRed(t))
+ 
\bHVbarVbar\;
\stateVbar(\stateRed(t))
\right),
           \end{split}
       \end{align}
       with different projections of $\bH$ as
       \begin{gather}
            \bHVyref
:= \rT\VV \bH \yref \in \real^{2r},\quad
\bHVV
:= \rT\VV \bH \VV \in \real^{2r \times 2r},\quad
\bHVVbar
:= \rT\VV \bH \VbarVbar \in\real^{2r \times r(r+1)},\nonumber\\
\label{eq:jac_based_proj_ops}
            \bHVbaryref
:= \rT\VbarVbar \bH \yref \in\real^{r(r+1)},\quad
\bHVbarVbar
:= \rT\VbarVbar \bH \VbarVbar \in\real^{r(r+1) \times r(r+1)}.
       \end{gather}
       Note that all multiplications in \eqref{eq:jac_based_lin_Ham_quad_gamma} are independent of the full dimension $n$ and only depend on the reduced dimension $r$. Thus, this model is offline--online separable if $\bH$ is parameter-separable as then the reduced operators \eqref{eq:jac_based_proj_ops} can be computed offline following standard techniques (see e.g.\ \cite[Section~3]{Haasdonk2011}). In the online stage, the SMG-QMCL approach requires $\bigO(r^5)$ multiplication operations to evaluate the right-hand side of the nonlinear ROM and $\bigO(r^5)$ multiplication operations to evaluate the Jacobian of this right-hand side.
\end{remark}

The derived equations hold for all $\Vq \in \real^{n \times r}$ and $\Vbar \in \real^{n \times r(r+1)/2}$ which fulfill \eqref{eq:bqgammaQMCL_emb_matrices}. In our approach, we compute $\Vq=\bPhi$ through the cotangent lift algorithm from Algorithm~\ref{alg:cotangent_lift} and the matrix $\Vqbar$ for the quadratic part is computed as in \Cref{sec:quadratic} by only using the snapshot matrix $\bQ_{\bmu}$. {\Ra{This approach minimizes the reconstruction errors in $\bq$ based on $\bQ_{\bmu}$ snapshot data, which leads to a linear least-squares problem for computing $\Vqbar$. Alternatively, one could think about construction $\Vqbar$ based on reconstruction errors for $\bq$ and $\bp$. However, since $\bpgammaQMCL$ depends on the Jacobian of $\bqgammaQMCL$, this would lead to challanging nonlinear least-squares problem.}} The SMG-QMCL approach is summarized in Algorithm~\ref{alg:jacobian}.

Since the SMG-QMCL approach is based on the SMG projection,
we can directly transfer results on the stability in the sense of Lyapunov:

\begin{proposition}[\cite{buchfink2021symplecticmanifold}]
    If the Hamiltonian FOM~\eqref{eq:FOM_general} has Lyapunov-stable states and these states are included in the image of the symplectic approximation mapping $\bgammaQMCL$~\eqref{eq:quad_mnf_co_tan}, then these states are also Lyapunov-stable states of the SMG-QMCL-ROM~\eqref{eq:jac_based_lin_Ham_quad_gamma}.
\end{proposition}

It remains to show that the proposed nonlinear approximation map $\bgammaQMCL$ is (i) continuously differentiable and (ii) a symplectic map. The proof is based on a more general approximation mapping, which we refer to as the manifold cotangent lift. We begin by introducing immersions and the Moore--Penrose inverse:
    \begin{definition}[Immersion]
        An \emph{immersion} $\bf \in \cC^1(\real^k, \real^m)$ is a mapping for which the Jacobian $\eval[\bx]{\bD \bf} \in \real^{m \times k}$ is of rank $k$ at every point $\bx \in \real^m$.
    \end{definition}
    \begin{lemma}[Moore--Penrose inverse and immersions]\label{lem:mpinv}
        For an immersion $\bf \in \cC^1(\real^k, \real^m)$ with $k \leq m$,
        the \emph{Moore--Penrose inverse} of the Jacobian
        \begin{equation}\label{eq:mpinv}
            \mpinv{ \eval[\bx]{\bD\bf} }
:= \invb{
    \rT{\eval[\bx]{\bD\bf}}
    \eval[\bx]{\bD\bf}
}
\rT{\eval[\bx]{\bD\bf}}
        \end{equation}
        is well-defined for all $\bx \in \real^k$.
        If $\bf \in \cC^2(\real^k, \real^m)$,
        we have that $\mpinv{ \eval[(\cdot)]{\bD\bf} } \in \cC^1(\real^k, \real^{m \times k})$.
    \end{lemma}
    \begin{proof}
        Since $\bf \in \cC^1(\real^k, \real^m)$ is an immersion with $k \leq m$,
        the Jacobian $\eval[\bx]{\bD \bf} \in \real^{m \times k}$ is of full column rank $k$ for all $\bx \in \real^k$.
        Thus, the product $\rTb{\eval[\bx]{\bD \bf}} \eval[\bx]{\bD \bf}$ is of full rank and thus \eqref{eq:mpinv} is well-defined.
        For an immersion $\bf \in \cC^2(\real^k, \real^m)$,
        we know that $\eval[(\cdot)]{\bD \bf} \in \cC^1(\real^k, \real^{m \times k})$ with constant rank $k$ which is enough to show $\mpinv{ \eval[(\cdot)]{\bD\bf} } \in \cC^1(\real^k, \real^{m \times k})$ (see \cite[Theorem~4.3]{Golub1973}).
    \end{proof}

We propose a nonlinear generalization of the {\Ra{linear symplectic lift employed in}} the cotangent lift {\Ra{algorithm}}, which we refer to as the \emph{manifold cotangent lift}.
\begin{definition}[Manifold cotangent lift]
        For a given immersion $\bqgammaMCL \in \cC^2(\real^r, \real^n)$,
        we define the manifold cotangent lift (MCL) embedding as 
        \begin{align}\label{eq:map_mnf_co_tan}
            \bgammaMCL(\stateqRed, \statepRed)
:=& \begin{pmatrix}
    \bqgammaMCL(\stateqRed)\\
    \bpgammaMCL(\stateqRed, \statepRed)
\end{pmatrix}&
            &\text{with}&
            &\bpgammaMCL(\stateqRed, \statepRed)
:= \pref + \rTb{\mpinv{ \eval[\stateqRed]{\bD\bqgammaMCL} }} \statepRed,
        \end{align}
        for some $\pref \in \real^{r}$.
\end{definition}

We show that an MCL embedding is (i) continuously differentiable and (ii) a symplectic map.
\begin{theorem}
        The MCL embedding $\bgammaMCL$ from \eqref{eq:map_mnf_co_tan} is continuously differentiable.
\end{theorem}
\begin{proof}
        We argue about the differentiability of the 
 two components of $\bgammaMCL$~\eqref{eq:map_mnf_co_tan} in three steps.
 First, we know that the first component $\bqgammaMCL \in \cC^2(\real^{r}, \real^n)$ is continuously differentiable by assumption.
 Second, $\bpgammaMCL$ is smooth in $\statepRed$ since the function is linear in $\statepRed$.
 Finally, since $\bqgammaMCL \in \cC^2(\real^r, \real^n)$ is an immersion, \Cref{lem:mpinv} applies and we know that the Moore--Penrose pseudo inverse is continuously differentiable.
 Then, we know that $\bpgammaMCL$ is continuously differentiable in $\stateqRed$ by composition of continuously differentiable functions.
\end{proof}%
\begin{theorem}%
    The MCL embedding $\bgammaMCL$ from \eqref{eq:map_mnf_co_tan} is a symplectic map \eqref{eq:sympl_map}.
\end{theorem}
\begin{proof}
        For all $\stateRed = (\stateqRed^\top, \statepRed^\top)^\top \in \real^{2r}$, the derivative of the nonlinear mapping $\bgammaMCL$ can be written as
        \begin{align*}
            \eval[\stateRed]{ \bD \bgammaMCL }
= \begin{pmatrix}
    \eval[\stateqRed]{\bD \bqgammaMCL} & \bzero_{n \times r}\\
    \eval[\stateRed]{\Dq \bpgammaMCL}
        & \rTb{\mpinv{ \eval[\stateqRed]{\bD\bqgammaMCL} }}
\end{pmatrix}
\in \real^{2n \times 2r},
        \end{align*}
        where $\bzero_{n \times r} \in \real^{n \times r}$ is the matrix of all zeros
        and $\Dq(\cdot)$ denotes the partial derivative w.r.t.\ $\stateqRed$. 
        Using the above expression for $\eval[\stateRed]{ \bD \bgammaMCL }$, we can write the left-hand side of~\eqref{eq:sympl_map} as
        \begin{align}\label{eq:proof_sympl_intermed_term}
            \rT{\eval[\stateRed]{ \bD \bgammaMCL }} \Jn \eval[\stateRed]{ \bD \bgammaMCL }
= \begin{pmatrix}
    2 \skewPart \left(
        \rT{\eval[\stateqRed]{\bD \bqgammaMCL}}
\eval[\stateRed]{\Dq \bpgammaMCL}
    \right)
    & \rT{\eval[\stateqRed]{\bD \bqgammaMCL}}
\rTb{\mpinv{ \eval[\stateqRed]{\bD\bqgammaMCL} }}\\
    -\mpinv{ \eval[\stateqRed]{\bD\bqgammaMCL} } 
\eval[\stateqRed]{\bD \bqgammaMCL}
    & \bzero_{r \times r}
\end{pmatrix},
        \end{align}
        where $\skewPart(\bA) := \frac{1}{2} (\bA - \rT\bA)$ denotes the skew-symmetric part of $\bA \in \real^{r \times r}$.
        To prove the symplecticity, we need to show that \eqref{eq:proof_sympl_intermed_term} is equal to $\Jr$. We observe that the off-diagonal terms already match due to
        $$\bI_{r}
= \mpinv{ \eval[\stateqRed]{\bD\bqgammaMCL} } \eval[\stateqRed]{\bD \bqgammaMCL}
= \rTb{\mpinv{ \eval[\stateqRed]{\bD\bqgammaMCL} } \eval[\stateqRed]{\bD \bqgammaMCL}}
= \rT{\eval[\stateqRed]{\bD \bqgammaMCL}} \rTb{\mpinv{ \eval[\stateqRed]{\bD\bqgammaMCL} }}.$$
        Thus, it remains to show that the term with the skew-symmetric part in \eqref{eq:proof_sympl_intermed_term} equals zero for which we use that
        \begin{align}\label{eq:proof_sympl_statepRed_term}
            \statepRed
= \bI_{r} \statepRed
= \rT{\eval[\stateqRed]{\bD \bqgammaMCL}}
\rTb{\mpinv{ \eval[\stateqRed]{\bD\bqgammaMCL} }} \statepRed
= \rT{\eval[\stateqRed]{\bD \bqgammaMCL}} \bpgammaMCL(\stateqRed, \statepRed).
        \end{align}
        For a better presentation, we reformulate this equation to index notation and use the notation
$(\bqgammaMCL)_{i,j} = \left( \eval[\stateqRed]{\bD \bqgammaMCL} \right)_{ij}$ for the Jacobian with respect to $\stateqRed$.
        Then, we can write \eqref{eq:proof_sympl_statepRed_term} in index notation for $1\leq i \leq r$ as
        \begin{align*}
            (\statepRed)_i
= \sum_{k=1}^r (\bqgammaMCL)_{k,i}
\left( \bpgammaMCL \right)_k.
        \end{align*}
        Deriving this equation w.r.t.\ $\stateqRed$ yields by using the product rule for $1\leq i,j \leq r$
        \begin{align*}
            0
= \sum_{k=1}^r (\bqgammaMCL)_{k,ij}
\left( \bpgammaMCL \right)_k
+ \underbrace{
    \sum_{k=1}^r
    (\bqgammaMCL)_{k,i} \left( 
        \bpgammaMCL
    \right)_{k,j}
}_{
    = \left(
        \rTb{\eval[\stateqRed]{\bD \bqgammaMCL}}
\eval[\stateRed]{\Dq \bpgammaMCL}
    \right)_{ij}
},
        \end{align*}
        where we identify in the underbrace the term which we want to make an assertion about for \eqref{eq:proof_sympl_intermed_term}.
        We move this term to the left-hand side and obtain for $1 \leq i,j \leq r$
        \begin{align}\label{eq:proof_sympl_sym_term}
            \left(
                \rT{\eval[\stateqRed]{\bD \bqgammaMCL}}
                \eval[\stateRed]{\Dq \bpgammaMCL}
            \right)_{ij}
 = - \sum_{k=1}^r (\bqgammaMCL)_{k,ij} \left( \bpgammaMCL \right)_k.
        \end{align}
        Since $\bqgammaMCL \in \cC^2(\real^{2r}, \real^{2n})$, we know by the Lemma of Schwarz that the Hessian of each component $(\bqgammaMCL)_i$ is symmetric, i.e.\ in the presented index notation $(\bqgammaMCL)_{k,ij} = (\bqgammaMCL)_{k,ji}$.
        Thus we know that the terms in \eqref{eq:proof_sympl_sym_term} are symmetric in $i$ and $j$ and the skew-symmetric part vanishes.
        This concludes the proof that \eqref{eq:proof_sympl_intermed_term} equals $\Jr$ and thus $\bgammaMCL$ is a symplectic map.
\end{proof}

We can show that  the cotangent lift and the QMCL can be interpreted as an MCL embedding. For this purpose, we have to show (a) that these methods can be formulated with \eqref{eq:map_mnf_co_tan} via a specific choice for $\bqgammaMCL$ and (b) that this specific $\bqgammaMCL$ is indeed an immersion.

\begin{proposition}[{\Ra{Cotangent lift algorithm}} as special case of MCL]
    For a matrix $\Vq \in \real^{n \times r}$ with orthonormal columns, $\rT\Vq \Vq = \bI_r$,
    the choice $\bqgammaMCL(\stateqRed) ={\Ra{\qref +}} \Vq \stateqRed$ in the MCL embedding \eqref{eq:map_mnf_co_tan} recovers the cotangent lift with $\Vq = \bPhi$.
\end{proposition}

\begin{proof}
    The Jacobian $\bD\bqgammaMCL \equiv \Vq$ is of full column rank by the assumption of orthonormal columns in $\Vq$ and thus this linear choice for $\bqgammaMCL$ is an immersion.
    Since $\bD\bqgammaMCL \equiv \Vq$ with $\rT\Vq \Vq = \bI_{r}$, we know that $\bpgammaMCL(\stateqRed, \statepRed) ={\Ra{ \pref +}} \invb{\rT\Vq \Vq} \Vq \statepRed ={\Ra{\pref +}}  \Vq \statepRed$.
    Thus, the approximation mapping from \eqref{eq:map_mnf_co_tan} equals the approximation mapping in the {\Ra{linear symplectic lift}} \eqref{eq:linear_approximation} with $\bV$ from \eqref{eq:basis_cotan_lift}.
\end{proof}

\begin{proposition}[QMCL as special case of MCL]
    For the quadratic function $\bqgammaMCL \equiv \bqgammaQMCL$ from \eqref{eq:bqgammaQMCL}, the MCL recovers the QMCL from the beginning of this section.
\end{proposition}

\begin{proof}
    With the requirements \eqref{eq:bqgammaQMCL_emb_matrices} on $\Vq$ and $\Vqbar$, we know
    \begin{align*}
        \rT{\eval[\stateqRed]{\bD\bqgammaMCL}} \eval[\stateqRed]{\bD\bqgammaMCL}
= \Ir + \rTb{\kronFun(\stateqRed)} \rT\Vqbar \Vqbar \kronFun(\stateqRed),
    \end{align*}
    and thus $\bpgammaMCL(\stateqRed, \statepRed) = \bpgammaQMCL(\stateqRed, \statepRed)$.
    Again, using \eqref{eq:bqgammaQMCL_emb_matrices}, we can show that $\bqgammaQMCL$ is an immersion since
        \begin{align*}
            &\eval[\stateqRed]{\bD\bqgammaMCL} = \Vq + \Vqbar \kronFun(\stateqRed)&
            &\text{and thus}&
            &\rT\Vq \eval[\stateqRed]{\bD\bqgammaMCL} = \Ir.
    \end{align*}
    Then, with a rank argument $\eval[\stateqRed]{\bD\bqgammaMCL}$ is of full column rank $r$ for all $\stateqRed \in \real^r$.
\end{proof}

MCL embeddings offer a blueprint for generating (nonlinear) mappings which are guaranteed to be symplectic. This is a clear advantage over \cite{buchfink2021symplecticmanifold} which uses a weakly symplectic approach in the numerical experiments to determine nonlinear approximations of the state. Moreover, the proposed QMCL approximation mapping allows an offline--online separation for parameter-separable, linear Hamiltonian systems (which is another missing key piece in~\cite{buchfink2021symplecticmanifold}) for symplectic model reduction using nonlinear approximation mappings.

\subsection{Blockwise-quadratic (BQ) state approximation and the Galerkin-BQ-ROM} \label{sec:symp_quad_pod}

The SMG-QMCL approach in Section~\ref{sec:symp_quad_jacobian} employs a state-dependent projection~\eqref{eq:quad_mnf_co_tan} which leads to a nonlinear Hamiltonian ROM~\eqref{eq:jac_based_lin_Ham_quad_gamma} with computational complexity $\bigO(r^5)$ in the online stage {\Rc{for linear problems}}. In this section, we propose an alternative approach to derive approximately Hamiltonian ROMs that are computationally more efficient in online computations. This approach also retains the physical interpretation of the state variables at the reduced level by choosing the following quadratic approximation for $\stateRed = (\stateqRed^\top, \statepRed^\top)^\top \in \real^{2r}$
\begin{align}
\label{eq:gammqp}
    \bgammaBQ(\stateRed) :=& \begin{pmatrix}
        \bqgammaBQ(\stateqRed)\\
        \bpgammaBQ(\statepRed)
    \end{pmatrix} := 
    \begin{pmatrix}
       \bq_{\text{ref}} + \bV_q \stateqRed + \Vbar_q ( \stateqRed \otimes \stateqRed) \\
       \bp_{\text{ref}} + \bV_p \statepRed + \Vbar_p ( \statepRed \otimes \statepRed )
    \end{pmatrix},
\end{align}
which we refer to as \emph{blockwise-quadratic} (BQ) approximation mapping.
We then choose $\bV_q=\bV_p=\bPhi$ in accord with cotangent lift algorithm~\ref{alg:cotangent_lift}. The matrices $\Vbar_q$ and $\Vbar_p$ are obtained from centered snapshots $\{ \state_j \}_{j=1}^{\ns}$ through the numerical solution of the following pair of optimization problems:
\begin{equation}
    \argmin_{\overline{\bV}_q} \left( J_q(\bV_q,\overline{\bV}_q,\{\stateqRed_j\}^{n_s}_{j=1}) + \gamma_q \left\| \overline{\bV}_q \right\|_F^2 \right);
    \quad
    \argmin_{\overline{\bV}_p} \left( J_p(\bV_p,\overline{\bV}_p,\{\statepRed_j\}^{n_s}_{j=1}) + \gamma_p \left\| \overline{\bV}_p \right\|_F^2 \right),
      \label{eq:optim_problem}
\end{equation}
where
\begin{equation}
\begin{aligned}
    J_q(\bV_q,\overline{\bV}_q,\{\stateqRed_j\}^{n_s}_{j=1}) &:= \sum_{j=1}^{\ns} \left\| \bq_j - \bPhi\stateqRed_j -\Vbar_q ( \stateqRed_j  \otimes \stateqRed_j ) \right\|_2^2, \\
    J_p(\bV_p,\overline{\bV}_p,\{\statepRed_j\}^{n_s}_{j=1}) &:= \sum_{j=1}^{\ns} \left\| \bp_j - \bPhi\statepRed_j -\Vbar_p ( \statepRed_j  \otimes \statepRed_j ) \right\|_2^2,
\end{aligned}
\end{equation}
with the reduced snapshots $\stateqRed_j = \rT\bPhi \stateq_j$ and $\statepRed_j = \rT\bPhi \statep_j$ from $\state_j = (\stateq_j^\top,\statep_j^\top)^\top \in \real^{2n}$ for $j = 1, \dots, \ns$ and the scalar regularization parameters $\gamma_q,\gamma_p \geq 0$.

{\Rb{Since the BQ state approximation is not a symplectic map, we can not use the SMG projection for deriving ROMs. Instead, we use a Galerkin projection}} that requires the residual \eqref{eq:mor_residual} to be orthogonal to the columns of $\bV=  \begin{pmatrix}
        \bPhi & \bzero \\
        \bzero & \bPhi
    \end{pmatrix}$, that is $\bV^\top \br(t; \bmu) = \bzero$.
    Fusing BQ approximation \eqref{eq:gammqp} with the Galerkin projection, we obtain the \emph{Galerkin-BQ-ROM}
    \begin{equation}\label{pod_based_nonlinear_Ham}
        \dot{\stateRed}(t) = \bV^\top\Jn\nabla_{\by} \dH(\bgammaBQ(\stateRed(t)))=\Jr\bV^\top\nabla_{\by}\dH(\bgammaBQ(\stateRed(t))).
    \end{equation}
For linear Hamiltonian systems of the form~\eqref{eq:quad_sep_Ham}, the Galerkin-BQ-ROM simplifies to
\begin{equation}
\begin{aligned}
    \dot{\stateqRed}(t)
&= \phantom{-}\bPhi^\top \bH_p \bp_{\text{ref}} + \bPhi^\top \bH_p \bPhi \statepRed(t) + \bPhi^\top\bH_p \Vbar_p \left(     \statepRed(t) \otimes \statepRed(t) \right), \\
    \dot{\statepRed}(t)
&=- \bPhi^\top \bH_q \bq_{\text{ref}} - \bPhi^\top \bH_q \bPhi \stateqRed(t) - \bPhi^\top \bH_q \Vbar_q \left( \stateqRed(t) \otimes \stateqRed(t)
\right),
\end{aligned}
\end{equation}
where we used the orthogonality properties of the basis matrices $\bPhi^\top \Vbar_q  = \bPhi^\top \Vbar_p = \bzero$. This system of equations can be rewritten as
\begin{equation}\label{eq:pod_Ham_rom}
    \begin{pmatrix}
        \dot{\stateqRed}(t) \\
        \dot{\statepRed}(t)
    \end{pmatrix}
    = \begin{pmatrix}
        \phantom{-}\tbHrefp\\
        -\tbHrefq
    \end{pmatrix}
    +\begin{pmatrix}
        \bzero & \Ir \\
        -\Ir & \bzero 
    \end{pmatrix}
    \begin{pmatrix}
        \tbHq & \bzero \\
        \bzero & \tbHp
    \end{pmatrix}
    \begin{pmatrix}
        \stateqRed(t) \\
        \statepRed(t)
    \end{pmatrix}\\
    + \begin{pmatrix}
     \bzero & \tbHpbar  \\
     -\tbHqbar & \bzero 
    \end{pmatrix}  \begin{pmatrix}
    \stateqRed(t)  \otimes \stateqRed(t) \\
    \statepRed(t)  \otimes \statepRed(t)
    \end{pmatrix},
\end{equation}
where the reduced operators are
\begin{align} \label{eq:pod_rom_operators}
    \begin{aligned}
        \tbHrefq &:= \bPhi^\top \bH_q \qref \in \real^{r},&
        \tbHq &:=\bPhi^\top \bH_q \bPhi\in \real^{r\times r},&
        \tbHqbar &:=\bPhi^\top \bH_q \Vbar_q\in \real^{r\times r(r+1)/2}, \\
        \tbHrefp &:= \bPhi^\top \bH_p \pref\in \real^{r},&
        \tbHp &:=\bPhi^\top \bH_p \bPhi\in \real^{r\times r},&
        \tbHpbar &:=\bPhi^\top \bH_p \Vbar_p\in \real^{r\times r(r+1)/2}.
    \end{aligned}
\end{align}

Approximation~\eqref{eq:quad_mnf_co_tan} of the SMG-QMCL-ROM transforms the high-dimensional problem~\eqref{eq:lin_sep_Ham_sys} into a reduced nonlinear Hamiltonian system. In sharp contrast, the use of the Galerkin-BQ-ROM leads to a reduced linear Hamiltonian system with a quadratic perturbation term. Nonetheless, \eqref{eq:gammqp} exploits knowledge about the canonical structure of the FOM, leading to an approximately Hamiltonian ROM~\eqref{eq:pod_Ham_rom} that retains the intrinsic coupled structure of the FOM and has interpretable states. The numerical results in Section~\ref{sec:results} demonstrate that this is indeed a very accurate strategy. The Galerkin-BQ-ROM is summarized in Algorithm~\ref{alg:pod}.

The Galerkin-BQ-ROM requires $\bigO(r^3)$ online multiplication operations to evaluate the right-hand side of the ROM and $\bigO(r^4)$ multiplication operations to evaluate the Jacobian of this right-hand side. The SMG-QMCL-ROM, on the other hand, requires $\bigO(r^5)$ multiplications for both of these tasks. In other words, by allowing an approximate Hamiltonian structure we achieve a computationally more efficient online phase compared to the SMG-QMCL-ROM from Section~\ref{sec:symp_quad_jacobian}. Importantly, we do not compromise the improved accuracy from data-driven quadratic manifolds.

\begin{algorithm}[tbp]
\caption{Offline phase for the SMG-QMCL-ROM}
\begin{algorithmic}[1]
\Require Centered snapshot data matrix $\bY_{\bmu}=(\bQ_{\bmu}^\top,\bP_{\bmu}^\top)^\top \in \real^{2n \times \ns}$ arranged as in~\eqref{eq:cotangent_snapshot}, reference state $\by_{\text{ref}}=(\bq_{\text{ref}}^\top,\bp_{\text{ref}}^\top)^\top$, Hamiltonian FOM operators $\bH_q$ and $\bH_p$~\eqref{eq:lin_sep_Ham_sys}, and reduced dimension $r$. 
\Ensure Reduced operators $\bHVyref$, $\bHVbaryref$, $\bHVV$, $\bHVbarV, \bHVbarVbar$ for Hamiltonian ROM~\eqref{eq:jac_based_lin_Ham_quad_gamma} and basis matrices $\Vq$ and $\Vqbar$ \eqref{eq:quad_mnf_co_tan}.
    \State $\Vq = \bPhi$ $\leftarrow$ Compute symplectic basis matrix using Algorithm \ref{alg:cotangent_lift} \Comment{Basis computation}
    \State $\widetilde{\bq}_j=\bPhi^\top\bq_j$ for $j=1,\dots,\ns$ \Comment{Represent ${\Rc{\bq}}_j$ in reduced-dim.\ coordinate system}
    \State $\Vqbar$ $\leftarrow$ Solve linear least-squares problem \eqref{eq:regularized2} \Comment{Representation learning problem} 
    \State $\bHVyref$, $\bHVbaryref$, $\bHVV, \bHVVbar, \bHVbarVbar$ $\leftarrow$ Compute matrix operators from \eqref{eq:jac_based_proj_ops} \Comment{ROM operators}
\end{algorithmic}
\label{alg:jacobian}
\end{algorithm}


 \begin{algorithm}[tbp]
\caption{Offline phase for the Galerkin-BQ-ROM}
\begin{algorithmic}[1]
\Require Centered snapshot data matrix $\bY_{\bmu}=(\bQ_{\bmu}^\top,\bP_{\bmu}^\top)^\top \in \real^{2n \times \ns}$ arranged as in~\eqref{eq:cotangent_snapshot}, reference state $\by_{\text{ref}}=(\bq_{\text{ref}}^\top,\bp_{\text{ref}}^\top)^\top$, Hamiltonian FOM operators $\bH_q$ and $\bH_p$~\eqref{eq:lin_sep_Ham_sys}, and reduced dimension $r$. 
\Ensure Reduced operators $\tbHrefq$, $\tbHrefp$, $\tbHq$, $\tbHp$, $\tbHqbar$, $\tbHpbar$ for Hamiltonian ROM~\eqref{eq:pod_Ham_rom}, and basis matrices $\bV_q$, $\bV_p$, $\Vbar_q$, and $\Vbar_p$~\eqref{eq:gammqp}. 
    \State $\bV_q=\bV_p=\bPhi$ $\leftarrow$ Compute symplectic basis matrix using Algorithm \ref{alg:cotangent_lift} \Comment{Basis computation}
    \State $\widetilde{\bq}_j=\bPhi^\top\bq_j$ for $j=1,\dots,\ns$ \Comment{Represent ${\Rc{\bq}}_j$ in reduced-dim.\ coordinate system}
    \State $\widetilde{\bp}_j=\bPhi^\top\bp_j$ for $j=1,\dots,\ns$ \Comment{Represent ${\Rc{\bp}}_j$ in reduced-dim.\ coordinate system}
    \State $\Vbar_q, \Vbar_p$ $\leftarrow$ Solve linear least-squares problems \eqref{eq:optim_problem} \Comment{Representation learning problem}
    \State $\tbHrefq, \tbHrefp, \tbHq, \tbHp, \tbHqbar, \tbHpbar$ $\leftarrow$ Compute matrix operators from~\eqref{eq:pod_rom_operators} \Comment{ROM operators}
\end{algorithmic}
\label{alg:pod}
\end{algorithm}

\section{Numerical results} \label{sec:results}
In this section, the proposed model reduction methods are applied to {\Rb{two parametrized wave equations, which are prototypical of transport problems}} for which slowly-decaying Kolmogorov $N$-widths have been observed in similar settings in~\cite{greif2019decay}. Section~\ref{sec:err_measures} provides details about the numerical implementations and the reported error measures. {\Rb{In Section~\ref{sec:wave_1d}, we demonstrate that the proposed approaches yield more accurate ROMs than the linear symplectic ROMs for a model that generalizes in the parameter and extrapolates in time. Finally, in Section~\ref{sec:wave_2d}, we demonstrate the parameter extrapolation capabilities of the proposed approaches on a two-dimensional nonlinear wave equation. }}

\subsection{Practical considerations \& error measures}
\label{sec:err_measures}

All numerical experiments in this paper are conducted using MATLAB version 2022a. For time integration we use the implicit midpoint rule for all FOM and ROM simulations. For time-continuous dynamical systems the corresponding time-marching equations are
\begin{equation*}
    \frac{\by_{k+1}-\by_k}{\Delta t}=\bf\left(\frac{\by_k + \by_{k+1}}{2}\right),
\end{equation*}
where $\Delta t$ is the fixed time step. The implicit midpoint rule is a second-order symplectic integrator that exhibits bounded energy error for nonlinear Hamiltonian systems \cite{hairer2006geometric}. The Galerkin-BQ-ROM \eqref{eq:pod_Ham_rom} is solved at every time step using Newton's method. The SMG-QMCL-ROM \eqref{eq:jac_based_lin_Ham_quad_gamma}, on the other hand, uses a quasi-Newton scheme  that neglects the second-order derivatives of $\bgammaQMCL$ for an improvement of the online runtime, following the procedures from~\cite{LEE2020108973}.
The regularization parameters $\gamma_q, \gamma_p$ in~\eqref{eq:optim_problem} are chosen to ensure accurate and stable ROMs throughout the range of time integration. They are found by means of a two-dimensional grid search across a sufficiently wide range of parameter values. {\Rb{For both the linear and the nonlinear wave equation example, we use the same $\gamma_q$ value for the SMG-QMCL-ROMs and the Galerkin-BQ-ROMs which leads to the same projection error in $\bq$ for both approaches. Furthermore, all the numerical examples in this work use the same regularization parameter values for different basis sizes.}}

In the following, we introduce the error measures used in the remainder of the paper with a general approximation $\bqgamma$, which can either be $\bgammaLSL$~\eqref{eq:linear_approximation}, $\bgammaQMCL$~\eqref{eq:quad_mnf_co_tan}, or $\bgammaBQ$~\eqref{eq:gammqp}. With $\bqgamma_q$, we denote the restriction of the mapping to the variable $\bq$. The average \emph{relative projection error} in $\bq(t;\bmu)$ is computed as
\begin{equation}\label{eq:err_rel_proj_train}
  {\Rb{\text{err}_{\text{proj,q}} =}}  \frac{1}{M}\sum_{i=1}^M \frac{\lVert \bQ(\bmu_i)- {\color{black}\bqgamma_q} (\bV^\top \bQ(\bmu_i)) \rVert^2_F}{\lVert \bQ(\bmu_i) \rVert^2_F},
\end{equation}
where $\bmu_1, \dots, \bmu_M \in \mathcal{P}$ denotes a set of parameters. This error is valid for both the training and the test parameters. 
An equivalent metric ${\Rb{\text{err}_{\text{proj,p}}}}$ is defined for $\bp(t;\bmu)$. The average \emph{relative state error} is computed as 
\begin{equation}\label{eq:err_rel_sim_train}
   {\Rb{\text{err}_{\text{sim}}=}}  \frac{1}{M}\sum_{j=1}^M \frac{\lVert \bY(\bmu_j)- {\color{black}\bqgamma} (\widetilde{\bY}(\bmu_j)) \rVert^2_F}{\lVert \bY(\bmu_j)\rVert^2_F},
\end{equation}
where $\widetilde{\bY}(\bmu_{j})$ is obtained from the ROM simulations and $\bqgamma(\widetilde{\bY}(\bmu_j))$ is the reconstruction of a trajectory at parameter value $\bmu_j$ in the original state space. {\Rb{We introduce two separate energy error measures for comparing the preservation of the Hamiltonian in the different ROMs. Since the FOM state trajectories for the linear wave equation preserve the energy exactly, the \emph{error in the Hamiltonian} for the linear wave equation example in Section~\ref{sec:wave_1d} is computed as
\begin{equation}\label{eq:err_hamiltonian}
  {\Rb{\Delta H_{\text{lin}}(t)=}}  \left| \dH(\bqgamma(\stateRed(t; \bmu)); \bmu)-\dH(\bqgamma(\stateRed(0; \bmu)); \bmu) \right|,
\end{equation}
where $\dH(\bqgamma(\stateRed(t; \bmu)); \bmu)$ is the FOM energy approximation. For nonlinear Hamiltonian FOMs, the FOM state trajectories do not preserve the energy exactly because the the implicit midpoint scheme only preserves invariances up to the quadratic order. As a result, the \emph{error in the Hamiltonian} for the nonlinear wave equation example in Section~\ref{sec:wave_2d} is computed as  
\begin{equation}\label{eq:err_hamiltonian_non}
  {\Rb{\Delta H_{\text{nonlin}}(t)=}}  \left| \dH(\bqgamma(\stateRed(t; \bmu)); \bmu)-\dH(\state(t; \bmu); \bmu) \right|.
\end{equation}
}}


\subsection{Parametrized linear wave equation} \label{sec:wave_1d}
We revisit the linear wave example from Section~\ref{sec:motivation}. This model problem is similar to the parametric linear wave example in~\cite{afkham2017structure}. Let $\Omega=\left(-0.5,0.5\right) \subset \real $ be the spatial domain and consider the parametrized one-dimensional wave equation 
     \begin{equation}
    \frac{\partial^2 }{\partial t^2}\varphi(x,t;\mu)=\mu^2 \frac{\partial^2 }{\partial x^2}\varphi(x,t;\mu), 
    \label{eq:wave1d_param}
\end{equation}
with the state $\varphi(x,t;\mu)$ at the spatial coordinate $x \in \Omega$, time $t \in  (0,T]$, and the scalar parameter $\mu \in \mathcal {P}=[5/12,5/7]$. Homogenous Dirichlet boundary conditions
     \begin{equation}
        \varphi(-0.5,t;\mu)=\varphi(0.5,t;\mu)=0, 
     \end{equation}
are imposed for $t \in (0,T]$ and $\mu \in \mathcal{P}$. We consider a parametric initial condition of the form $\varphi_0(x;\mu):=h(s(x;\mu))$ based on the spline function 
     \begin{equation}\label{eq:parameteric_ic}
     h(s(x;\mu)):= \begin{cases} 
      1-\frac{3}{2}s(x;\mu)^2 + \frac{3}{4}s(x;\mu)^3& 0\leq s\leq 1, \\
      \frac{1}{4}(2-s(x;\mu))^3 & 1\leq s\leq 2, \\
      0 & s>2
   \end{cases}
\end{equation}
with $s(x;\mu):=\left(4\left|x + \frac{1}{2}-\frac{\mu}{2} \right|\right)/\mu$ for which the exact solution is given by $\varphi(x,t;\mu)=\varphi_0(x-\mu t;\mu)$. 
\subsubsection{Hamiltonian PDE formulation and FOM implementation details}
We rewrite \eqref{eq:wave1d_param} as an infinite-dimensional Hamiltonian system with $q(x,t;\mu):=\varphi(x,t;\mu)$ and $p(x,t;\mu):=\partial \varphi(x,t;\mu)/\partial t$. The associated Hamiltonian functional is
    \begin{equation}
        \mathcal{H}(q(x,t;\mu), p(x,t;\mu); \mu)
=\int_{\Omega} 
\left[
    \frac{1}{2}p(x,t;\mu)^2
    + \frac{1}{2}\mu^2 \left(
        \frac{\partial}{\partial x}q(x,t;\mu)
    \right)^2 
\right] \ \text{d}x,
\label{eq:functional}
    \end{equation}
     and the original PDE can be recast as a Hamiltonian PDE 
    \begin{equation}
   \frac{\partial }{\partial t}q(x,t;\mu)
=\frac{\delta \cH}{\delta p}{\Rc{(q,p;\mu)}}
=p(x,t;\mu), \qquad 
    \frac{\partial }{\partial t}p(x,t;\mu)
=-\frac{\delta \cH}{\delta q}{\Rc{(q,p;\mu)}}
=\mu^2\frac{\partial^2 }{\partial x^2}q(x,t;\mu).
    \end{equation}
We discretize the spatial domain $\Omega$ with $n=2048$ equally spaced grid points leading to {\Ra{a Hamiltonian FOM of dimension $2n=4096$.}} {\Rc{Using a finite difference scheme, we obtain the following space-discretized Hamiltonian}} 
    \begin{equation}
    \dH( \bq(t;\mu), \bp(t;\mu))
=\Delta x\sum^n_{i=1} \left[
    \frac{p_i(t;\mu)^2}{2}
    + \frac{\mu^2 (q_{i+1}(t;\mu)-q_i(t;\mu))^2}{4\Delta x^2} + \frac{\mu^2 (q_{i}(t;\mu)-q_{i-1}(t;\mu))^2}{4\Delta x^2} 
\right] ,
\end{equation}
{\Rc{with}}
\begin{align*}
    q_i(t;\mu)&:=\varphi(x_i,t;\mu),&
\bq(t;\mu)&=(q_1(t;\mu), \ldots, q_n(t;\mu) )^\top \in \real^n,\\
    p_i(t;\mu)&:=\frac{\partial}{\partial t} \varphi(x_i,t:\mu),&
\bp(t;\mu)&=(p_1(t;\mu), \ldots, p_n(t;\mu) )^\top \in \real^n.
\end{align*}
The corresponding parametrized Hamiltonian FOM equals
\begin{equation}
    \dot{\by}(t;\mu) =\begin{pmatrix}
       \dot{\bq}(t;\mu)  \\
       \dot{\bp}(t;\mu) 
     \end{pmatrix}=\Jn \nabla_{ \by }\dH(\by(t;\mu))
     =\begin{pmatrix}
       \bzero & \In \\
       \mu^2\mathbf{D}_{\text{fd}} & \bzero
     \end{pmatrix}
     \begin{pmatrix}
       \bq(t;\mu) \\
       \bp(t;\mu)
     \end{pmatrix}.
     \label{eq:wave_FOM}
\end{equation}
\subsubsection{{\Rb{Generalization to unseen parameters}}}\label{sec:results_lw}

Let $\mu_1=0.417$, $\mu_2=0.476$, $\mu_3=0.536$, $\mu_4=0.595$, $\mu_5=0.655$, $\mu_6=0.714 \in \mathcal P$ be $M=6$ parameters equidistantly distributed in $\mathcal{P}$. In this study, a training dataset is built by integrating the Hamiltonian FOM \eqref{eq:wave_FOM} for each training parameter with the implicit midpoint method until final time ${\Rc{T}}=1$. We use a fixed time step of $\Delta t=2.5 \times 10^{-4}$. For this study, we do not shift the trajectory snapshot data, i.e., $\bq_{\text{ref}}=\bp_{\text{ref}}=\bzero$. From these six trajectories, we construct the nonlinear approximation functions $\bgammaQMCL$~\eqref{eq:quad_mnf_co_tan} and the corresponding SMG-QMCL-ROM, as well as $\bgammaBQ$~\eqref{eq:gammqp} and the corresponding Galerkin-BQ-ROM. {\Rb{For this study, we found $\gamma_q=\gamma_p=10^2$ to be a robust choice for SMG-QMCL and Galerkin-BQ ROMs of dimension $4\leq2r\leq 20$. We provide more details on this selection in Section~\ref{sec:gamma}.}} We consider $M_{\text{test}}=2$ test parameters $\mu_{\text{test},1}=0.51 \in \mathcal{P}$ and $\mu_{\text{test},2}=0.625\in \mathcal{P}$ to evaluate how these data-driven ROMs generalize for parameter values outside the training dataset.



\begin{figure}[tbp]
\small
\captionsetup[subfigure]{oneside,margin={1.8cm,0 cm}}
\begin{subfigure}{.45\textwidth}
       \setlength\fheight{6 cm}
        \setlength\fwidth{\textwidth}
%
%
\begin{tikzpicture}
\small

\begin{axis}[%
width=0.976\fheight,
height=0.59\fheight,
at={(0\fheight,0\fheight)},
scale only axis,
xmin=2,
xmax=20,
xlabel style={font=\color{white!15!black}},
xlabel={Reduced dimension $2r$},
ymode=log,
ymin=0.001,
ymax=1,
yminorticks=true,
ylabel style={font=\color{white!15!black}},
ylabel={Relative projection error},
axis background/.style={fill=white},
xmajorgrids,
ymajorgrids,
legend style={at={(0.65,1.1)}, anchor=south west, legend cell align=left, align=left, draw=white!15!black},
legend style={font=\footnotesize}
]
\addplot [color=blue!80, dotted, line width=2.0pt, mark size=4.0pt, mark=o, mark options={solid, blue!80}]
  table[row sep=crcr]{%
4	0.781335443314784\\
6	0.702325190426793\\
8	0.619901012251458\\
10	0.601444775447315\\
12	0.492038252838663\\
14	0.0462115626426013\\
16	0.0389621419279202\\
18	0.0061035094804675\\
20	0.00426149219795537\\
};
\addlegendentry{LSL-ROM}

\addplot [color=red!80, dashed, line width=2.0pt, mark size=4pt, mark=square, mark options={solid, red!80}]
  table[row sep=crcr]{%
4	0.401007931083306\\
6	0.259336993336295\\
8	0.0742634922009393\\
10	0.0260448893615363\\
12	0.0111208229868014\\
14	0.0062312032461161\\
16	0.00345767554055899\\
18	0.00244266582013545\\
20	0.00206962541806258\\
};
\addlegendentry{SMG-QMCL-ROM}

\addplot [color=green!80, dashdotted, line width=2.0pt, mark size=4.0pt, mark=triangle, mark options={solid, green!80}]
  table[row sep=crcr]{%
4	0.401007931083306\\
6	0.259336993336295\\
8	0.0742634922009393\\
10	0.0260448893615363\\
12	0.0111208229868014\\
14	0.0062312032461161\\
16	0.00345767554055899\\
18	0.00244266582013545\\
20	0.00206962541806258\\
};
\addlegendentry{Galerkin-BQ-ROM}

\end{axis}

\begin{axis}[%
width=1.259\fheight,
height=0.743\fheight,
at={(-0.164\fheight,-0.097\fheight)},
scale only axis,
xmin=0,
xmax=1,
ymin=0,
ymax=1,
axis line style={draw=none},
ticks=none,
axis x line*=bottom,
axis y line*=left
]
\end{axis}
\end{tikzpicture}%
\caption{Relative projection error in $\bq$}
\label{fig:LW_proj_q}
    \end{subfigure}
    \hspace{0.9cm}
    \begin{subfigure}{.45\textwidth}
           \setlength\fheight{6 cm}
           \setlength\fwidth{\textwidth}
\raisebox{-60mm}{
%
%
\begin{tikzpicture}

\begin{axis}[%
width=0.976\fheight,
height=0.59\fheight,
at={(0\fheight,0\fheight)},
scale only axis,
xmin=2,
xmax=20,
xlabel style={font=\color{white!15!black}},
xlabel={Reduced dimension $2r$},
ymode=log,
ymin=0.01,
ymax=1,
yminorticks=true,
ylabel style={font=\color{white!15!black}},
ylabel={Relative projection error},
axis background/.style={fill=white},
xmajorgrids,
ymajorgrids,
legend style={at={(0.486,0.727)}, anchor=south west, legend cell align=left, align=left, draw=white!15!black}
]
\addplot [color=blue!80, dotted, line width=2.0pt, mark size=4.0pt, mark=o, mark options={solid, blue!80}]
  table[row sep=crcr]{%
4	0.563258970495644\\
6	0.471692452854397\\
8	0.269281656654424\\
10	0.146212780514712\\
12	0.0999402702027812\\
14	0.0831703425700218\\
16	0.0480138789621741\\
18	0.0315730186207208\\
20	0.0234918802477253\\
};

\addplot [color=red!80, dashed, line width=2.0pt, mark size=4pt, mark=square, mark options={solid, red!80}]
  table[row sep=crcr]{%
4	0.585524628337326\\
6	0.416928083477077\\
8	0.186812135528262\\
10	0.0837559442907907\\
12	0.0545122196975472\\
14	0.0307944483323163\\
16	0.0216304325819216\\
18	0.018506438480999\\
20	0.0168199434378351\\
};

\addplot [color=green!80, dashdotted, line width=2.0pt, mark size=4.0pt, mark=triangle, mark options={solid, green!80}]
  table[row sep=crcr]{%
4	0.374268491097728\\
6	0.221599066097063\\
8	0.145423882092875\\
10	0.0802882205907461\\
12	0.0432235315385919\\
14	0.0285122695375385\\
16	0.0187673527001748\\
18	0.0154405541392028\\
20	0.0117576574134853\\
};

\end{axis}

\begin{axis}[%
width=1.259\fheight,
height=0.743\fheight,
at={(-0.164\fheight,-0.097\fheight)},
scale only axis,
xmin=0,
xmax=1,
ymin=0,
ymax=1,
axis line style={draw=none},
ticks=none,
axis x line*=bottom,
axis y line*=left
]
\end{axis}
\end{tikzpicture}%
}
\caption{Relative projection error in $\bp$}
\label{fig:LW_proj_p}
    \end{subfigure}
\caption{ {\Rb{One-dimensional linear wave equation (generalization to unseen parameters).}} Data-driven approximations based on quadratic manifolds yield lower relative projection error~\eqref{eq:err_rel_proj_train} than the linear symplectic subspaces for $\bq$ and $\bp$ variables. The regularization factors are chosen to be $\gamma_q=\gamma_p=10^{2}$. }
 \label{fig:LW_proj}
\end{figure}
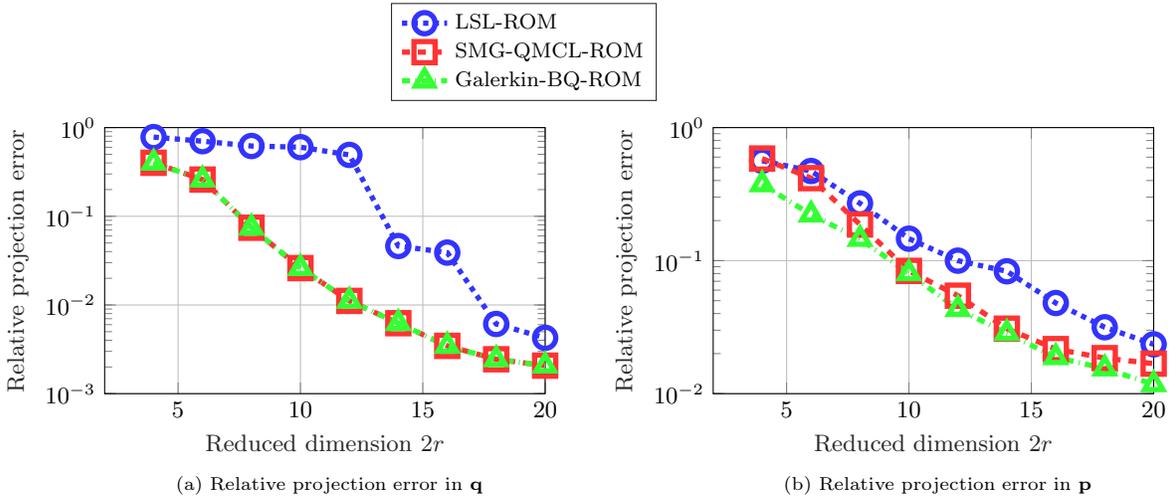

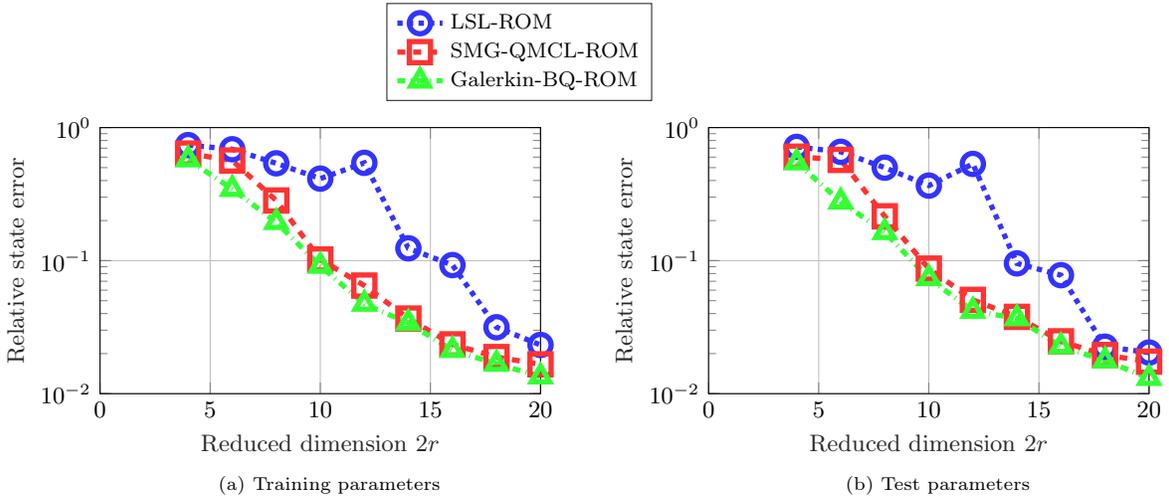
\begin{figure}[tbp]
\small
\captionsetup[subfigure]{oneside,margin={1.8cm,0 cm}}
\begin{subfigure}{.45\textwidth}
       \setlength\fheight{6 cm}
        \setlength\fwidth{\textwidth}
%
%
\begin{tikzpicture}
\small 
\begin{axis}[%
width=0.976\fheight,
height=0.59\fheight,
at={(0\fheight,0\fheight)},
scale only axis,
xmin=0,
xmax=20,
xlabel style={font=\color{white!15!black}},
xlabel={Reduced dimension $2r$},
ymode=log,
ymin=0.01,
ymax=1,
yminorticks=true,
ylabel style={font=\color{white!15!black}},
ylabel={Relative state error},
axis background/.style={fill=white},
xmajorgrids,
ymajorgrids,
legend style={at={(0.65,1.1)}, anchor=south west, legend cell align=left, align=left, draw=white!15!black},
legend style={font=\footnotesize}
]
\addplot [color=blue!80, dotted, line width=2.0pt, mark size=4.0pt, mark=o, mark options={solid, blue!80}]
  table[row sep=crcr]{%
4	0.738699338441775\\
6	0.687211308291726\\
8	0.538956455100131\\
10	0.413466329461469\\
12	0.543273881629925\\
14	0.123720005473688\\
16	0.0925823737485729\\
18	0.0315724881983109\\
20	0.0232319437182599\\
};
\addlegendentry{LSL-ROM}

\addplot [color=red!80, dashed, line width=2.0pt, mark size=4pt, mark=square, mark options={solid, red!80}]
  table[row sep=crcr]{%
4	0.63819686365414\\
6	0.565833762139555\\
8	0.283089333314099\\
10	0.101926836448981\\
12	0.0649537575094531\\
14	0.0367048156950361\\
16	0.0234751209292224\\
18	0.0188447963208537\\
20	0.01661225307388\\
};
\addlegendentry{SMG-QMCL-ROM}

\addplot [color=green!80, dashdotted, line width=2.0pt, mark size=4.0pt, mark=triangle, mark options={solid, green!80}]
  table[row sep=crcr]{%
4	0.575584704420031\\
6	0.344986885417358\\
8	0.194680783098774\\
10	0.0912799290285405\\
12	0.0472966983693054\\
14	0.0339337985157404\\
16	0.0212504079001584\\
18	0.0167534429042319\\
20	0.0134499381257655\\
};
\addlegendentry{Galerkin-BQ-ROM}

\end{axis}
\end{tikzpicture}%
\caption{Training parameters}
\label{fig:LW_state_train}
    \end{subfigure}
    \hspace{0.9cm}
    \begin{subfigure}{.45\textwidth}
           \setlength\fheight{6 cm}
           \setlength\fwidth{\textwidth}
\raisebox{-59mm}{
%
%
\begin{tikzpicture}

\begin{axis}[%
width=0.976\fheight,
height=0.59\fheight,
at={(0\fheight,0\fheight)},
scale only axis,
xmin=0,
xmax=20,
xlabel style={font=\color{white!15!black}},
xlabel={Reduced dimension $2r$},
ymode=log,
ymin=0.01,
ymax=1,
yminorticks=true,
ylabel style={font=\color{white!15!black}},
ylabel={Relative state error},
axis background/.style={fill=white},
xmajorgrids,
ymajorgrids,
legend style={at={(0.678,0.701)}, anchor=south west, legend cell align=left, align=left, draw=white!15!black}
]
\addplot [color=blue!80, dotted, line width=2.0pt, mark size=4.0pt, mark=o, mark options={solid, blue!80}]
  table[row sep=crcr]{%
4	0.715978822217625\\
6	0.657632022650425\\
8	0.49979771199157\\
10	0.366424249081413\\
12	0.533265140609842\\
14	0.0953500901264656\\
16	0.077736958463625\\
18	0.0225511671801405\\
20	0.020376983524802\\
};

\addplot [color=red!80, dashed, line width=2.0pt, mark size=4pt, mark=square, mark options={solid, red!80}]
  table[row sep=crcr]{%
4	0.604612244528034\\
6	0.571794871168971\\
8	0.215303735823801\\
10	0.0871174172293378\\
12	0.0506879124059949\\
14	0.0377946285902017\\
16	0.0247467888162811\\
18	0.0195282632420229\\
20	0.01740166864186\\
};

\addplot [color=green!80, dashdotted, line width=2.0pt, mark size=4.0pt, mark=triangle, mark options={solid, green!80}]
  table[row sep=crcr]{%
4	0.546692847339913\\
6	0.277398552177696\\
8	0.163517634787314\\
10	0.0736689207109806\\
12	0.0418683551712522\\
14	0.036818333667255\\
16	0.0226624771470306\\
18	0.0177358360493477\\
20	0.0130947835616436\\
};


\end{axis}

\begin{axis}[%
width=1.273\fheight,
height=0.705\fheight,
at={(-0.165\fheight,-0.098\fheight)},
scale only axis,
xmin=0,
xmax=1,
ymin=0,
ymax=1,
axis line style={draw=none},
ticks=none,
axis x line*=bottom,
axis y line*=left
]
\end{axis}
\end{tikzpicture}%
}
\caption{Test parameters}
\label{fig:LW_state_test}
    \end{subfigure}
\caption{{\Rb{One-dimensional linear wave equation (generalization to unseen parameters).}} Galerkin-BQ-ROMs and SMG-QMCL-ROMs achieve lower state error~\eqref{eq:err_rel_sim_train} than the linear symplectic ROMs for both training and test parameters.}
 \label{fig:LW_state}
\end{figure}

In Figure~\ref{fig:LW_proj}, we compare the relative projection error from \eqref{eq:err_rel_proj_train} for the training data for different values of the reduced dimension. For both $\bq$ and $\bp$ the quadratic manifold approximations yield higher accuracy compared to the linear symplectic subspaces for all reduced dimension values. {\Rb{For this study}}, we use the same quadratic approximations for $\bq$ in both the SMG-QMCL and the Galerkin-BQ approach, and therefore the relative projection error for $\bq$ in Figure~\ref{fig:LW_proj_q} is the same. The projection error comparison for $\bp$ is shown in Figure~\ref{fig:LW_proj_p} where we observe that the Galerkin-BQ yields marginally lower projection error than the SMG-QMCL.

The comparison of the relative state error \eqref{eq:err_rel_sim_train} between nonlinear SMG-QMCL-ROM and quadratic Galerkin-BQ-ROM is shown in Figure~\ref{fig:LW_state_train} and Figure~\ref{fig:LW_state_test} for the training and the test parameters, respectively. Compared to the linear symplectic ROMs, we observe that SMG-QMCL-ROM and Galerkin-BQ-ROM obtain lower relative state error in both training and testing data regimes. The plots in Figure~\ref{fig:LW_state_train} and Figure~\ref{fig:LW_state_test} show that the linear symplectic subspace approach produces the least accurate ROMs whereas the Galerkin-BQ-ROM admits the highest accuracy in both training and test data regimes. 

\begin{figure}[tbp]
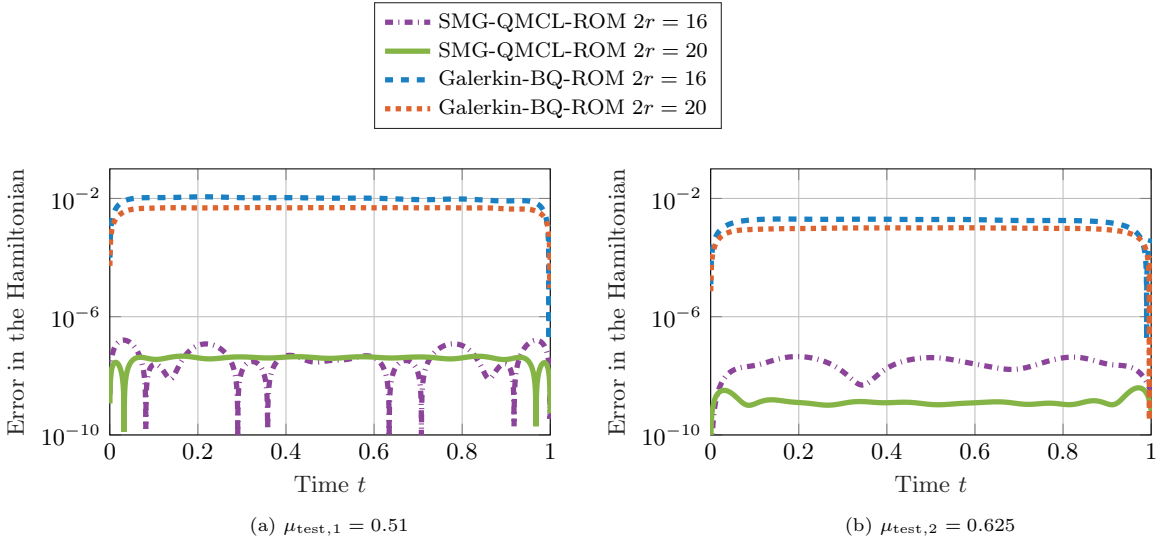

\small
\captionsetup[subfigure]{oneside,margin={1.8cm,0 cm}}
\begin{subfigure}{.45\textwidth}
       \setlength\fheight{6 cm}
        \setlength\fwidth{\textwidth}
\input{figures/LW2/LW_FOM_mu1.tex}
\caption{$\mu_{\text{test},1}=0.51$}
\label{fig:LW_FOM_mu1}
    \end{subfigure}
    \hspace{0.8cm}
    \begin{subfigure}{.45\textwidth}
           \setlength\fheight{6 cm}
           \setlength\fwidth{\textwidth}
\raisebox{-66mm}{\input{figures/LW2/LW_FOM_mu2.tex}}
\caption{$\mu_{\text{test},2}=0.625$}
\label{fig:LW_FOM_mu2}
    \end{subfigure}
\caption{{\Rb{One-dimensional linear wave equation (generalization to unseen parameters).}} {\Rc{The Hamiltonian error~\eqref{eq:err_hamiltonian} for the Galerkin-BQ-ROMs remains below}} $10^{-2}$ for $\mu_{\text{test},1}=0.51$ and $\mu_{\text{test},2}=0.625$ whereas {\Rc{the Hamiltonian error for the SMG-QMCL-ROMs remains below }}$10^{-6}$ for both test parameter values.}
 \label{fig:LW_energy}
\end{figure}

{\Rc{The energy error}}~\eqref{eq:err_hamiltonian} plots in Figures~\ref{fig:LW_FOM_mu1} and \ref{fig:LW_FOM_mu2} show that the approximately Hamiltonian Galerkin-BQ-ROMs exhibit bounded {\Rc{error in the Hamiltonian}} for both $\mu_{\text{test},1}=0.51$ and $\mu_{\text{test},2}=0.625$. The Hamiltonian SMG-QMCL-ROMs demonstrate a substantially lower {\Rc{Hamiltonian error}} due to the usage of the exactly symplectic mapping $\bgammaQMCL$ and the structure-preserving SMG projection. For both approaches, we observe a decrease in the {\Rc{Hamiltonian error}} when we increase the ROM size from $2r=16$ to $2r=20$.

\begin{figure}[tbp]
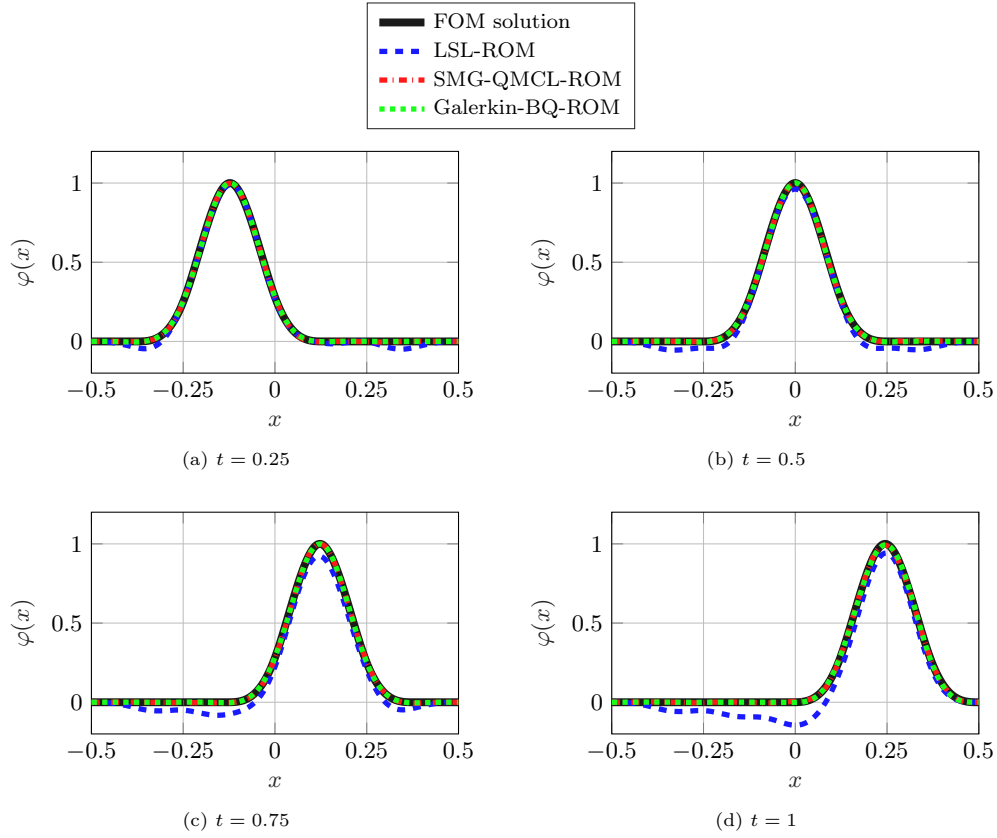

\centering
\small
\begin{subfigure}{.4\textwidth}
       \setlength\fheight{5cm}
        \setlength\fwidth{.8\textwidth}
\input{figures/LW2/LW_T_25_mu1.tex}
\subcaption{$t=0.25$}
 \end{subfigure}
\hspace{0.5cm}
\begin{subfigure}{.4\textwidth}
           \setlength\fheight{5cm}
           \setlength\fwidth{.8\textwidth}
\raisebox{-6.25cm}{\input{figures/LW2/LW_T_50_mu1.tex}}
\subcaption{$t=0.5$}
\end{subfigure}\\[0.5cm]

%

%
\begin{subfigure}{.4\textwidth}
           \setlength\fheight{5cm}
           \setlength\fwidth{\textwidth}
\input{figures/LW2/LW_T_75_mu1.tex}
\subcaption{$t=0.75$}
\end{subfigure}
\hspace{0.5cm}
 \begin{subfigure}{.4\textwidth}
           \setlength\fheight{5cm}
           \setlength\fwidth{.8\textwidth}
\input{figures/LW2/LW_T_1_mu1.tex}
\subcaption{$t=1$}
 \end{subfigure}
    \caption{{\Rb{One-dimensional linear wave equation (generalization to unseen parameters)}}. Plots show the numerical approximation of the solution of \eqref{eq:wave1d_param} for $\mu_{\text{test}}=0.51$ using low-dimensional ($2r=16$) LSL-ROM, Galerkin-BQ-ROM and SMG-QMCL-ROM at different $t$ values. The Galerkin-BQ-ROMs and SMG-QMCL-ROMs capture the correct wave shape at $t=1$ whereas the LSL-ROM $2r=16$ exhibits spurious oscillations from $t=0.25$ onwards.}
\label{fig:LW_traj}
\end{figure}

We compare the approximate numerical solution for $\mu_{\text{test},1}=0.51$ of the linear symplectic ROM, the SMG-QMCL-ROM and the Galerkin-BQ-ROM of size $2r=16$ in Figure~\ref{fig:LW_traj}. Even though the FOM solution snapshots for $\mu_{\text{test},1}=0.51$ are not included in the training data, the ROMs based on an quadratic approximation, SMG-QMCL-ROM and Galerkin-BQ-ROM, capture the correct wave shape at $t=1$. The approximate solutions of the linear symplectic ROM, on the other hand, suffer from spurious oscillations, which grow in magnitude as the solution evolves with time.
\subsubsection{{\Rb{Sensitivity to regularization prameters}}}
\label{sec:gamma}
\begin{figure}[tbp]
\small
\begin{subfigure}{.26\textwidth}
       \setlength\fheight{5 cm}
        \setlength\fwidth{\textwidth}
%
%
\begin{tikzpicture}

\begin{axis}[%
width=0.694\fheight,
height=0.59\fheight,
at={(0\fheight,0\fheight)},
scale only axis,
xmode=log,
xmin=0.001,
xmax=1000,
xminorticks=true,
xtick={1e-3,1e-1,1e1,1e3},
xlabel style={font=\color{white!15!black}},
xlabel={Regularization parameter $\gamma_q$},
ymode=log,
ymin=0.001,
ymax=0.1,
yminorticks=true,
ylabel style={font=\color{white!15!black}},
ylabel={Relative projection error},
axis background/.style={fill=white},
xmajorgrids,
ymajorgrids,
legend style={at={(1.403,1.2)}, anchor=south west, legend cell align=left, align=left, draw=white!15!black}
]
\addplot [color=blue!80, dotted, line width=2.0pt, mark size=4.0pt, mark=o, mark options={solid, blue!80}]
  table[row sep=crcr]{%
0.001	0.0389621419279202\\
0.01	0.0389621419279202\\
0.1	0.0389621419279202\\
1	0.0389621419279202\\
10	0.0389621419279202\\
100	0.0389621419279202\\
1000	0.0389621419279202\\
};
\addlegendentry{LSL-ROM}

\addplot [color=red!80, dashed, line width=2.0pt, mark size=4pt, mark=square, mark options={solid, red!80}]
  table[row sep=crcr]{%
0.001	0.00299588260127251\\
0.01	0.00299588778458346\\
0.1	0.00299605565888641\\
1	0.00300557747718311\\
10	0.00314677536771405\\
100	0.00345767554055898\\
1000	0.00394328091489954\\
};
\addlegendentry{SMG-QMCL-ROM}

\addplot [color=green!80, dashdotted, line width=2.0pt, mark size=4.0pt, mark=triangle, mark options={solid, green!80}]
  table[row sep=crcr]{%
0.001	0.00299588260127251\\
0.01	0.00299588778458346\\
0.1	0.00299605565888641\\
1	0.00300557747718311\\
10	0.00314677536771405\\
100	0.00345767554055898\\
1000	0.00394328091489954\\
};
\addlegendentry{Galerkin-BQ-ROM}

\end{axis}

\begin{axis}[%
width=1.282\fheight,
height=0.756\fheight,
at={(-0.167\fheight,-0.11\fheight)},
scale only axis,
xmin=0,
xmax=1,
ymin=0,
ymax=1,
axis line style={draw=none},
ticks=none,
axis x line*=bottom,
axis y line*=left
]
\end{axis}
\end{tikzpicture}%
\caption{Relative projection error in $\bq$}
\label{fig:gamma_q}
    \end{subfigure}
    \hspace{1cm}
    \begin{subfigure}{.26\textwidth}
           \setlength\fheight{5 cm}
           \setlength\fwidth{\textwidth}
\raisebox{-59mm}{
%
%
\begin{tikzpicture}

\begin{axis}[%
width=0.694\fheight,
height=0.59\fheight,
at={(0\fheight,0\fheight)},
scale only axis,
xmode=log,
xmin=0.001,
xmax=1000,
xminorticks=true,
xtick={1e-3,1e-1,1e1,1e3},
xlabel style={font=\color{white!15!black}},
xlabel={Regularization parameter $\gamma_q$},
ymode=log,
ymin=0.01,
ymax=0.1,
yminorticks=true,
ylabel style={font=\color{white!15!black}},
ylabel={Relative projection error},
axis background/.style={fill=white},
xmajorgrids,
ymajorgrids,
legend style={at={(0.31,0.327)}, anchor=south west, legend cell align=left, align=left, draw=white!15!black}
]
\addplot [color=blue!80, dotted, line width=2.0pt, mark size=4.0pt, mark=o, mark options={solid, blue!80}]
  table[row sep=crcr]{%
0.001	0.0480138789621739\\
0.01	0.0480138789621739\\
0.1	0.0480138789621739\\
1	0.0480138789621739\\
10	0.0480138789621739\\
100	0.0480138789621739\\
1000	0.0480138789621739\\
};

\addplot [color=red!80, dashed, line width=2.0pt, mark size=4pt, mark=square, mark options={solid, red!80}]
  table[row sep=crcr]{%
0.001	0.0196760073596324\\
0.01	0.0196766274122813\\
0.1	0.0196828588851075\\
1	0.019746639178449\\
10	0.0202321939688371\\
100	0.0216304325819215\\
1000	0.0243637273875739\\
};

\addplot [color=green!80, dashdotted, line width=2.0pt, mark size=4.0pt, mark=triangle, mark options={solid, green!80}]
  table[row sep=crcr]{%
0.001	0.0187673913742443\\
0.01	0.0187673913706068\\
0.1	0.0187673913342334\\
1	0.0187673909706555\\
10	0.0187673873504488\\
100	0.0187673527001747\\
1000	0.0187671559891852\\
};

\end{axis}

\begin{axis}[%
width=1.282\fheight,
height=0.756\fheight,
at={(-0.167\fheight,-0.11\fheight)},
scale only axis,
xmin=0,
xmax=1,
ymin=0,
ymax=1,
axis line style={draw=none},
ticks=none,
axis x line*=bottom,
axis y line*=left
]
\end{axis}
\end{tikzpicture}%
}
\caption{Relative projection error in $\bp$}
\label{fig:gamma_p}
    \end{subfigure}
       \hspace{1cm}
        \begin{subfigure}{.26\textwidth}
           \setlength\fheight{5 cm}
           \setlength\fwidth{\textwidth}
\raisebox{-59mm}{
%
%
\begin{tikzpicture}

\begin{axis}[%
width=0.694\fheight,
height=0.59\fheight,
at={(0\fheight,0\fheight)},
scale only axis,
xmode=log,
xmin=0.001,
xmax=1000,
xminorticks=true,
xtick={1e-3,1e-1,1e1,1e3},
xlabel style={font=\color{white!15!black}},
xlabel={Regularization parameter $\gamma_q$},
ymode=log,
ymin=0.01,
ymax=0.1,
yminorticks=true,
ylabel style={font=\color{white!15!black}},
ylabel={Relative state error},
axis background/.style={fill=white},
xmajorgrids,
ymajorgrids,
legend style={at={(0.38,0.383)}, anchor=south west, legend cell align=left, align=left, draw=white!15!black}
]
\addplot [color=blue!80, dotted, line width=2.0pt, mark size=4.0pt, mark=o, mark options={solid, blue!80}]
  table[row sep=crcr]{%
0.001	0.0925823737485728\\
0.01	0.0925823737485728\\
0.1	0.0925823737485728\\
1	0.0925823737485728\\
10	0.0925823737485728\\
100	0.0925823737485728\\
1000	0.0925823737485728\\
};

\addplot [color=red!80, dashed, line width=2.0pt, mark size=4pt, mark=square, mark options={solid, red!80}]
  table[row sep=crcr]{%
0.001	0.0213707276185599\\
0.01	0.0213707778523377\\
0.1	0.0213716926556343\\
1	0.0214072931878077\\
10	0.0219446716993026\\
100	0.0234751209292225\\
1000	0.0268668252481628\\
};

\addplot [color=green!80, dashdotted, line width=2.0pt, mark size=4.0pt, mark=triangle, mark options={solid, green!80}]
  table[row sep=crcr]{%
0.001	0.0217104998205157\\
0.01	0.0217101213393164\\
0.1	0.0217068581821765\\
1	0.0217096146896626\\
10	0.0219272916076884\\
100	0.0212504079001584\\
1000	0.0219776099708309\\
};

\end{axis}

\begin{axis}[%
width=1.227\fheight,
height=0.723\fheight,
at={(-0.16\fheight,-0.08\fheight)},
scale only axis,
xmin=0,
xmax=1,
ymin=0,
ymax=1,
axis line style={draw=none},
ticks=none,
axis x line*=bottom,
axis y line*=left
]
\end{axis}
\end{tikzpicture}%
}
\caption{Relative state error (training)}
\label{fig:gamma_state}
    \end{subfigure}
\caption{ {\Rb{One-dimensional linear wave equation (sensitivity to regularization parameters). The accuracy of the proposed ROMs with reduced dimension $2r=16$ does not change with the regularization parameter up to $\gamma_q=10^1$. For $\gamma_q>10^{1}$, we observe a marginal decrease in the accuracy for both the SMG-QMCL-ROM and the Galerkin-BQ-ROM.}} }
 \label{fig:LW_gamma}
\end{figure}
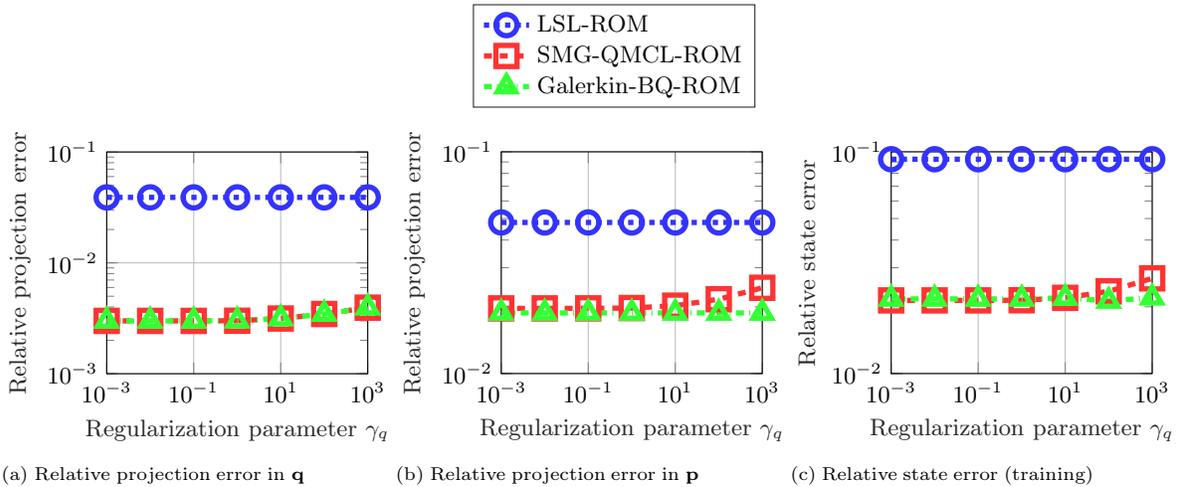

{\Rb{We study the effect of the regularization parameters $\gamma_q$ and $\gamma_p$ on the accuracy of the proposed ROMs in Figure~\ref{fig:LW_gamma}. There, we fix the reduced dimension to $2r=16$, set $\gamma_p=\gamma_q$, and vary the regularization parameter $\gamma_q$ from $10^{-3}$ to $10^3$. In Figure~\ref{fig:gamma_q}, we observe that the relative projection error in $\bq$ for both approaches remain approximately constant up to $\gamma_q=10^{1}$. However, as $\gamma_q$ increases further, the relative projection error increases in magnitude.  For the comparison in Figure~\ref{fig:gamma_p}, we observe that the projection error in $\bp$ increases marginally for the SMG-QMCL approach whereas the Galerkin-BQ approach demonstrates nearly constant projection error. In Figure~\ref{fig:gamma_state}, we observe that both the SMG-QMCL and the Galerkin-BQ approach yield ROMs with approximately same accuracy up to $\gamma_q=1$. For $\gamma_q>1$, we observe a marginal decrease in the accuracy with increasing $\gamma_q$ for both approaches. Overall, the relative state error values for SMG-QMCL and Galerkin-BQ ROMs of different dimensions remain largely unaffected by the regularization parameter in the investigated range and hence we choose $\gamma_q$ and $\gamma_p$ values that lead to accurate and stable ROMs over a range of reduced dimensions. The results from this sensitivity study illustrate that the proposed approaches provide higher accuracy than the linear symplectic subspace approach over a wide range of regularization parameters. Therefore, the regularization parameters used for the numerical results in Section~\ref{sec:results_lw} are quite robust and can be used to derive accurate and stable ROMs without extensive fine-tuning. }}
\subsubsection{Time extrapolation study}
\label{sec:time_extrapolation}
{\Rb{A key motivation for symplectic model reduction using data-driven quadratic manifolds is to obtain low-dimensional ROMs that can capture the periodic behavior in an accurate and stable manner.}} In this study, we fix the scalar parameter in~\eqref{eq:wave_FOM} to $\mu=0.5$ and build a training dataset of snapshots over one cycle by integrating the Hamiltonian FOM with the implicit midpoint method until time $t=4$ with a fixed time step $\Delta t=10^{-3}$. Note that we use the same training dataset with the same fixed time step $\Delta t=10^{-3}$ as in the motivational example in Section~\ref{sec:motivation}. {\Rb{For this time extrapolation study, we found $\gamma_q=\gamma_p=10^6$ to be a robust choice for the proposed approaches}}. For test data, we consider the FOM solution snapshots over the next nine cycles, i.e., from $t=4$ to $t=40$. For this time extrapolation study, we center the trajectory snapshot data about the initial condition with $\bq_{\text{ref}}=\bq(0)$ and $\bp_{\text{ref}}=\bp(0)$.

\begin{figure}[tbp]
\small
\captionsetup[subfigure]{oneside,margin={1.8cm,0 cm}}
\begin{subfigure}{.42\textwidth}
       \setlength\fheight{6 cm}
        \setlength\fwidth{\textwidth}
%
%
\begin{tikzpicture}

\begin{axis}[%
width=0.951\fheight,
height=0.536\fheight,
at={(0\fheight,0\fheight)},
scale only axis,
xmin=0,
xmax=40,
xlabel style={font=\color{white!15!black}},
xlabel={Reduced dimension $2r$},
ymode=log,
ymin=0.002,
ymax=4,
yminorticks=true,
ylabel style={font=\color{white!15!black}},
ylabel={Relative state error},
axis background/.style={fill=white},
xmajorgrids,
ymajorgrids,
legend style={at={(0.125,1.4)}, anchor=south west, legend cell align=left, align=left, draw=white!15!black}
]
\addplot [color=blue!80, dotted, line width=2.0pt, mark size=4.0pt, mark=o, mark options={solid, blue!80}]
  table[row sep=crcr]{%
4	0.958616195405763\\
8	0.815635119024444\\
12	0.501309800593042\\
16	0.383008509540221\\
20	0.0571150856934921\\
24	0.0382551168490792\\
28	0.01668784842278\\
32	0.0126639793195574\\
36	0.0105272831227202\\
40	0.00726266476603372\\
};
\addlegendentry{LSL-ROM}

\addplot [color=red!80, dashed, line width=2.0pt, mark size=4pt, mark=square, mark options={solid, red!80}]
  table[row sep=crcr]{%
4	1.39350505318364\\
8	1.31776761577814\\
12	1.26736031766801\\
16	0.64230990071716\\
20	0.0714320822789134\\
24	0.027303497376185\\
28	0.0182391683764995\\
32	0.0150586194520962\\
36	0.00998934433154718\\
40	0.00773150888634736\\
};
\addlegendentry{SMG-QMCL-ROM}

\addplot [color=green!80, dashdotted, line width=2.0pt, mark size=4.0pt, mark=triangle, mark options={solid, green!80}]
  table[row sep=crcr]{%
4	1.77892730418935\\
8	1.55737234408099\\
12	0.634260415458255\\
16	0.113410950590196\\
20	0.0216894700787247\\
24	0.0550534172417606\\
28	0.010682711906898\\
32	0.00774519411570462\\
36	0.00856635380229549\\
40	0.00647856968736228\\
};
\addlegendentry{Galerkin-BQ-ROM}

\end{axis}
\end{tikzpicture}%
\caption{State error (training data)}
\label{fig: extrapolation_state}
    \end{subfigure}
    \hspace{0.9cm}
    \begin{subfigure}{.42\textwidth}
           \setlength\fheight{6 cm}
           \setlength\fwidth{\textwidth}
\input{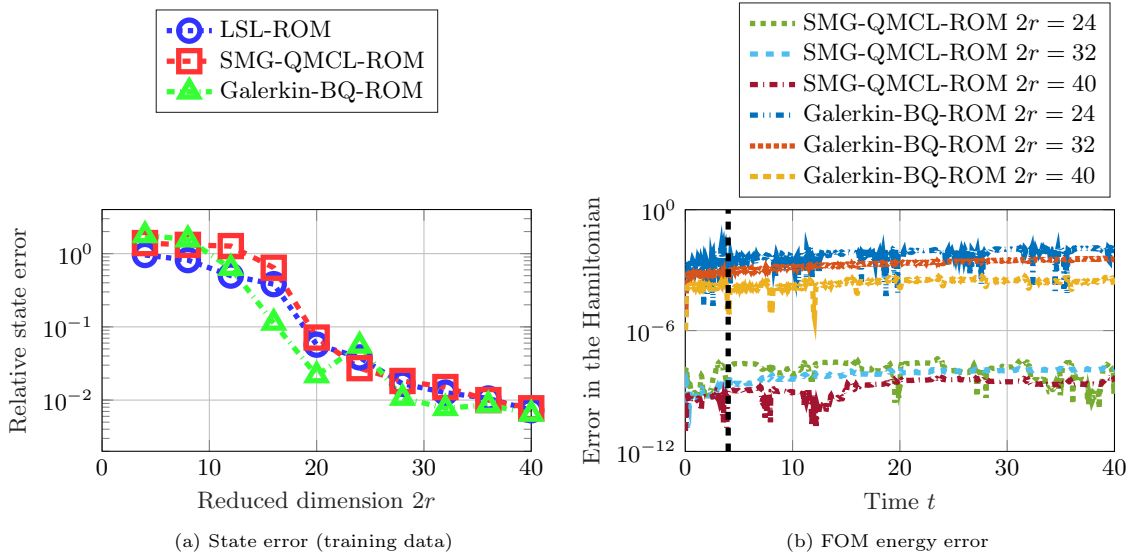}
\caption{FOM energy error}
\label{fig:LW_FOM_extrapolation}
    \end{subfigure}
   \caption{{\Ra{One-dimensional linear wave equation (time extrapolation study).  Plot (a) shows that LSL-ROMs, Galerkin-BQ-ROMs, and SMG-QMCL-ROMs achieve similar accuracy in the training data regime. Plot(b) compares the energy error behavior for ROMs of different sizes. The nonlinear SMG-QMCL-ROMs exhibit bounded error in the Hamiltonian~\eqref{eq:err_hamiltonian} outside the training data regime. In contrast, the Hamiltonian errors for the quadratic Galerkin-BQ-ROMs slowly grow with time.  The dashed black line indicates end of training time interval at $t=4$.}}}
 \label{fig:motivation_compare}
\end{figure}
{\Ra{We compare the accuracy and the energy error performance of the proposed ROMs in Figure~\ref{fig:motivation_compare}. In Figure~\ref{fig: extrapolation_state}, we observe that both the SMG-QMCL and the Galerkin-BQ approach yield ROMs with similar accuracy in the training data regime. Compared with Figure~\ref{fig: separate_state}, we observe that the proposed approaches yield stable ROMs that demonstrate a decrease in state error with increasing reduced dimension. The energy error plots in Figure~\ref{fig:LW_FOM_extrapolation} demonstrate stability outside the training data regime for both the SMG-QMCL-ROMs and the Galerkin-BQ-ROMs of different reduced dimensions. For both approaches, we observe that the error in the Hamiltonian decreases marginally with an increase in the reduced dimension $2r$. All three SMG-QMCL-ROMs exhibit bounded energy error and the error in the Hamiltonian for all three ROMs remain below an absolute value of $10^{-7}$ outside the training time interval.}} Despite the {\Rc{Hamiltonian error}} for the Galerkin-BQ-ROM remaining below an absolute value of $10^{-1}$ over the entire observation interval $t \in [0,40]$, we do observe a slow growth because of its approximately symplectic structure. 

\begin{figure}[tbp]
\centering
\begin{subfigure}{.22\textwidth}
\includegraphics[width=\linewidth]{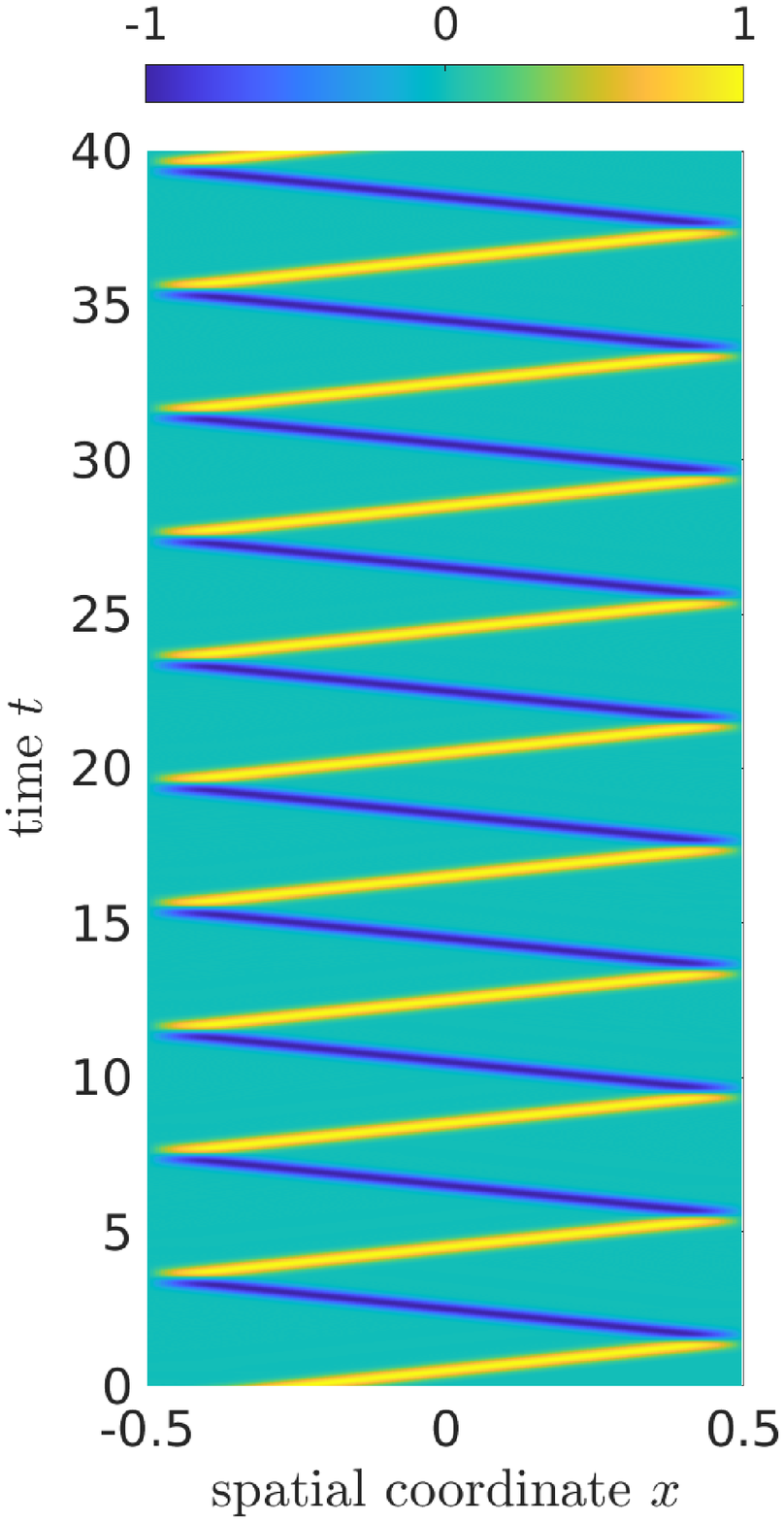}
\subcaption{FOM solution}
\label{fig:time_fom}
 \end{subfigure}
\begin{subfigure}{.22\textwidth}
\includegraphics[width=\linewidth]{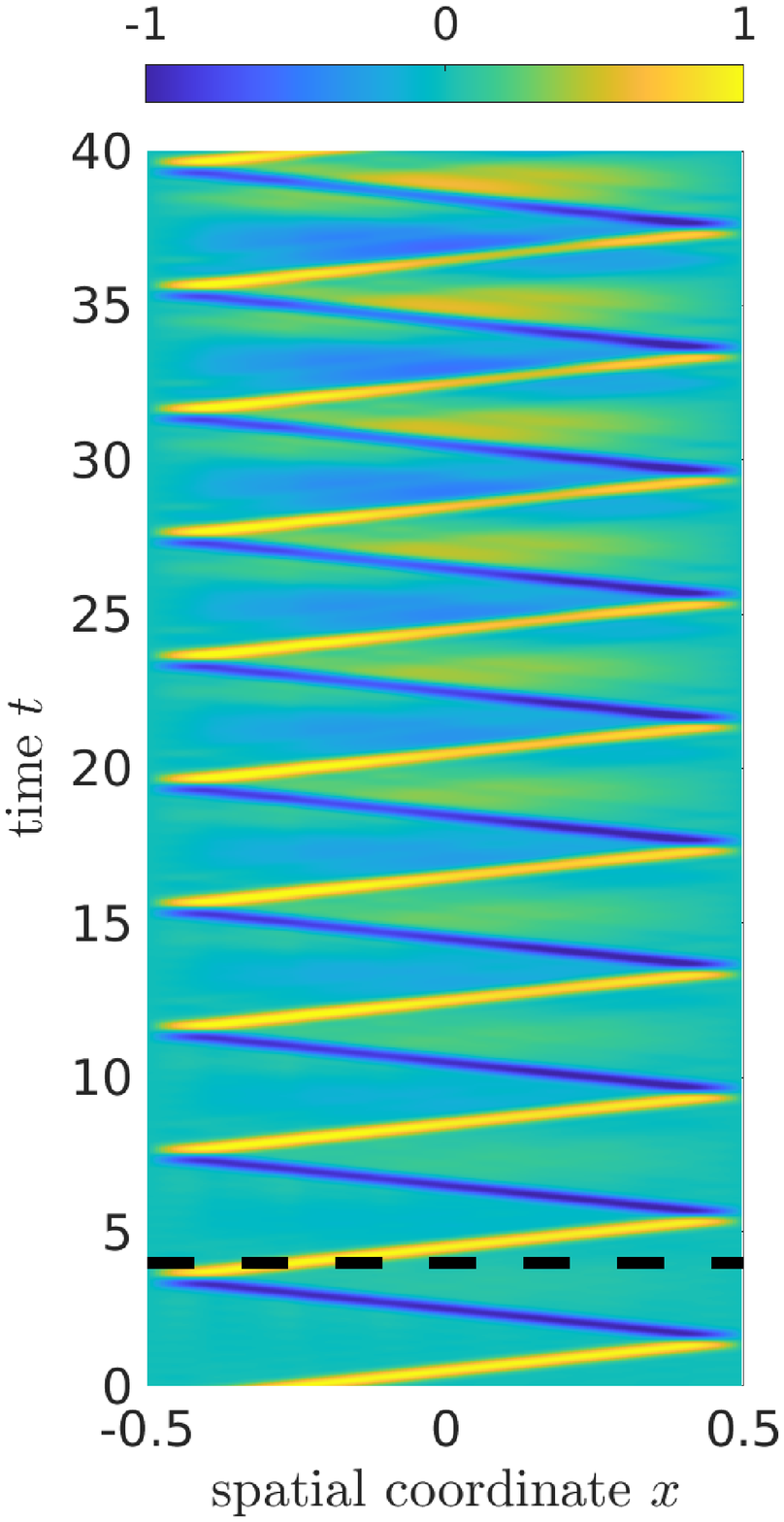}
\subcaption{LSL-ROM}
\label{fig:time_linear}
\end{subfigure}
\begin{subfigure}{.22\textwidth}
\includegraphics[width=\linewidth]{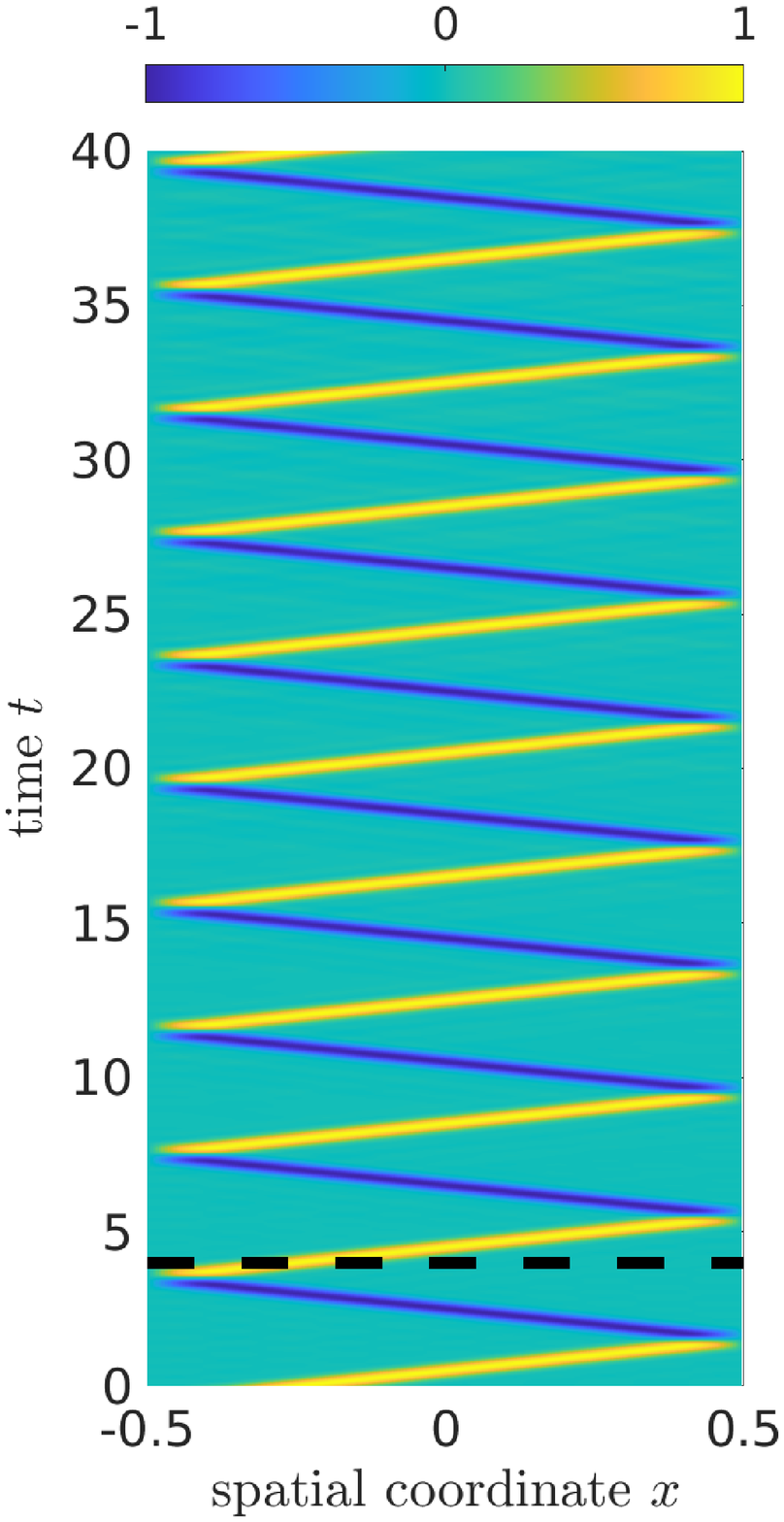}
\subcaption{SMG-QMCL-ROM}
\label{fig:time_jac}
\end{subfigure}
 \begin{subfigure}{.22\textwidth}
\includegraphics[width=\linewidth]{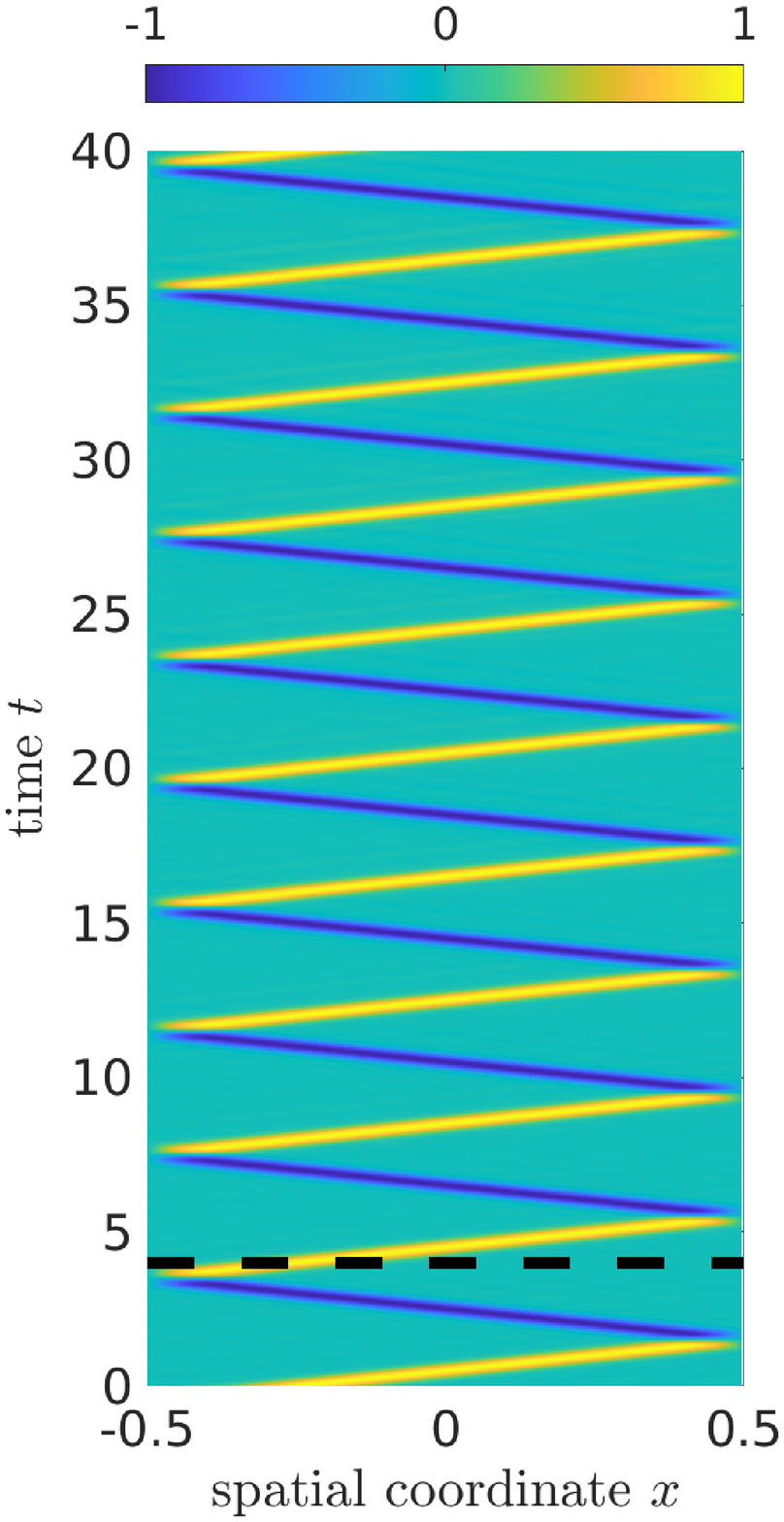}
\subcaption{Galerkin-BQ-ROM}
\label{fig:time_pod}
 \end{subfigure}
    \caption{{\Rb{One-dimensional linear wave equation (time extrapolation study).}} Plots show the numerical approximation of the 
 FOM solution using low-dimensional ($2r=24$) LSL-ROM, SMG-QMCL-ROM, and Galerkin-BQ-ROM for $t \in [0,40]$. The dashed black line indicates the end of the training time interval at $t=4$. Plot (b) shows that the LSL-ROM yields inaccurate solutions after $t=8$ whereas plot (c) and plot (d) show that SMG-QMCL-ROM and Galerkin-BQ-ROM provide accurate approximate solutions even at $t=40$.}
\label{fig:LW_spacetime}
\end{figure}

The space-time evolution of the FOM solution field is compared with the space-time evolution of the approximate solutions obtained using ROMs of dimension $2r=24$ in Figure~\ref{fig:LW_spacetime}. While all three ROMs approximate the FOM solution accurately in the training data regime, we observe that the linear symplectic ROM exhibits spurious oscillations that gradually become more pronounced as we march forward in time in Figure~\ref{fig:time_linear}. In Figure~\ref{fig:time_jac} and Figure~\ref{fig:time_pod}, we observe that both proposed approaches provide stable and accurate approximate solutions at time $t=40$, which is $900\%$ outside the training time interval. {\Rb{Unlike the conventional quadratic manifolds approach in Section~\ref{sec:motivation} which can lead to unstable ROMs, the long-time stable and accurate predictions of future cycles in this time extrapolation study highlights the utility of the proposed methods for high-dimensional systems with periodic behavior.}}  


\subsection{{\Rb{Parametrized two-dimensional nonlinear wave equation}}} \label{sec:wave_2d}
{\Rb{ For a second numerical example, we consider the two-dimensional nonlinear wave equation with parametric nonlinearity. This parametric model problem is similar to the two-dimensional nonlinear wave equation example in~\cite{jiang2020linearly}. Let $\Omega=(-10,10) \times (-10,10) \subset \real^2$ be the spatial domain and consider the parametrized two-dimensional nonlinear wave equation 
\begin{equation}\label{eq:non_pde}
\frac{\partial^2 \varphi}{\partial t^2}(x,y,t;\mu)=\frac{\partial^2 \varphi}{\partial x^2} (x,y,t;\mu)+ \frac{\partial^2 \varphi}{\partial y^2}(x,y,t;\mu) - \mu\varphi(x,y,t;\mu)^3, 
\end{equation}
with the state $\varphi(x,y,t;\mu)$ at $(x,y)\in \Omega$, time $t\in(0,T]$, and the scalar parameter $\mu \in \mathcal P=[0.1,3]$. Periodic boundary conditions
\begin{equation}
\varphi(-10,y,t;\mu)=\varphi(10,y,t;\mu), \qquad \varphi(x,-10,t;\mu)=\varphi(x,10,t;\mu),
\end{equation}
are imposed for $t \in (0,T]$ and $\mu \in \mathcal P$. 
\subsubsection{Hamiltonian PDE formulation and FOM implementation details}
We rewrite~\eqref{eq:non_pde} as an infinite-dimensional canonical Hamiltonian system with $q(x,y,t;\mu):=\varphi(x,y,t;\mu)$ and $p(x,y,t;\mu):=\partial \varphi(x,y,t;\mu)/\partial t$. The associated space-time continuous Hamiltonian functional is
\begin{multline}\label{eq:fham_non}
\mathcal H(q(x,y,t;\mu),p(x,y,t;\mu);\mu)=\frac{1}{2}\int_{\Omega} \left[ p(x,y,t;\mu)^2 + \left(\frac{\partial}{\partial x}q(x,y,t;\mu)\right)^2 +  \left(\frac{\partial}{\partial y}q(x,y,t;\mu)\right)^2 \right] \mathrm{d}x\mathrm{d}y \\
+ \int_{\Omega} \frac{\mu}{4}q(x,y,t;\mu)^4  \mathrm{d}x\mathrm{d}y
\end{multline}
and the original nonlinear PDE can be recast as a canonical Hamiltonian PDE
\begin{align*}
\frac{\partial }{\partial t}q(x,y,t;\mu)&=\frac{\textcolor{black}{\delta} \mathcal H}{\textcolor{black}{\delta} p}(q,p;\mu)=p(x,y,t;\mu),\\
\frac{\partial }{\partial t}p(x,y,t;\mu)&=-\frac{\textcolor{black}{\delta} \mathcal H}{\textcolor{black}{\delta} q}(q,p;\mu)=\frac{\partial ^2}{\partial x^2}q(x,y,t;\mu)+\frac{\partial ^2}{\partial y^2}q(x,y,t;\mu)-\mu q(x,y,t;\mu)^3.
\end{align*}
We discretize the two-dimensional spatial domain $\Omega$ with $n_x=n_y=100$ equally spaced grid points in both spatial directions leading to system states of dimension $2n$ with $n=n_xn_y=10,000$. Using a finite difference scheme in both spatial directions, we obtain the following space-discretized Hamiltonian
\begin{multline}\label{eq:dham_non}
H(\bq(t;\mu),\bp(t;\mu);\mu)=\Delta x \Delta y \sum_{i,j=1}^{n_x,n_y} \left[ p_{i,j}(t;\mu)^2 + \frac{ (q_{i+1,j}(t;\mu)-q_{i,j}(t;\mu))^2}{4\Delta x^2} + \frac{ (q_{i,j}(t;\mu)-q_{i-1,j}(t;\mu))^2}{4\Delta x^2}\right] \\
+\Delta x \Delta y  \sum_{i,j=1}^{n_x,n_y} \left[ \frac{ (q_{i,j+1}(t;\mu)-q_{i,j}(t;\mu))^2}{4\Delta y^2} + \frac{ (q_{i,j}(t;\mu)-q_{i,j-1}(t;\mu))^2}{4\Delta y^2} +  \frac{\mu}{4}q_{i,j}(t;\mu)^4 \right],
\end{multline}
with
\begin{align*}
    q_{i,j}(t;\mu)&:=\varphi(x_i,y_j,t;\mu),&
\bq(t;\mu)&=(q_{1,1}(t;\mu), \ldots, q_{n_x,n_y}(t;\mu) )^\top \in \real^n,\\
    p_{i,j}(t;\mu)&:=\frac{\partial}{\partial t} \varphi(x_i,y_j,t;\mu),&
\bp(t;\mu)&=(p_{1,1}(t;\mu), \ldots, p_{n_x,n_y}(t;\mu) )^\top \in \real^n.
\end{align*}
The corresponding parametrized nonlinear Hamiltonian FOM is
 \begin{equation}
   \dot{\by}(t;\mu) =\begin{pmatrix}
       \dot{\bq}(t;\mu)  \\
       \dot{\bp}(t;\mu) 
     \end{pmatrix}=\Jn \nabla_{ \by }\dH(\by(t;\mu))
     =\begin{pmatrix}
       \bzero & \In \\
       \mathbf{D}_{\text{fd},2d} & \bzero
     \end{pmatrix}
     \begin{pmatrix}
       \bq(t;\mu) \\
       \bp(t;\mu)
     \end{pmatrix} - \mu\begin{pmatrix}
       \bzero \\
       \bq(t;\mu)^3
     \end{pmatrix} ,
     \label{eq:wave_2d_FOM}
 \end{equation}
where $\mathbf{D}_{\text{fd},2d}$ denotes the finite difference approximation in the two-dimensional setting and the vector $ \bq(t;\mu)^3\in \real^{n}$ contains as components the entry-wise cubic exponential of the generalized state vector $\bq(t;\mu)$ .}}
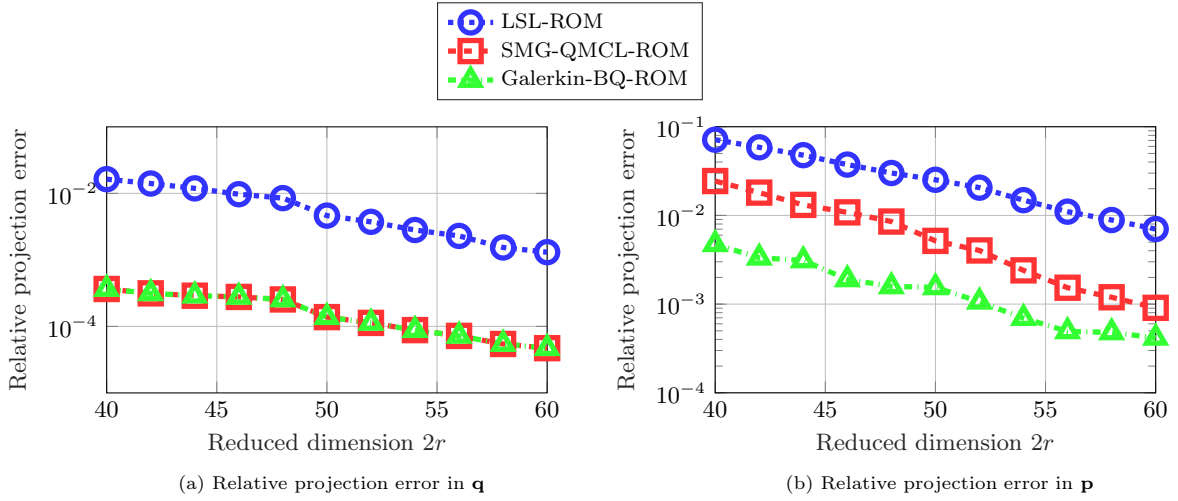
\begin{figure}[tbp]
\small
\captionsetup[subfigure]{oneside,margin={1.8cm,0 cm}}
\begin{subfigure}{.45\textwidth}
       \setlength\fheight{6 cm}
        \setlength\fwidth{\textwidth}
%
%
\begin{tikzpicture}

\begin{axis}[%
width=0.976\fheight,
height=0.59\fheight,
at={(0\fheight,0\fheight)},
scale only axis,
xmin=40,
xmax=60,
xlabel style={font=\color{white!15!black}},
xlabel={Reduced dimension $2r$},
ymode=log,
ymin=1e-05,
ymax=0.1,
yminorticks=true,
ylabel style={font=\color{white!15!black}},
ylabel={Relative projection error},
axis background/.style={fill=white},
xmajorgrids,
ymajorgrids,
yminorgrids,
legend style={at={(0.75,1.1)}, anchor=south west, legend cell align=left, align=left, draw=white!15!black},
legend style={font=\footnotesize}
]
\addplot [color=blue!80, dotted, line width=2.0pt, mark size=4.0pt, mark=o, mark options={solid, blue!80}]
  table[row sep=crcr]{%
40	0.0163367057803846\\
42	0.0140581634446982\\
44	0.0118868897599276\\
46	0.00971173158159623\\
48	0.00847705440508853\\
50	0.00462755451014946\\
52	0.00372989580202568\\
54	0.00281142070679676\\
56	0.00229143575761105\\
58	0.00153662427279414\\
60	0.00130568875104336\\
};
\addlegendentry{LSL-ROM}

\addplot [color=red!80, dashed, line width=2.0pt, mark size=4pt, mark=square, mark options={solid, red!80}]
  table[row sep=crcr]{%
40	0.000365021081994925\\
42	0.000313883668213957\\
44	0.000289543934524435\\
46	0.000277409935389811\\
48	0.000254197861535119\\
50	0.000137851218219085\\
52	0.000112827838344225\\
54	8.79180661171069e-05\\
56	7.19893604029802e-05\\
58	5.44179956718145e-05\\
60	4.70503832326577e-05\\
};
\addlegendentry{SMG-QMCL-ROM}

\addplot [color=green!80, dashdotted, line width=2.0pt, mark size=4.0pt, mark=triangle, mark options={solid, green!80}]
  table[row sep=crcr]{%
40	0.000365021081874528\\
42	0.000313883667707765\\
44	0.000289543934307879\\
46	0.000277409935237766\\
48	0.000254197861434462\\
50	0.000137851217885768\\
52	0.000112827838350184\\
54	8.79180663167833e-05\\
56	7.19893607310725e-05\\
58	5.44179956006357e-05\\
60	4.70503831266323e-05\\
};
\addlegendentry{Galerkin-BQ-ROM}

\end{axis}

\begin{axis}[%
width=1.259\fheight,
height=0.743\fheight,
at={(-0.164\fheight,-0.097\fheight)},
scale only axis,
xmin=0,
xmax=1,
ymin=0,
ymax=1,
axis line style={draw=none},
ticks=none,
axis x line*=bottom,
axis y line*=left
]
\end{axis}
\end{tikzpicture}%
\caption{Relative projection error in $\bq$}
\label{fig:NW_proj_q}
    \end{subfigure}
    \hspace{0.9cm}
    \begin{subfigure}{.45\textwidth}
           \setlength\fheight{6 cm}
           \setlength\fwidth{\textwidth}
\raisebox{-60mm}{
%
%
\begin{tikzpicture}

\begin{axis}[%
width=0.976\fheight,
height=0.59\fheight,
at={(0\fheight,0\fheight)},
scale only axis,
xmin=40,
xmax=60,
xlabel style={font=\color{white!15!black}},
xlabel={Reduced dimension $2r$},
ymode=log,
ymin=0.0001,
ymax=0.1,
yminorticks=true,
ylabel style={font=\color{white!15!black}},
ylabel={Relative projection error},
axis background/.style={fill=white},
xmajorgrids,
ymajorgrids,
legend style={at={(0.171,0.16)}, anchor=south west, legend cell align=left, align=left, draw=white!15!black}
]
\addplot [color=blue!80, dotted, line width=2.0pt, mark size=4.0pt, mark=o, mark options={solid, blue!80}]
  table[row sep=crcr]{%
40	0.0714280925845315\\
42	0.0586249276358492\\
44	0.0478033590830446\\
46	0.0374315752591774\\
48	0.0300817691677405\\
50	0.025227355677965\\
52	0.0205682760703789\\
54	0.0148632384754037\\
56	0.0110306048927208\\
58	0.00887621864000561\\
60	0.00702214607979646\\
};

\addplot [color=green!80, dashdotted, line width=2.0pt, mark size=4.0pt, mark=triangle, mark options={solid, green!80}]
  table[row sep=crcr]{%
40	0.00472731425145413\\
42	0.00333082907438541\\
44	0.0031081942507089\\
46	0.00189999614802119\\
48	0.00158404968785584\\
50	0.00154436567396206\\
52	0.00106364015777339\\
54	0.000695356999683315\\
56	0.000492822214627649\\
58	0.000479299187347268\\
60	0.000416497824039841\\
};

\addplot [color=red!80, dashed, line width=2.0pt, mark size=4pt, mark=square, mark options={solid, red!80}]
  table[row sep=crcr]{%
40	0.0244197659954281\\
42	0.018081222855358\\
44	0.0131415986409003\\
46	0.0108247092473821\\
48	0.00853186887757204\\
50	0.00516624396360892\\
52	0.00402757160326894\\
54	0.00240069420460724\\
56	0.0015299776012956\\
58	0.00119457683692629\\
60	0.000908110620881029\\
};

\end{axis}
\end{tikzpicture}%
}
\caption{Relative projection error in $\bp$}
\label{fig:NW_proj_p}
    \end{subfigure}
\caption{{\Rb{Two-dimensional nonlinear wave equation. Both the QMCL and the BQ state approximation yield a substantially lower relative projection error~\eqref{eq:err_rel_proj_train} than the linear symplectic subspaces for $\bq$ and $\bp$ variables. The regularization factors are chosen to be $\gamma_q=\gamma_p=10^{-1}$.}}}
 \label{fig:NW_proj}
\end{figure}
\begin{figure}[tbp]
\small
\captionsetup[subfigure]{oneside,margin={1.8cm,0 cm}}
\begin{subfigure}{.45\textwidth}
       \setlength\fheight{6 cm}
        \setlength\fwidth{\textwidth}
%
%
\begin{tikzpicture}

\begin{axis}[%
width=0.976\fheight,
height=0.59\fheight,
at={(0\fheight,0\fheight)},
scale only axis,
xmin=40,
xmax=60,
xlabel style={font=\color{white!15!black}},
xlabel={Reduced dimension $2r$},
ymode=log,
ymin=0.0001,
ymax=0.5,
yminorticks=true,
ylabel style={font=\color{white!15!black}},
ylabel={Relative state error},
axis background/.style={fill=white},
xmajorgrids,
ymajorgrids,
legend style={at={(0.75,1.1)}, anchor=south west, legend cell align=left, align=left, draw=white!15!black},
legend style={font=\footnotesize}
]
\addplot [color=blue!80, dotted, line width=2.0pt, mark size=4.0pt, mark=o, mark options={solid, blue!80}]
  table[row sep=crcr]{%
40	0.0725373475987685\\
42	0.0638150624346438\\
44	0.0489444892482522\\
46	0.038480448293741\\
48	0.0274193363549455\\
50	0.0251054530210967\\
52	0.019718584914568\\
54	0.0138235100453645\\
56	0.00978055369336385\\
58	0.00814266247787982\\
60	0.00621276449482052\\
};
\addlegendentry{LSL-ROM}

\addplot [color=red!80, dashed, line width=2.0pt, mark size=4pt, mark=square, mark options={solid, red!80}]
  table[row sep=crcr]{%
40	0.0181676490697177\\
42	0.0109559356255111\\
44	0.00857475975604421\\
46	0.00745542044084839\\
48	0.00592354648319891\\
50	0.00517150437335792\\
52	0.0034464198212887\\
54	0.00176155847505971\\
56	0.00111744475042563\\
58	0.00108754824334084\\
60	0.000815095493950784\\
};
\addlegendentry{SMG-QMCL-ROM}

\addplot [color=green!80, dashdotted, line width=2.0pt, mark size=4.0pt, mark=triangle, mark options={solid, green!80}]
  table[row sep=crcr]{%
40	0.00527122500828252\\
42	0.00433218379302466\\
44	0.0030414721815242\\
46	0.00216139918802086\\
48	0.00170952267634643\\
50	0.0018781877920258\\
52	0.00111826572144721\\
54	0.000728024332289327\\
56	0.000464630785873502\\
58	0.000455024551389656\\
60	0.000378643233744798\\
};
\addlegendentry{Galerkin-BQ-ROM}

\end{axis}
\end{tikzpicture}%
\caption{Training parameters}
\label{fig:NW_state_train}
    \end{subfigure}
    \hspace{0.9cm}
    \begin{subfigure}{.45\textwidth}
           \setlength\fheight{6 cm}
           \setlength\fwidth{\textwidth}
\raisebox{-59mm}{
%
%
\begin{tikzpicture}

\begin{axis}[%
width=0.976\fheight,
height=0.59\fheight,
at={(0\fheight,0\fheight)},
scale only axis,
xmin=40,
xmax=60,
xlabel style={font=\color{white!15!black}},
xlabel={Reduced dimension $2r$},
ymode=log,
ymin=1e-3,
ymax=.5,
yminorticks=true,
ylabel style={font=\color{white!15!black}},
ylabel={Relative state error},
axis background/.style={fill=white},
xmajorgrids,
ymajorgrids,
legend style={at={(0.529,0.724)}, anchor=south west, legend cell align=left, align=left, draw=white!15!black}
]
\addplot [color=blue!80, dotted, line width=2.0pt, mark size=4.0pt, mark=o, mark options={solid, blue!80}]
  table[row sep=crcr]{%
40	0.106053345219061\\
42	0.101276443157002\\
44	0.0957658295129943\\
46	0.0915642259818487\\
48	0.0868812873510938\\
50	0.0710586511426693\\
52	0.0409343313705578\\
54	0.0276265960780708\\
56	0.0242270255816793\\
58	0.0142188408557695\\
60	0.0129390835037293\\
};

\addplot [color=green!80, dashdotted, line width=2.0pt, mark size=4.0pt, mark=triangle, mark options={solid, green!80}]
  table[row sep=crcr]{%
40	0.0389516029112314\\
42	0.0322282942237123\\
44	0.0233987182984175\\
46	0.0217877271706315\\
48	0.0170733930203921\\
50	0.0240153509845101\\
52	0.0101496572010465\\
54	0.00639331010434201\\
56	0.00510468940838075\\
58	0.003114245087014\\
60	0.00281735452537096\\
};

\addplot [color=red!80, dashed, line width=2.0pt, mark size=4pt, mark=square, mark options={solid, red!80}]
  table[row sep=crcr]{%
40	0.0420799170840403\\
42	0.0367032238643159\\
44	0.0317152040818589\\
46	0.0301180713517823\\
48	0.0271658504601051\\
50	0.0255573542040006\\
52	0.0138430921288641\\
54	0.00984074048876972\\
56	0.00894929741975595\\
58	0.00488128784019317\\
60	0.00457544150517699\\
};

\end{axis}
\end{tikzpicture}%
}
\caption{Test parameters}
\label{fig:NW_state_test}
    \end{subfigure}
\caption{{\Rb{Two-dimensional nonlinear wave equation. The proposed SMG-QMCL and Galerkin-BQ ROMs  achieve lower state error~\eqref{eq:err_rel_sim_train} than the LSL-ROMs for both training and test parameters.}}}
 \label{fig:NW_state}
\end{figure}
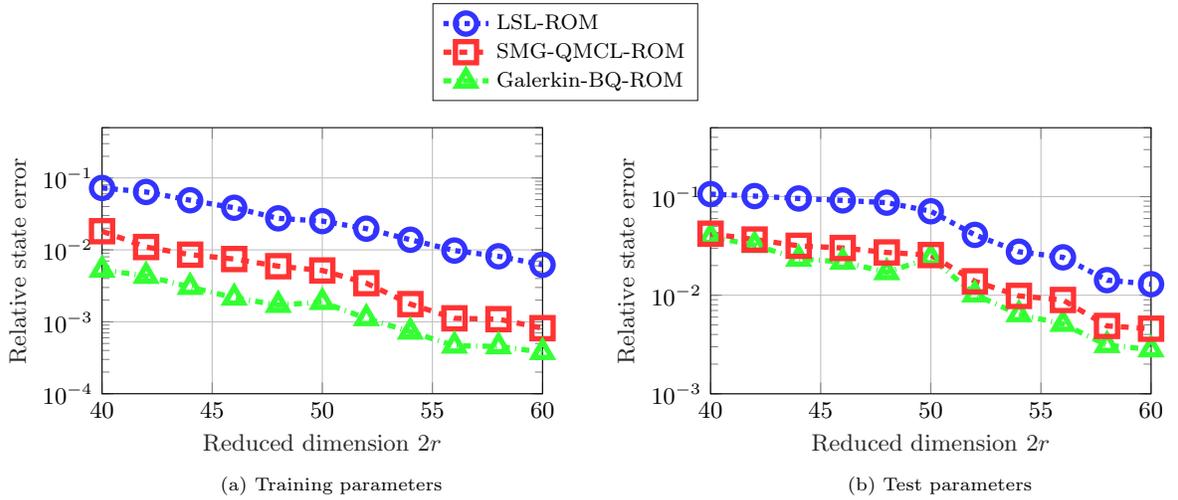
\subsubsection{{\Rb{Parameter extrapolation study}}}

\begin{figure}[tbp]
\small
\captionsetup[subfigure]{oneside,margin={1.8cm,0 cm}}
\begin{subfigure}{.45\textwidth}
       \setlength\fheight{6 cm}
        \setlength\fwidth{\textwidth}
\input{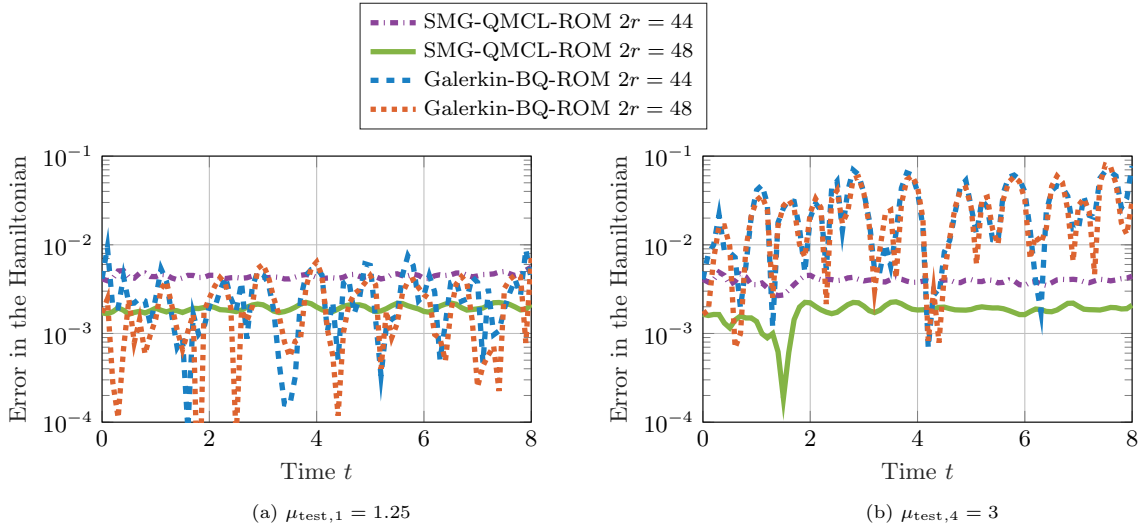}
\caption{$\mu_{\text{test},1}=1.25$}
\label{fig:NW_energy_mu1}
    \end{subfigure}
    \hspace{0.8cm}
    \begin{subfigure}{.45\textwidth}
           \setlength\fheight{6 cm}
           \setlength\fwidth{\textwidth}
\raisebox{-65mm}{
%
%
\definecolor{mycolor1}{rgb}{0.00000,0.44700,0.74100}%
\definecolor{mycolor2}{rgb}{0.85000,0.32500,0.09800}%
\definecolor{mycolor3}{rgb}{0.49400,0.18400,0.55600}%
\definecolor{mycolor4}{rgb}{0.46600,0.67400,0.18800}%
\begin{tikzpicture}

\begin{axis}[%
width=0.951\fheight,
height=0.59\fheight,
at={(0\fheight,0\fheight)},
scale only axis,
xmin=0,
xmax=8,
xlabel style={font=\color{white!15!black}},
xlabel={Time $t$},
ymode=log,
ymin=1e-4,
ymax=0.1,
yminorticks=true,
ylabel style={font=\color{white!15!black}},
ylabel={Error in the Hamiltonian},
axis background/.style={fill=white},
xmajorgrids,
ymajorgrids,
legend style={at={(0.363,0.251)}, anchor=south west, legend cell align=left, align=left, draw=white!15!black}
]

\addplot [color=mycolor3!90, dashdotted, line width=2.0pt]
  table[row sep=crcr]{%
0	0.00407562619962931\\
0.0999999999999996	0.00378848285877353\\
0.199999999999999	0.00405342280931186\\
0.300000000000001	0.00490235696542028\\
0.4	0.00451231145334318\\
0.5	0.00355076723885971\\
0.6	0.00387315475438044\\
0.699999999999999	0.0046713765970523\\
0.800000000000001	0.00434262737266698\\
0.9	0.00356319060999155\\
1	0.00355462461648414\\
1.1	0.00371123438011977\\
1.2	0.00326986090559877\\
1.3	0.00282315727589832\\
1.4	0.00267240954089054\\
1.5	0.00271845740587118\\
1.6	0.0030990764775769\\
1.7	0.00356782993994287\\
1.8	0.00404515331669231\\
1.9	0.00449453789768306\\
2	0.00452776083929863\\
2.1	0.0042083641543691\\
2.2	0.00403349591481472\\
2.3	0.00407741685561191\\
2.4	0.00405563041834966\\
2.5	0.00394646112767915\\
2.6	0.00400191292133023\\
2.7	0.00423259901436879\\
2.8	0.00429801115287453\\
2.9	0.00407638381546071\\
3	0.00391855858639899\\
3.1	0.00400577447388173\\
3.2	0.00401979139816788\\
3.3	0.00383571413244965\\
3.4	0.00380934446812376\\
3.5	0.00405273276089168\\
3.6	0.00420729306137524\\
3.7	0.00409061160426028\\
3.8	0.00388908811414268\\
3.9	0.00377258369669776\\
4	0.00378665252402754\\
4.1	0.00388550744959559\\
4.2	0.00389690415645205\\
4.3	0.00379036610052274\\
4.4	0.00381626181890481\\
4.5	0.0040054991961771\\
4.6	0.00400370534508987\\
4.7	0.00377055462531395\\
4.8	0.00374321588949833\\
4.9	0.00398907399282323\\
5	0.00405025282322846\\
5.1	0.00385873770899003\\
5.2	0.00384956719118179\\
5.3	0.00403292879626705\\
5.4	0.00397197806557692\\
5.5	0.00372599582332889\\
5.6	0.00377159681873706\\
5.7	0.00401371255170701\\
5.8	0.00391565901831293\\
5.9	0.00354619336472694\\
6	0.00348870123947265\\
6.1	0.00377722807818692\\
6.2	0.00386279109302911\\
6.3	0.00360931116929564\\
6.4	0.00348814484576288\\
6.5	0.00373474728475819\\
6.6	0.00401958289550031\\
6.7	0.00407506954331715\\
6.8	0.00405077743521786\\
6.9	0.0040826527024027\\
7	0.0040250655286762\\
7.1	0.00382310601788305\\
7.2	0.003734993090319\\
7.3	0.00391208099063078\\
7.4	0.00408673638534295\\
7.5	0.00401966606224668\\
7.6	0.00392050245846142\\
7.7	0.00401563423783813\\
7.8	0.00414513592854604\\
7.9	0.0042073447110127\\
8	0.00442475745860449\\
};

\addplot [color=mycolor4!90, line width=2.0pt]
  table[row sep=crcr]{%
0	0.0016467615727862\\
0.0999999999999996	0.00159550085968085\\
0.199999999999999	0.00163759730742008\\
0.300000000000001	0.00163446273358403\\
0.4	0.00130910096392078\\
0.5	0.00115575683319591\\
0.6	0.00142051200731252\\
0.699999999999999	0.00157147303514016\\
0.800000000000001	0.00149618941109253\\
0.9	0.00149758849729586\\
1	0.00136560928817353\\
1.1	0.000996392588403337\\
1.2	0.000893403123910729\\
1.3	0.000997672198753464\\
1.4	0.00062087339338035\\
1.5	0.000187058187476552\\
1.6	0.000610191748661392\\
1.7	0.00143364214116737\\
1.8	0.00198438049793559\\
1.9	0.00224167314455137\\
2	0.00221796522330919\\
2.1	0.00201540191615735\\
2.2	0.00184904109604567\\
2.3	0.00173993131894612\\
2.4	0.00167353491633371\\
2.5	0.00173360897634467\\
2.7	0.00211063377512843\\
2.8	0.00223778509782057\\
2.9	0.00221764567328542\\
3	0.00203142767147391\\
3.1	0.00180849637667179\\
3.2	0.0017215824756871\\
3.3	0.00183132467384616\\
3.4	0.00205845248682302\\
3.5	0.00224494911010176\\
3.6	0.0022758242154115\\
3.7	0.00216329613557775\\
3.8	0.00199899895425659\\
3.9	0.00187654953267469\\
4	0.00184064299753573\\
4.1	0.0018501578902832\\
4.2	0.00185663873405595\\
4.3	0.00189497623030199\\
4.4	0.00196984579671344\\
4.5	0.00198139288412413\\
4.6	0.00190624198792833\\
4.7	0.001848472187298\\
4.8	0.00183889232779855\\
4.9	0.00183334165364535\\
5	0.00186566829865953\\
5.1	0.00195288911690561\\
5.2	0.00199667557401639\\
5.3	0.00196278529140592\\
5.4	0.00193976070868196\\
5.5	0.00194249142513492\\
5.6	0.00190176355068772\\
5.7	0.00181752613351272\\
5.9	0.00163422077402764\\
6	0.00163233288626685\\
6.1	0.00176474913049595\\
6.2	0.00188712609652697\\
6.3	0.00186068967447817\\
6.4	0.00177282625259068\\
6.5	0.00175963629652997\\
6.6	0.00183980073949108\\
6.8	0.00218403775624231\\
6.9	0.00219814923914008\\
7	0.00202407403405047\\
7.1	0.00187819734486893\\
7.2	0.00185699965615527\\
7.3	0.00185010758688619\\
7.4	0.00183937997074501\\
7.5	0.00189562943287313\\
7.6	0.00196500570768521\\
7.7	0.00195271107764583\\
7.8	0.00189107466663607\\
7.9	0.00190226197593723\\
8	0.00208233164809179\\
};

\addplot [color=mycolor1!90, dashed, line width=2.0pt]
  table[row sep=crcr]{%
0	0.00407562619962931\\
0.0999999999999996	0.00793222676181494\\
0.199999999999999	0.00996824327603286\\
0.300000000000001	0.0218703792048505\\
0.4	0.00852579962749723\\
0.5	0.00817052763458468\\
0.6	0.00721111791156126\\
0.699999999999999	0.00193005585069272\\
0.800000000000001	0.00573330805516469\\
0.9	0.020678806515166\\
1	0.041258558538475\\
1.1	0.0471225049207131\\
1.2	0.0290412703237414\\
1.3	0.00117673625287118\\
1.4	0.0178555357631662\\
1.5	0.0258308149511031\\
1.6	0.0319607123442047\\
1.7	0.0308796013534938\\
1.8	0.0153585447058937\\
1.9	0.00489919600765809\\
2	0.0200792179322171\\
2.1	0.0293161941061583\\
2.2	0.0243664459247359\\
2.3	0.0052691862547256\\
2.4	0.043545354749586\\
2.5	0.0521469878800236\\
2.6	0.0144649338488225\\
2.7	0.0405710079544086\\
2.8	0.0697901707293039\\
2.9	0.0619434910286732\\
3	0.0386345048610951\\
3.1	0.019017291357704\\
3.2	0.00512181158612002\\
3.3	0.00510368377984882\\
3.4	0.010197955228227\\
3.5	0.00572453500585723\\
3.6	0.0137451427660764\\
3.7	0.0437579306811066\\
3.8	0.064812544260615\\
3.9	0.0613359011047229\\
4	0.0395168599487956\\
4.1	0.0160031544402308\\
4.2	0.000685956404570477\\
4.3	0.00298951578957407\\
4.4	0.00268013694778269\\
4.5	0.00831708789587181\\
4.6	0.0112149897828022\\
4.7	0.0229888172784119\\
4.8	0.0454728233715346\\
4.9	0.0532585134783794\\
5	0.0266177400732488\\
5.1	0.0131485604420511\\
5.2	0.0319533982077101\\
5.3	0.0266071890129797\\
5.4	0.0113729747597393\\
5.5	0.0109929673493116\\
5.6	0.0386330103266346\\
5.7	0.0578103030608155\\
5.8	0.0610545434825109\\
5.9	0.0547241898562962\\
6	0.0396659800817467\\
6.1	0.0145403132028985\\
6.2	0.00398532631113469\\
6.3	0.00175237206519618\\
6.4	0.0235164554274992\\
6.5	0.0407610651453733\\
6.6	0.0462409533310807\\
6.7	0.0441953486204092\\
6.8	0.0353972444887223\\
6.9	0.0150683191298793\\
7	0.0150320144977574\\
7.1	0.0351385719973401\\
7.2	0.0214649585328459\\
7.3	0.0203326363746718\\
7.4	0.056373866516811\\
7.5	0.0674775117517575\\
7.6	0.0663433824172282\\
7.7	0.0597242167358871\\
7.8	0.0260710843594651\\
7.9	0.0333256903149426\\
8	0.0772407714922386\\
};

\addplot [color=mycolor2!90, dotted, line width=2.0pt]
  table[row sep=crcr]{%
0	0.0016467615727862\\
0.0999999999999996	0.00182196549627101\\
0.199999999999999	0.00689106027309627\\
0.300000000000001	0.0162683920433757\\
0.4	0.016328255394427\\
0.5	0.00958194866368344\\
0.6	0.000691041553682937\\
0.699999999999999	0.000993889646863977\\
0.800000000000001	0.00901344041134422\\
0.9	0.0210165028604478\\
1	0.0317263989344152\\
1.1	0.0357726287010393\\
1.2	0.0248297700784264\\
1.3	0.00253933080017709\\
1.4	0.0162924879654209\\
1.5	0.0275055070337285\\
1.6	0.0320170798655966\\
1.7	0.0257879378523445\\
1.8	0.0101496544867387\\
1.9	0.00724352227294566\\
2	0.0228780872997145\\
2.1	0.0331321427002571\\
2.2	0.0261775910471744\\
2.3	0.00214473577648455\\
2.4	0.0331451964411666\\
2.5	0.0411106436647107\\
2.6	0.0173714207073665\\
2.7	0.0237295296145934\\
2.8	0.0552461151531355\\
2.9	0.0580725562250653\\
3	0.0380420249683784\\
3.1	0.0164551905968722\\
3.2	0.00306393868541122\\
3.3	0.00907331599435323\\
3.4	0.0225033418651961\\
3.5	0.024506979687045\\
3.6	0.00423606278696554\\
3.7	0.0288868537801062\\
3.8	0.0542925352623645\\
3.9	0.0597863545056754\\
4	0.0450268523832965\\
4.1	0.0197321784434212\\
4.2	0.000826293011090938\\
4.3	0.00302087307171916\\
4.4	0.000790862804759486\\
4.5	0.00461497435007808\\
4.6	0.0119187771640412\\
4.7	0.0263856628359154\\
4.8	0.0412487373840982\\
4.9	0.0440766893586416\\
5	0.0265781915811836\\
5.1	0.00604468140738489\\
5.2	0.032693049278608\\
5.3	0.034485319634678\\
5.4	0.0148480550776\\
5.5	0.0101401605316358\\
5.6	0.0326594863790497\\
5.7	0.0508258205038138\\
5.8	0.0586231892249013\\
5.9	0.0516046473615183\\
6	0.0323453897725461\\
6.1	0.00699286544792343\\
6.2	0.0125680093572829\\
6.3	0.0104658924187415\\
6.4	0.0153077984150407\\
6.5	0.0456400207120229\\
6.6	0.0595634713288325\\
6.7	0.050835338672707\\
6.8	0.0246921546248906\\
6.9	0.00654156314664965\\
7	0.0277892687788909\\
7.1	0.0332960777595215\\
7.2	0.0239405589636536\\
7.3	0.0046671898129027\\
7.4	0.0482731776418817\\
7.5	0.0792318228721438\\
7.6	0.0752497927635078\\
7.7	0.0485257528578025\\
7.8	0.0176885040047168\\
7.9	0.0131154396580761\\
8	0.0310344525313067\\
};

\end{axis}
\end{tikzpicture}%
}
\caption{$\mu_{\text{test},4}=3$}
\label{fig:NW_energy_mu2}
    \end{subfigure}
\caption{{\Rb{Two-dimensional nonlinear wave equation. The error in the Hamiltonian~\eqref{eq:err_hamiltonian_non} for SMG-QMCL-ROMs and Galerkin-BQ-ROMs remains below $10^{-2}$ for $\mu_{\text{test},1}=0.125$  whereas both approaches exhibit bounded error in the Hamiltonian below $10^{-1}$ for $\mu_{\text{test},4}=3$.}}}
 \label{fig:NW_energy}
\end{figure}

\begin{figure}[tbp]
\tiny{LSL-ROM}
\hspace{1.2cm}
\begin{subfigure}{.27\textwidth}
\includegraphics[width=\linewidth]{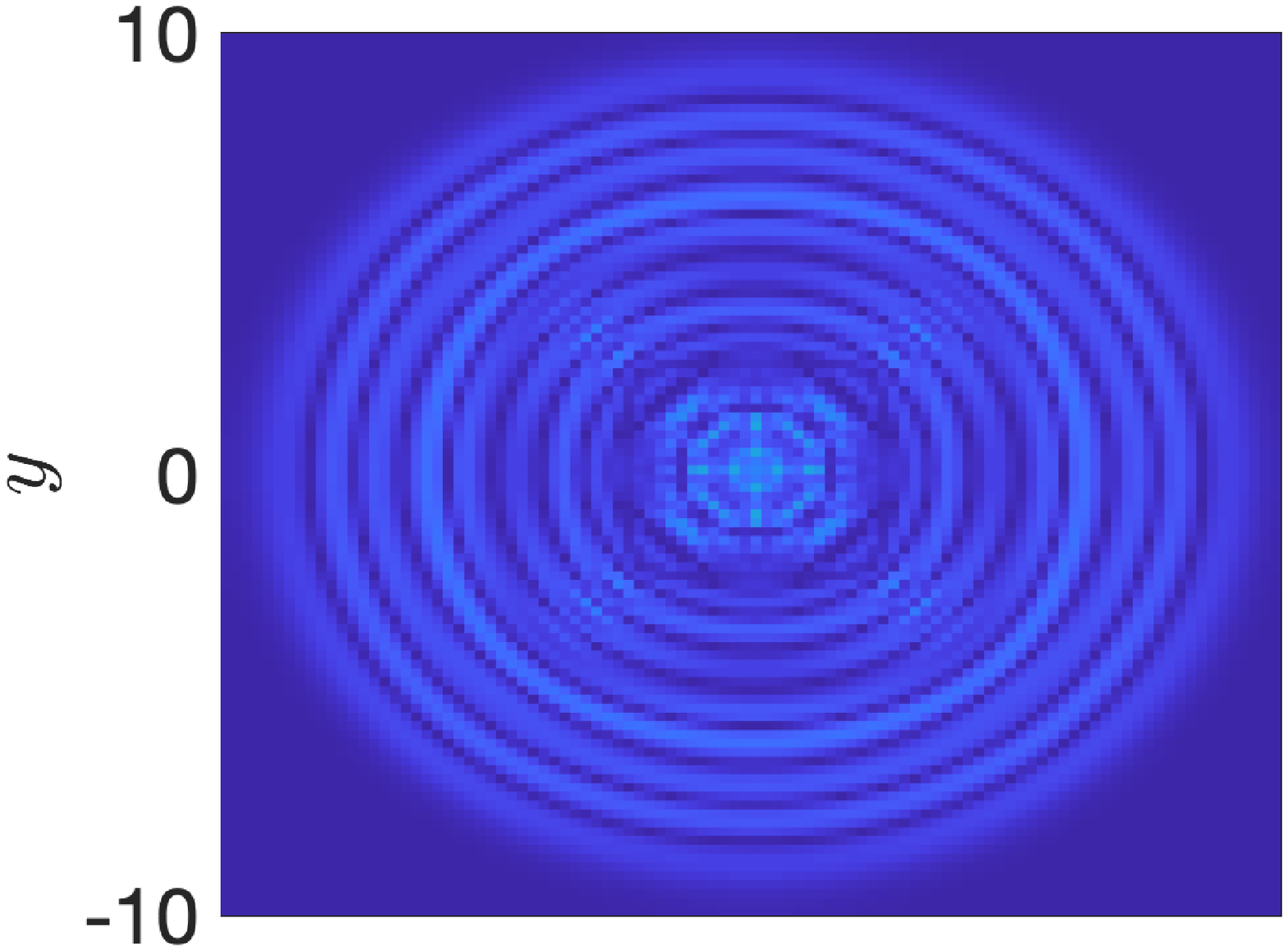}
\end{subfigure}
\begin{subfigure}{.27\textwidth}
\includegraphics[width=\linewidth]{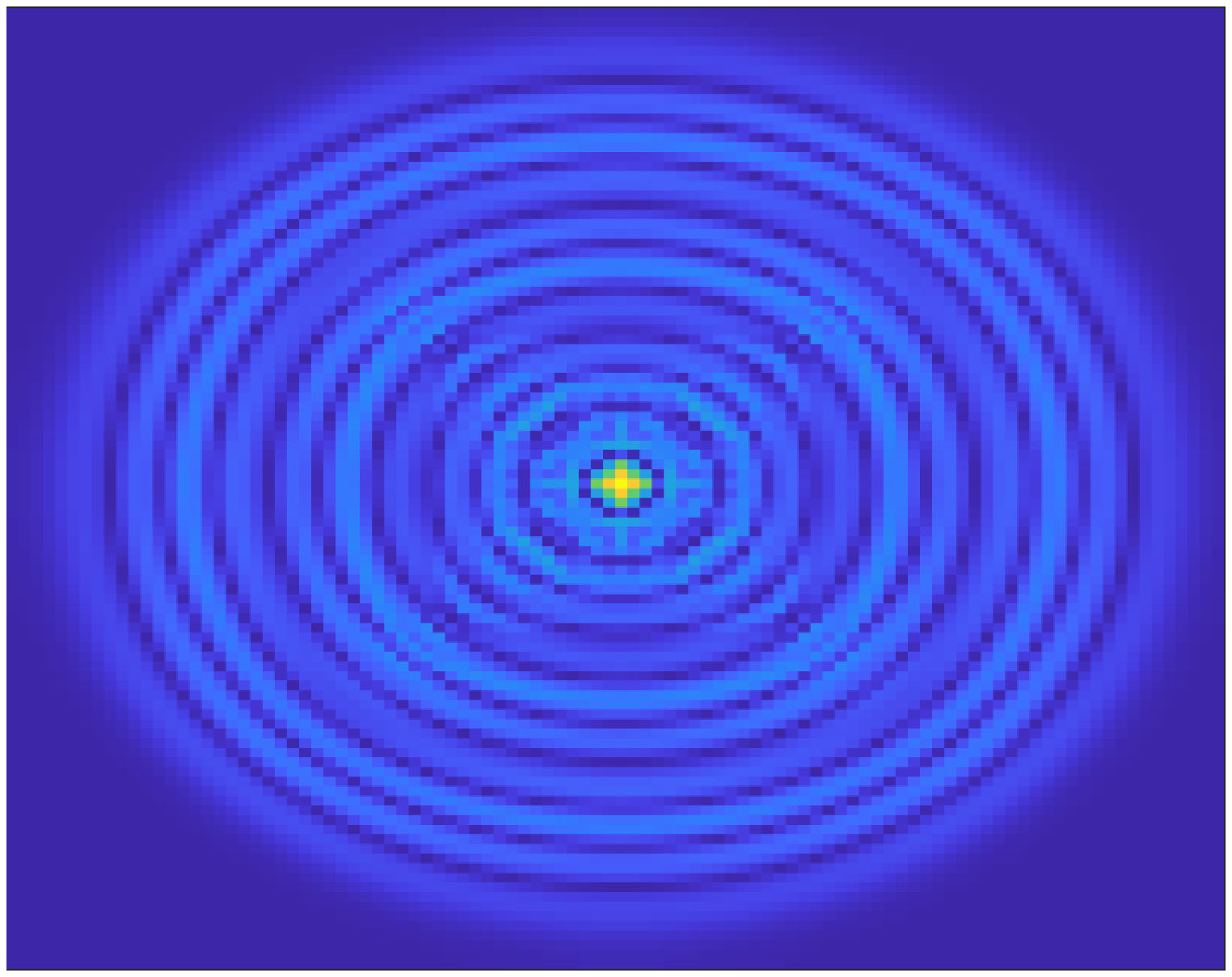}
\end{subfigure}
 \begin{subfigure}{.27\textwidth}
\includegraphics[width=\linewidth]{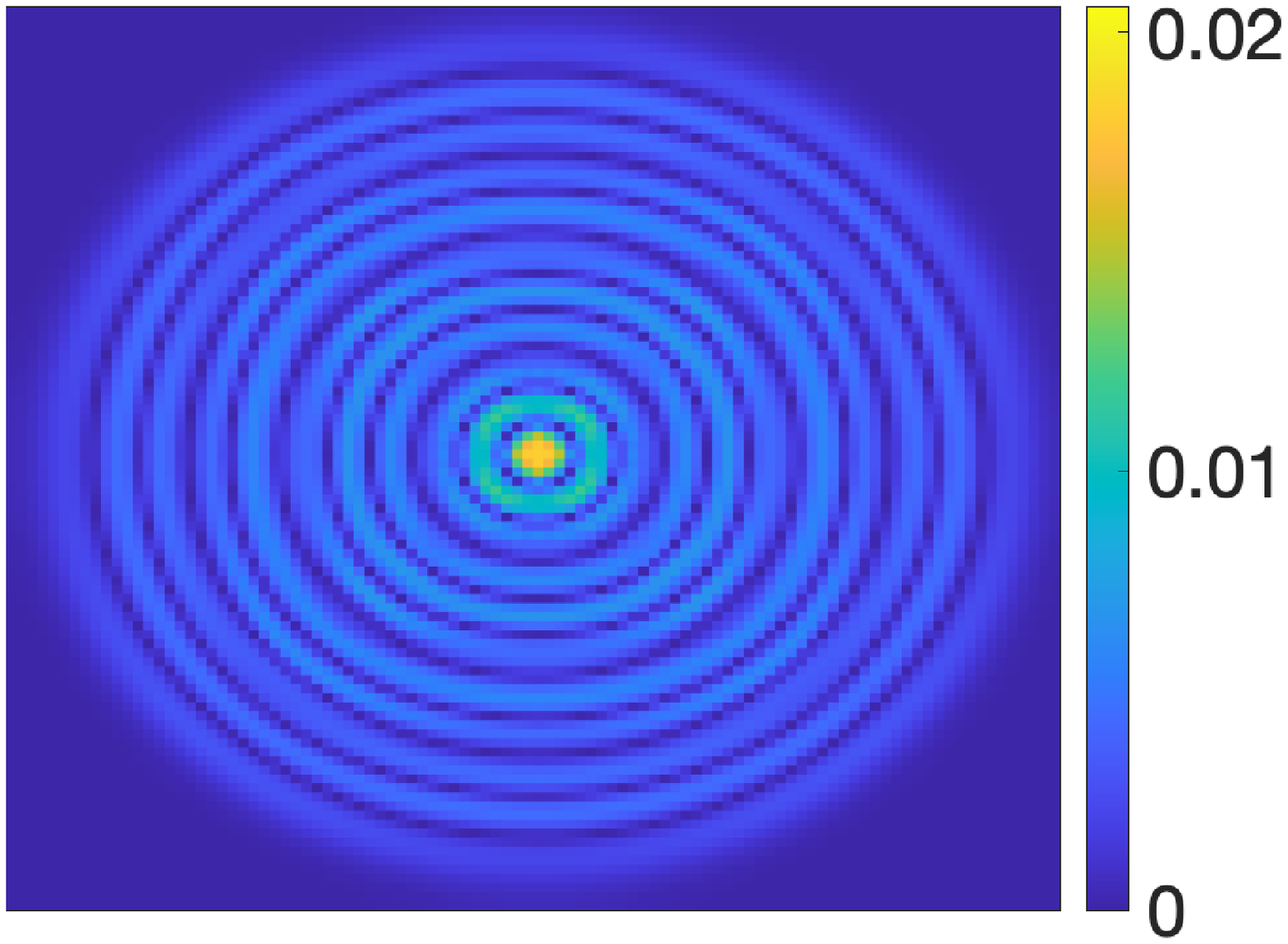}
 \end{subfigure}\\
 \tiny{SMG-QMCL-ROM}
 \hspace{0.3cm}
 \begin{subfigure}{.27\textwidth}
\includegraphics[width=\linewidth]{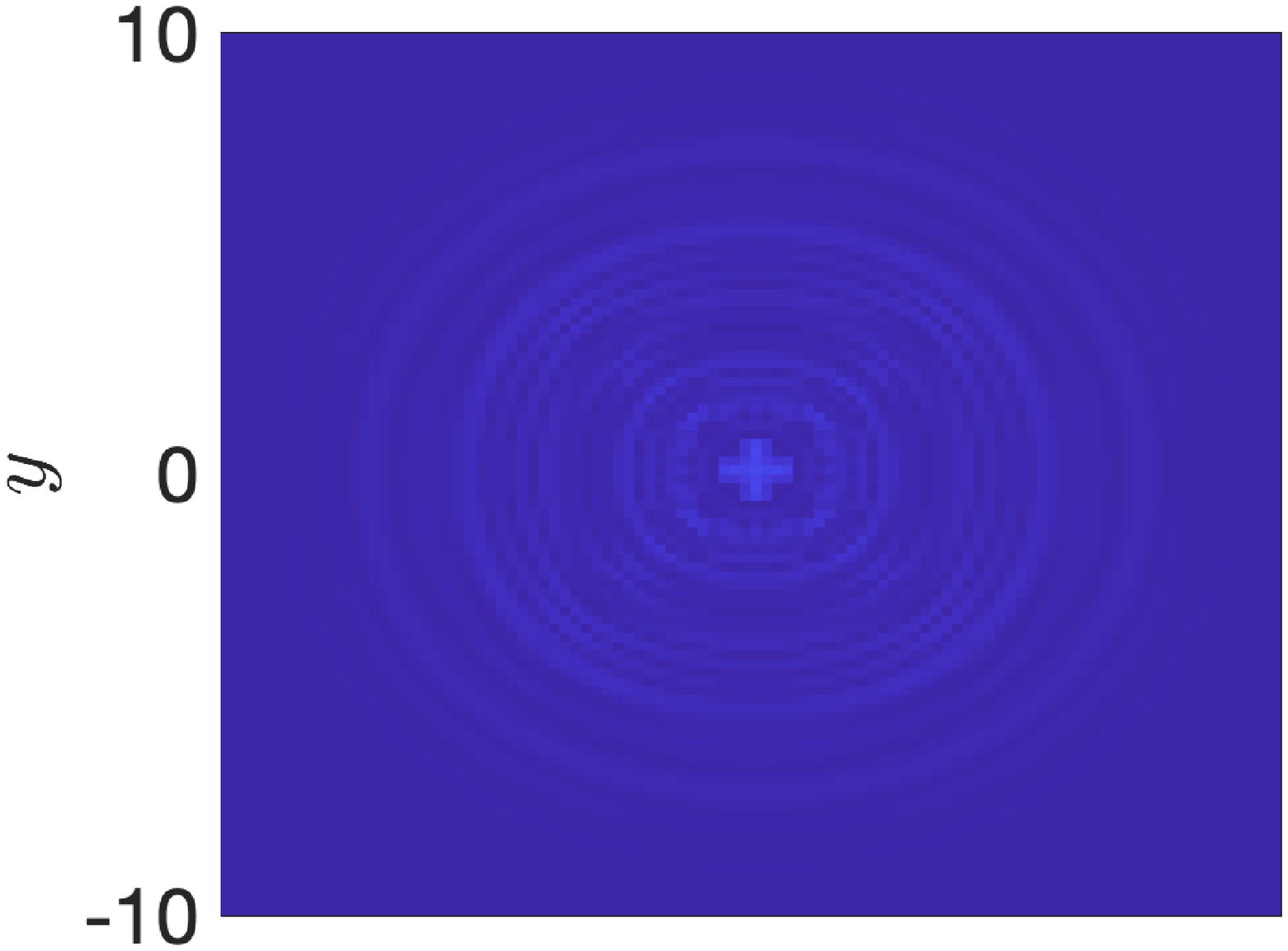}
\end{subfigure} 
\begin{subfigure}{.27\textwidth}
\includegraphics[width=\linewidth]{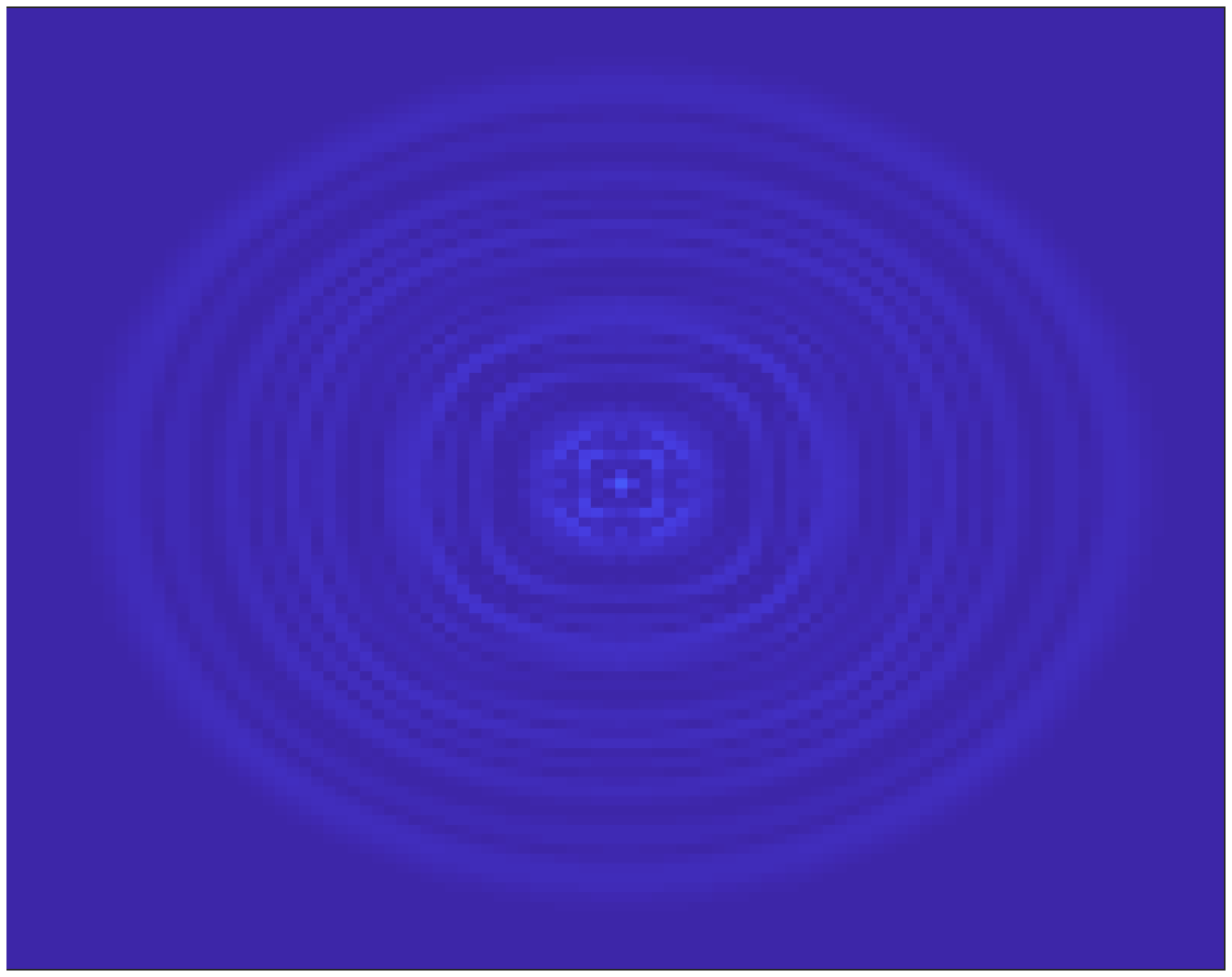}
\end{subfigure}
 \begin{subfigure}{.27\textwidth}
\includegraphics[width=\linewidth]{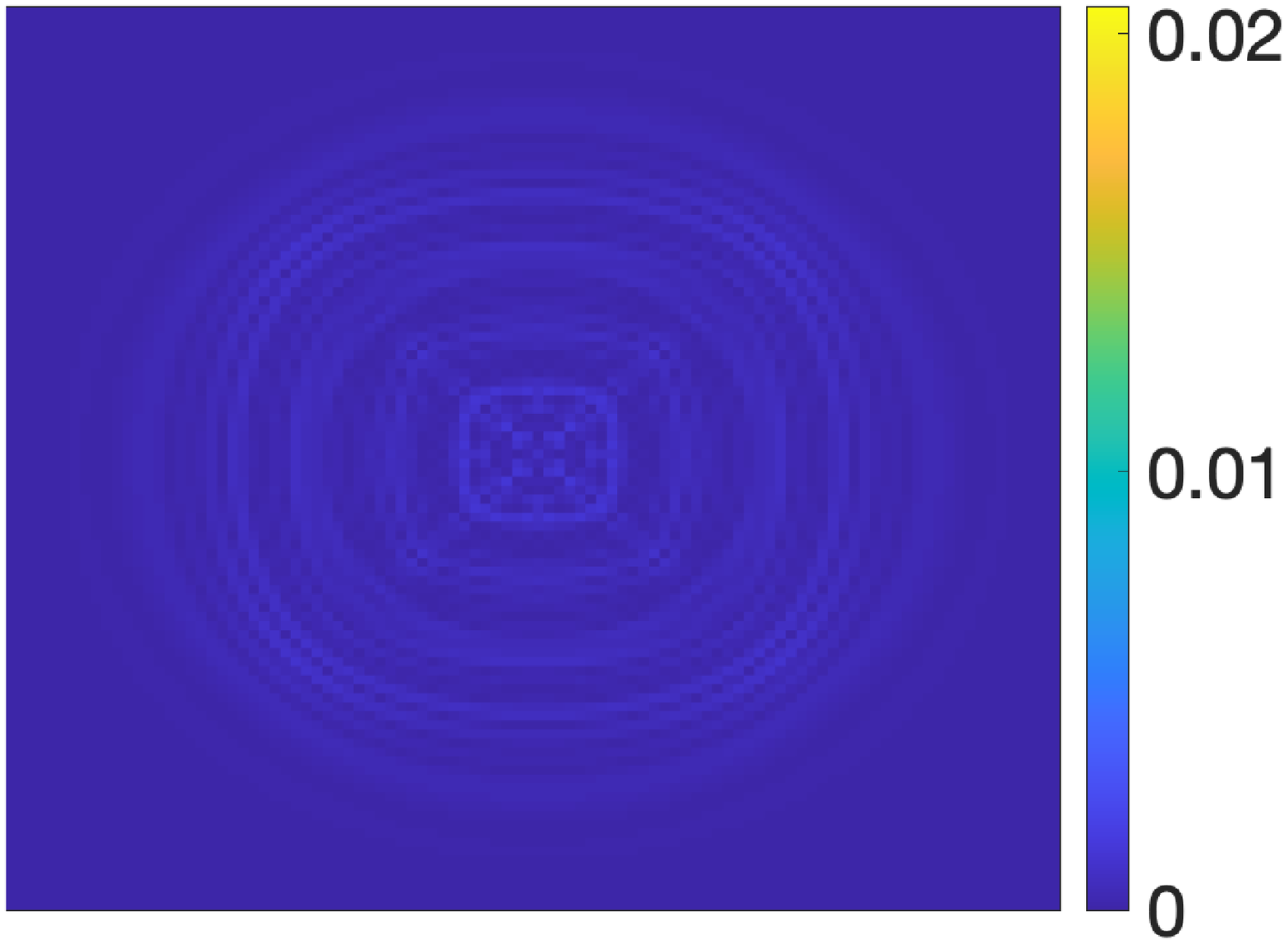}
 \end{subfigure}\\
 \tiny{Galerkin-BQ-ROM}
 \hspace{0.2cm}
 \begin{subfigure}{.27\textwidth}
\includegraphics[width=1.025\linewidth,height=0.8\linewidth]{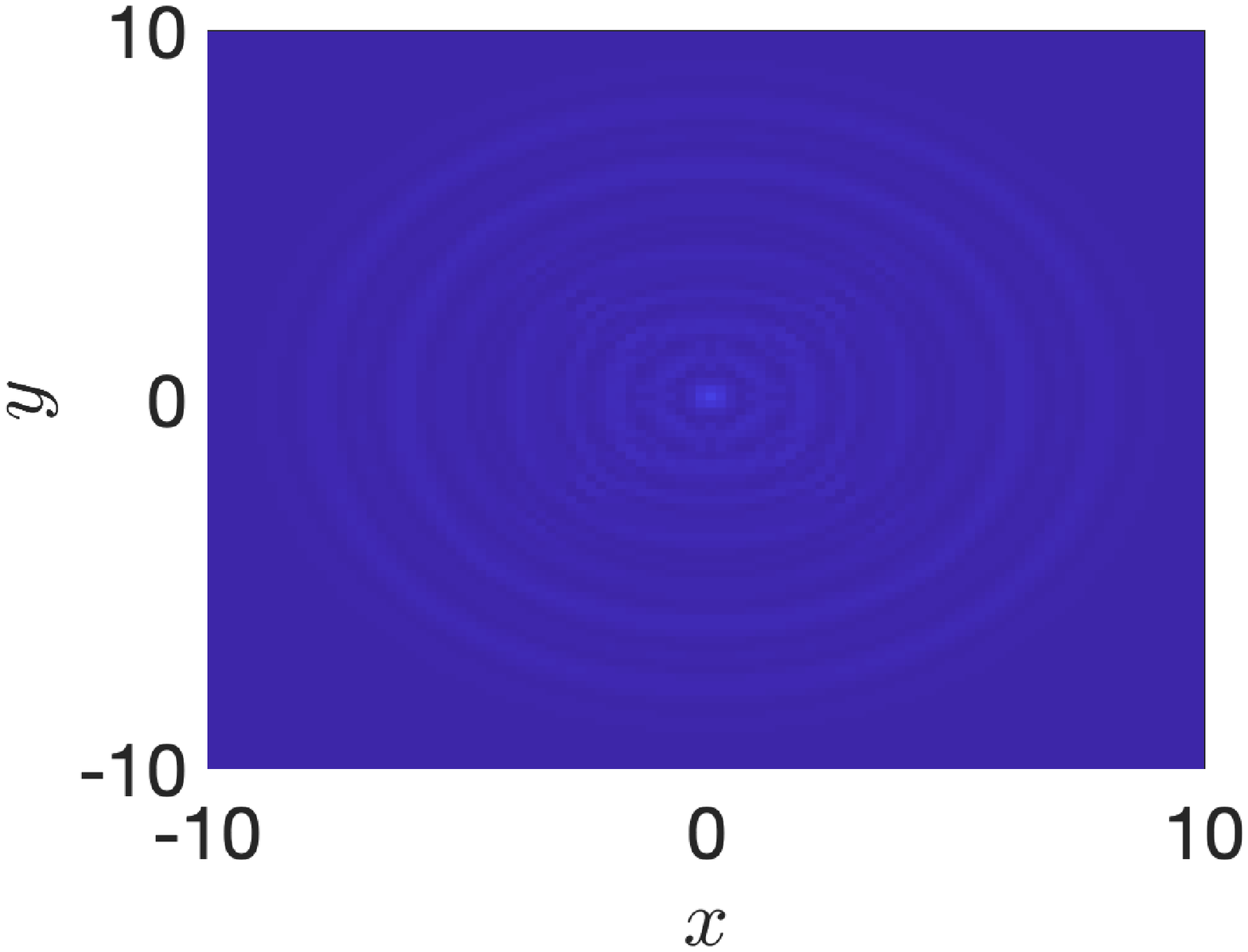}
\subcaption{$t=6$}
\end{subfigure}
\begin{subfigure}{.27\textwidth}
\includegraphics[width=1.025\linewidth,height=0.8\linewidth]{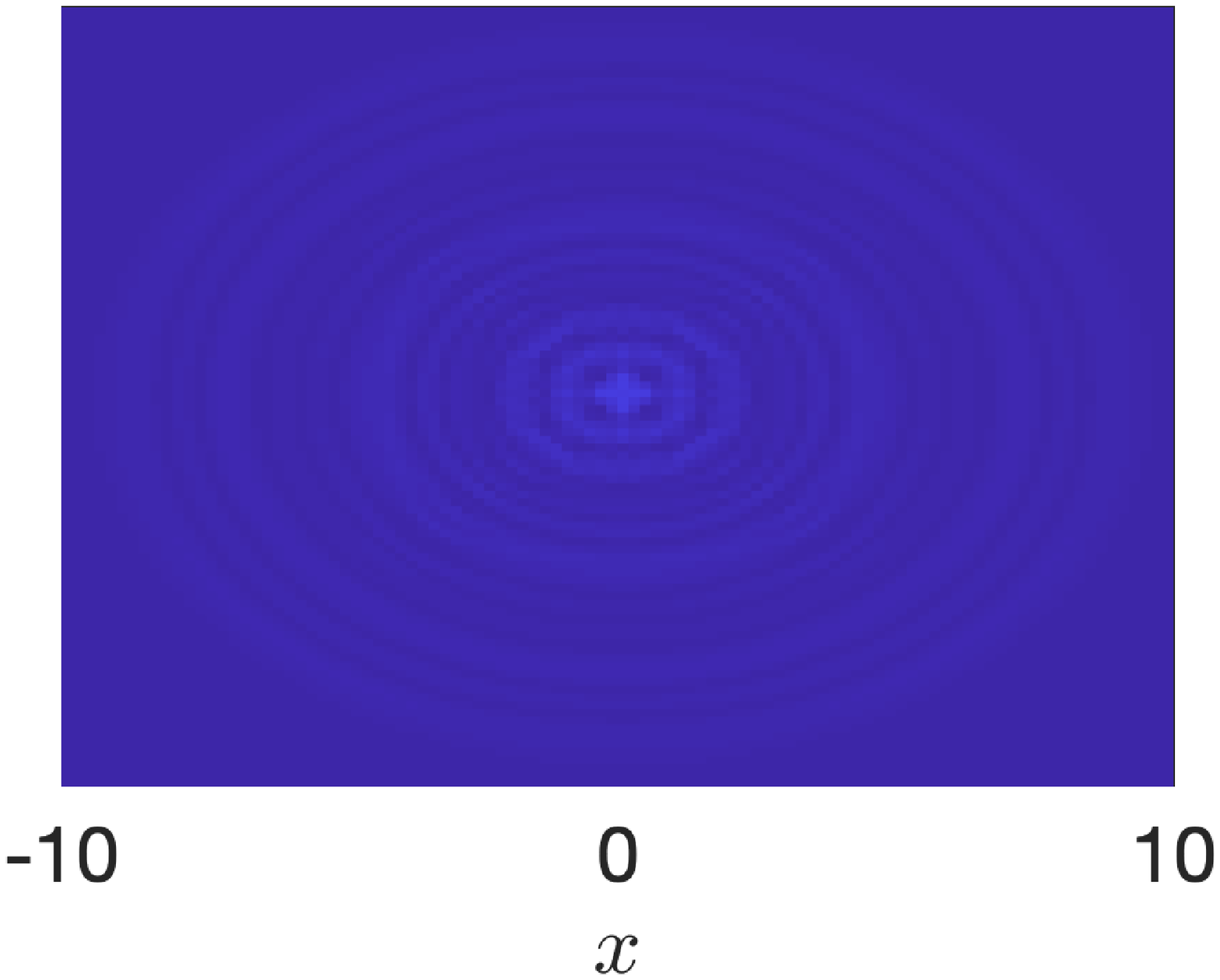}
\subcaption{$t=7$}
\end{subfigure} \  
 \begin{subfigure}{.27\textwidth}
\includegraphics[width=1.025\linewidth,height=0.8\linewidth]{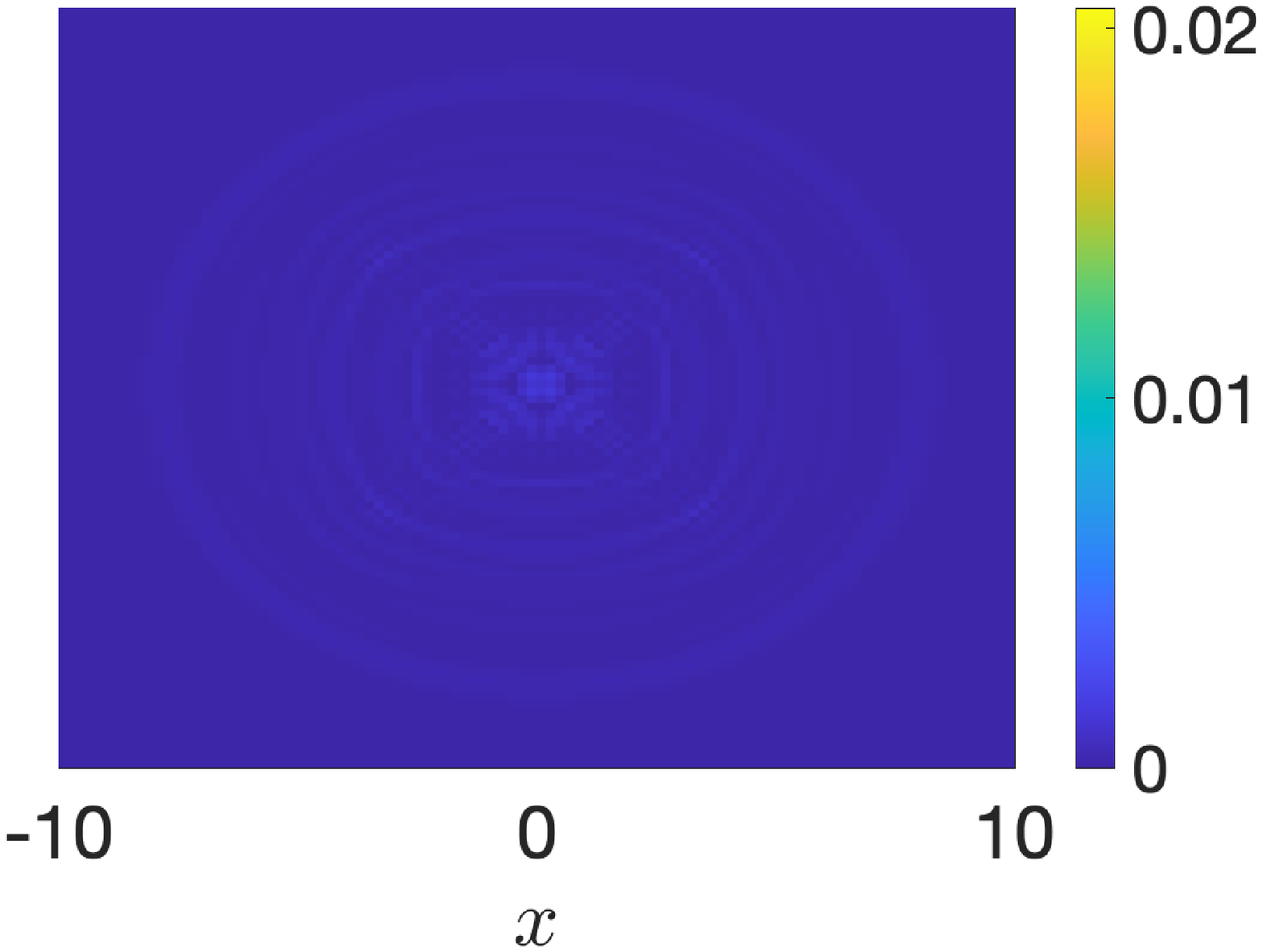}
\subcaption{$t=8$}
 \end{subfigure} 
    \caption{{\Rb{Two-dimensional nonlinear wave equation. Plots compare the pointwise error in $\bq$ between the FOM solution and the reconstructed solution $\bGamma_{\bq}(\widetilde{\bQ}(\mu))$ of the ROM solution $\widetilde{\bQ}(\mu)$ using low-dimensional ($2r=48$) LSL-ROM, SMG-QMCL-ROM, and Galerkin-BQ-ROM at selected time instances $t\in\{ 6,7,8\}$ for $\mu_{\text{test,1}}=1.25$. The pointwise error for the linear symplectic ROM increases as we march forward in time whereas the SMG-QMCL ROM and Galerkin-BQ-ROM yield accurate approximate solutions even at $t=8$.}}}
\label{fig:NW_gam_11}
\end{figure}

\begin{figure}[tbp]
\tiny{LSL-ROM}
\hspace{1.2cm}
\begin{subfigure}{.27\textwidth}
\includegraphics[width=\linewidth]{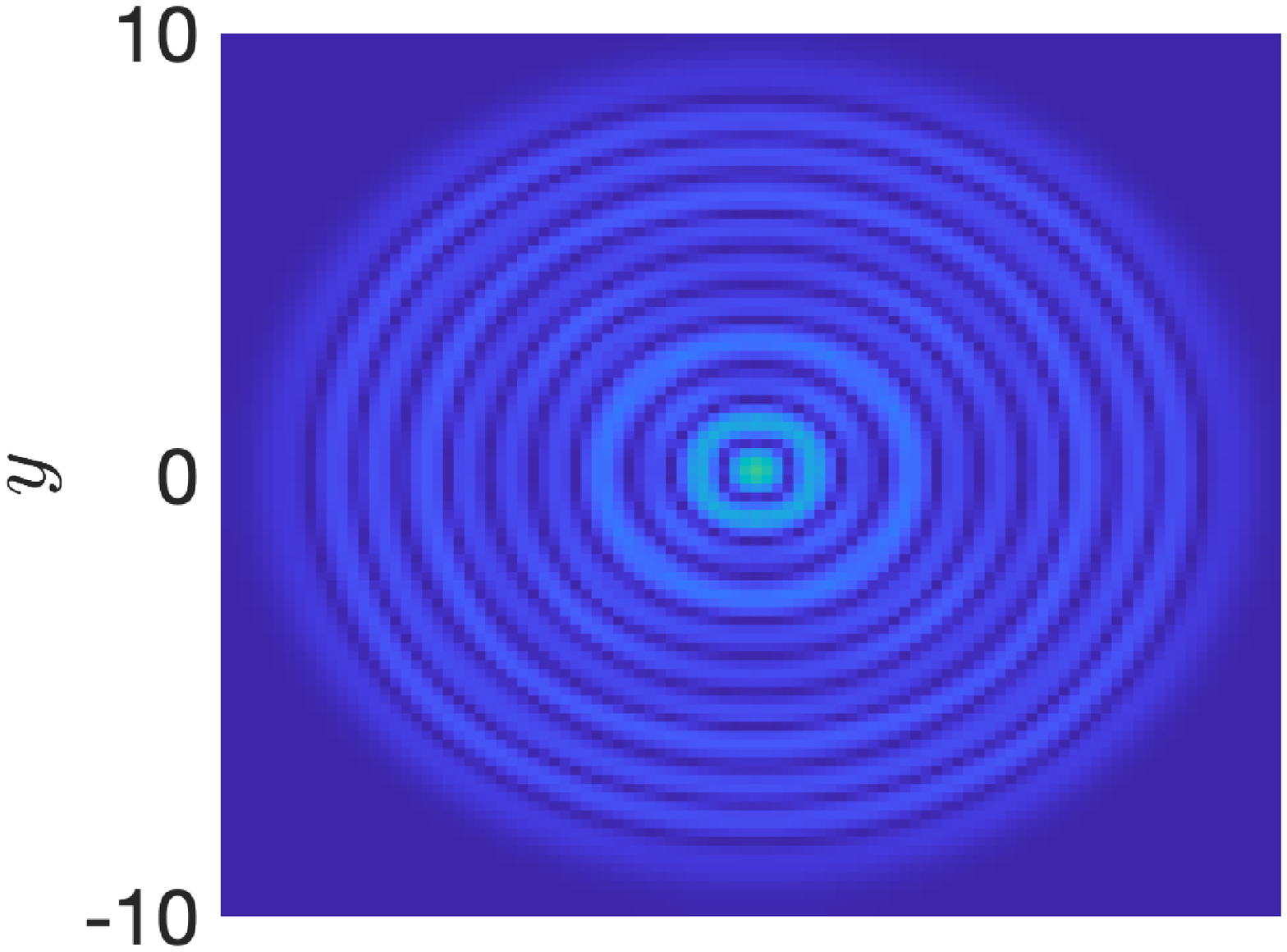}
\end{subfigure}
\begin{subfigure}{.27\textwidth}
\includegraphics[width=\linewidth]{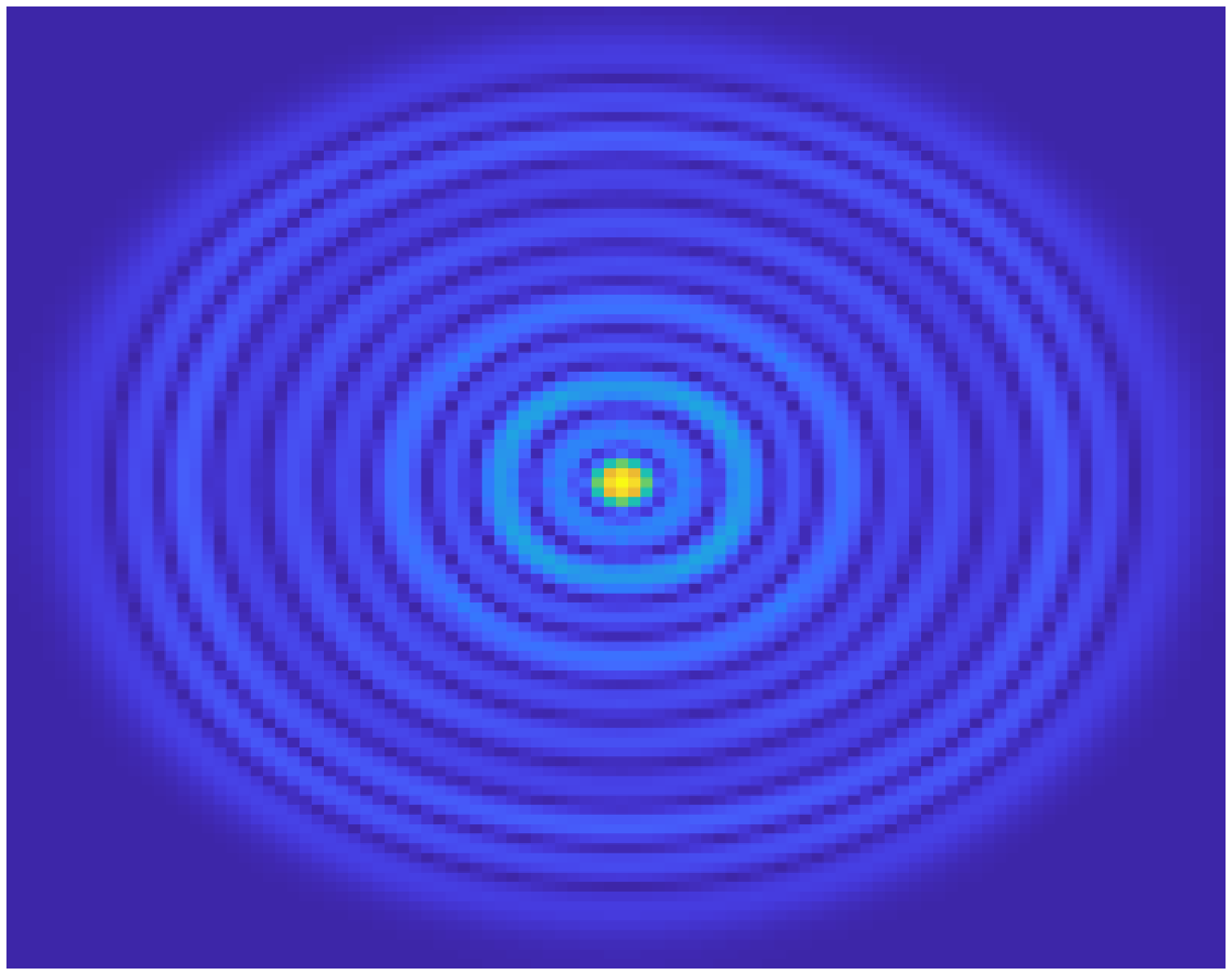}
\end{subfigure}
 \begin{subfigure}{.27\textwidth}
\includegraphics[width=\linewidth]{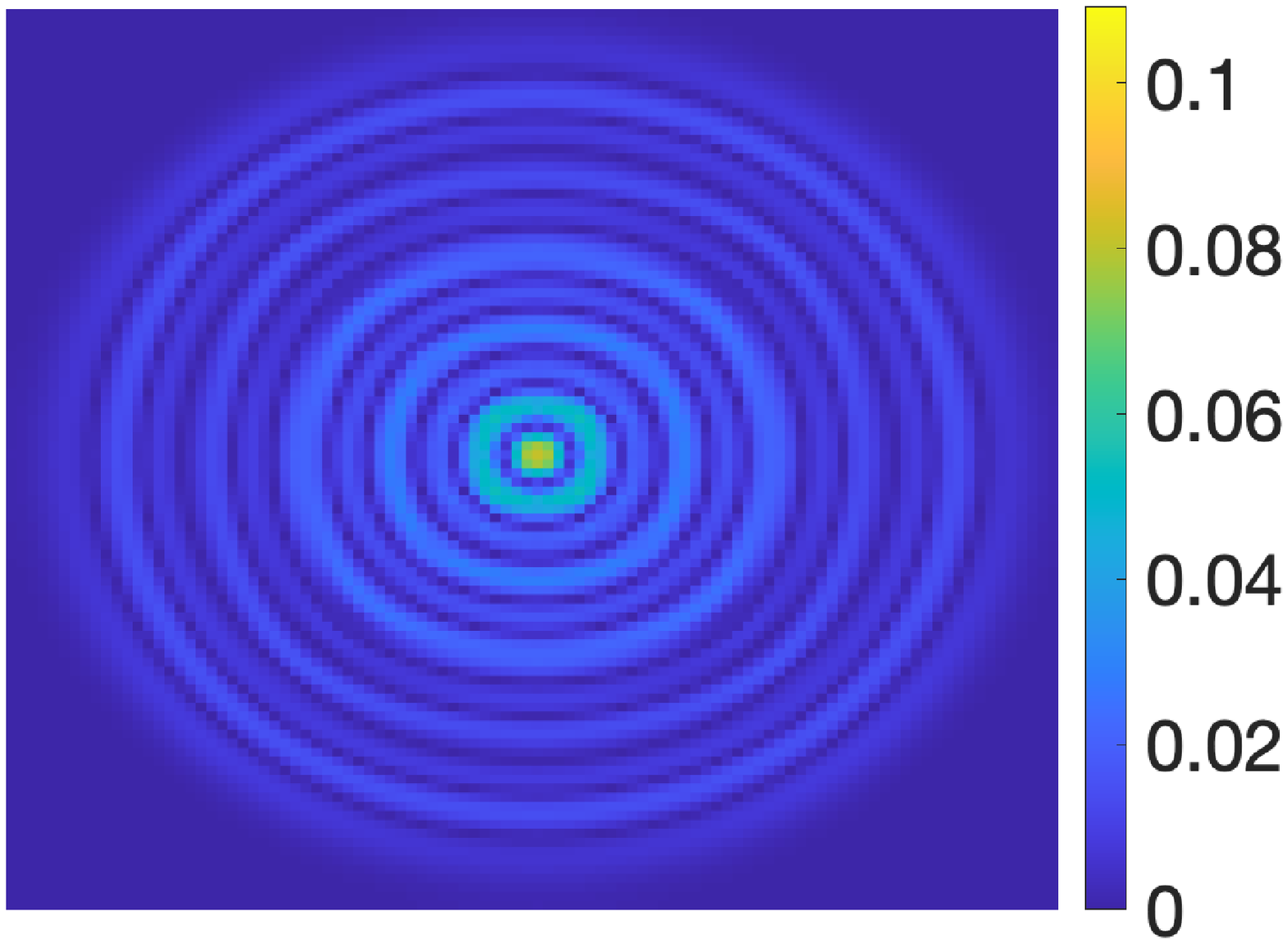}
 \end{subfigure}\\
 \tiny{SMG-QMCL-ROM}
 \hspace{0.3cm}
 \begin{subfigure}{.27\textwidth}
\includegraphics[width=\linewidth]{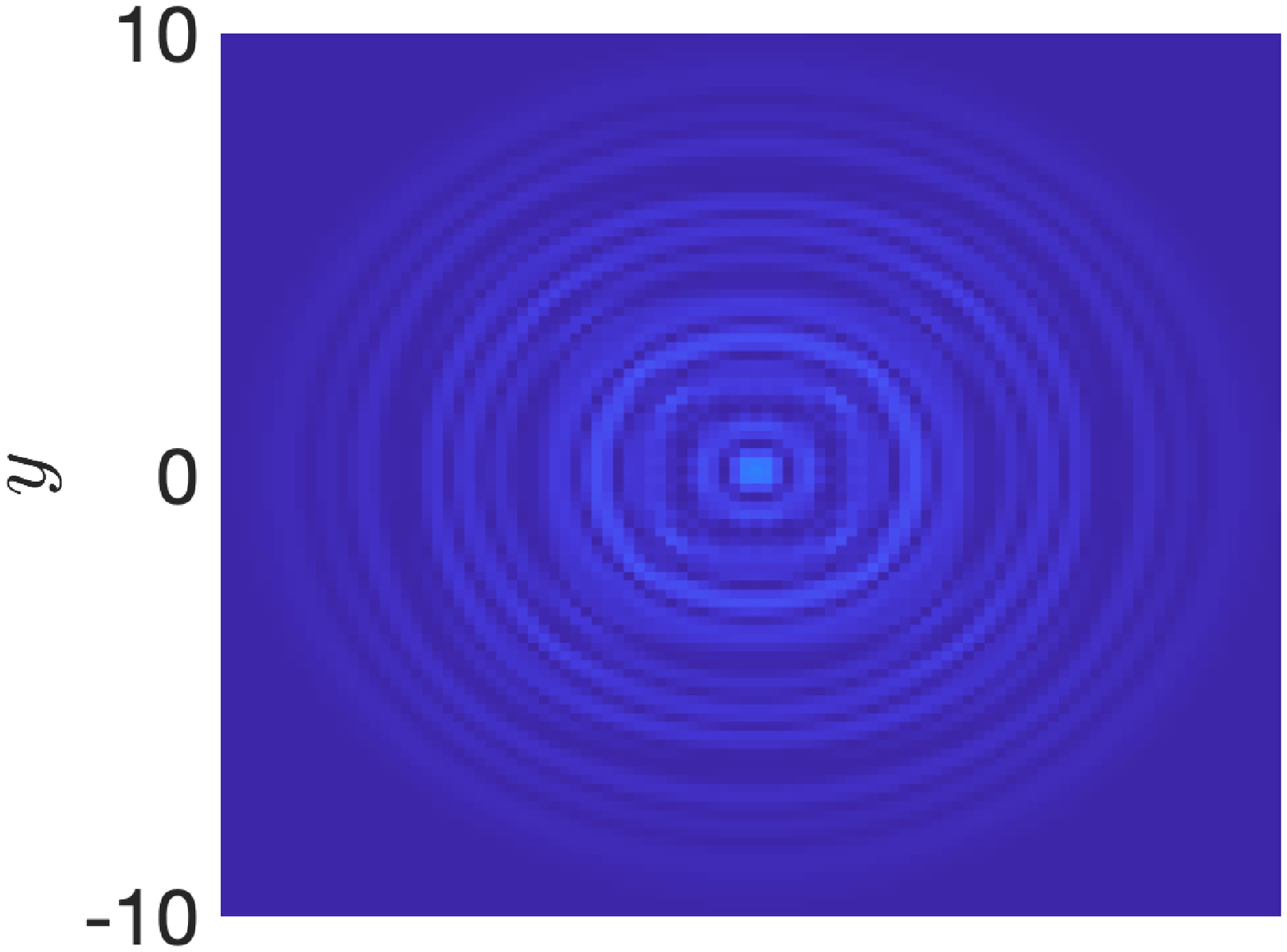}
\end{subfigure} 
\begin{subfigure}{.27\textwidth}
\includegraphics[width=\linewidth]{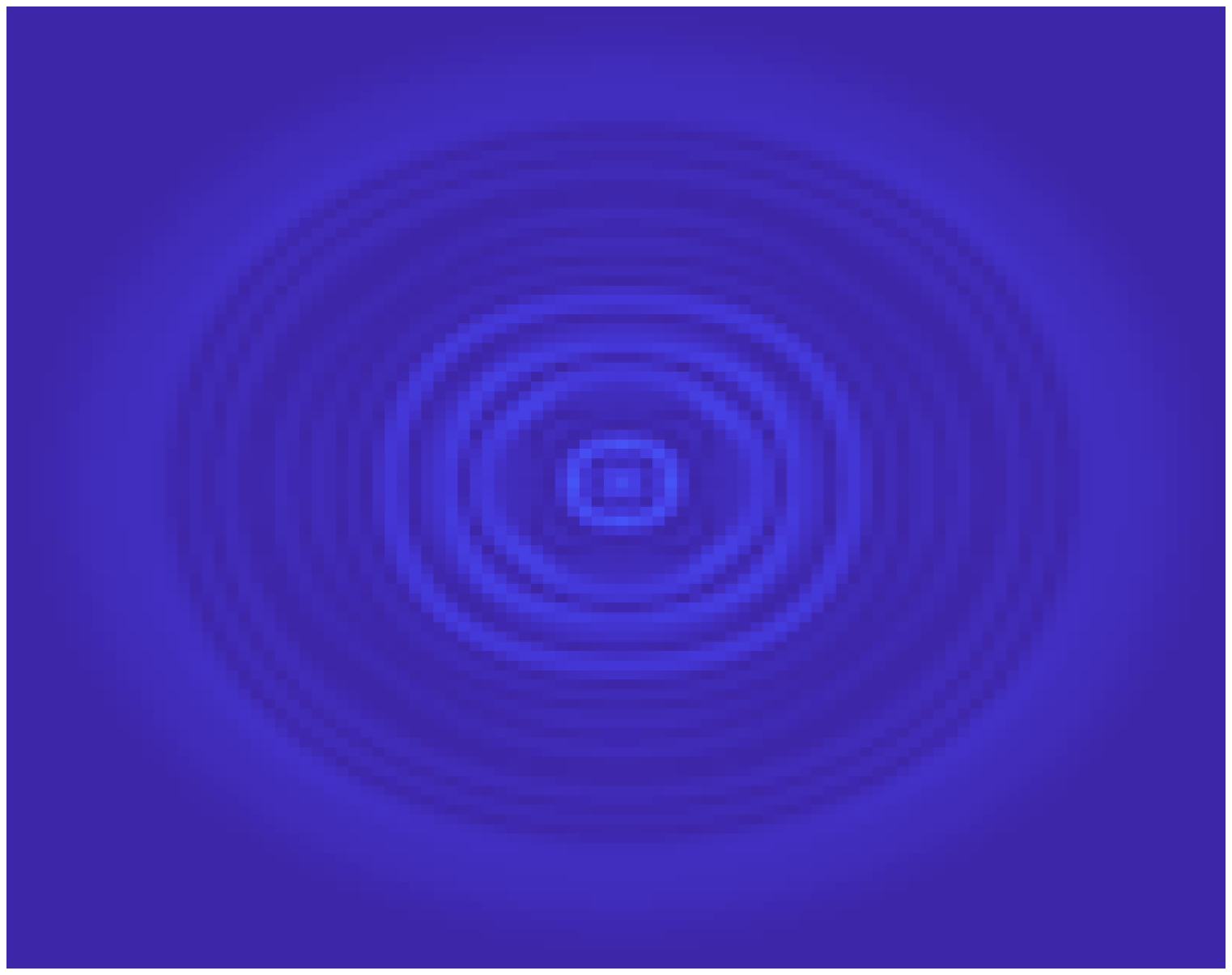}
\end{subfigure}
 \begin{subfigure}{.27\textwidth}
\includegraphics[width=\linewidth]{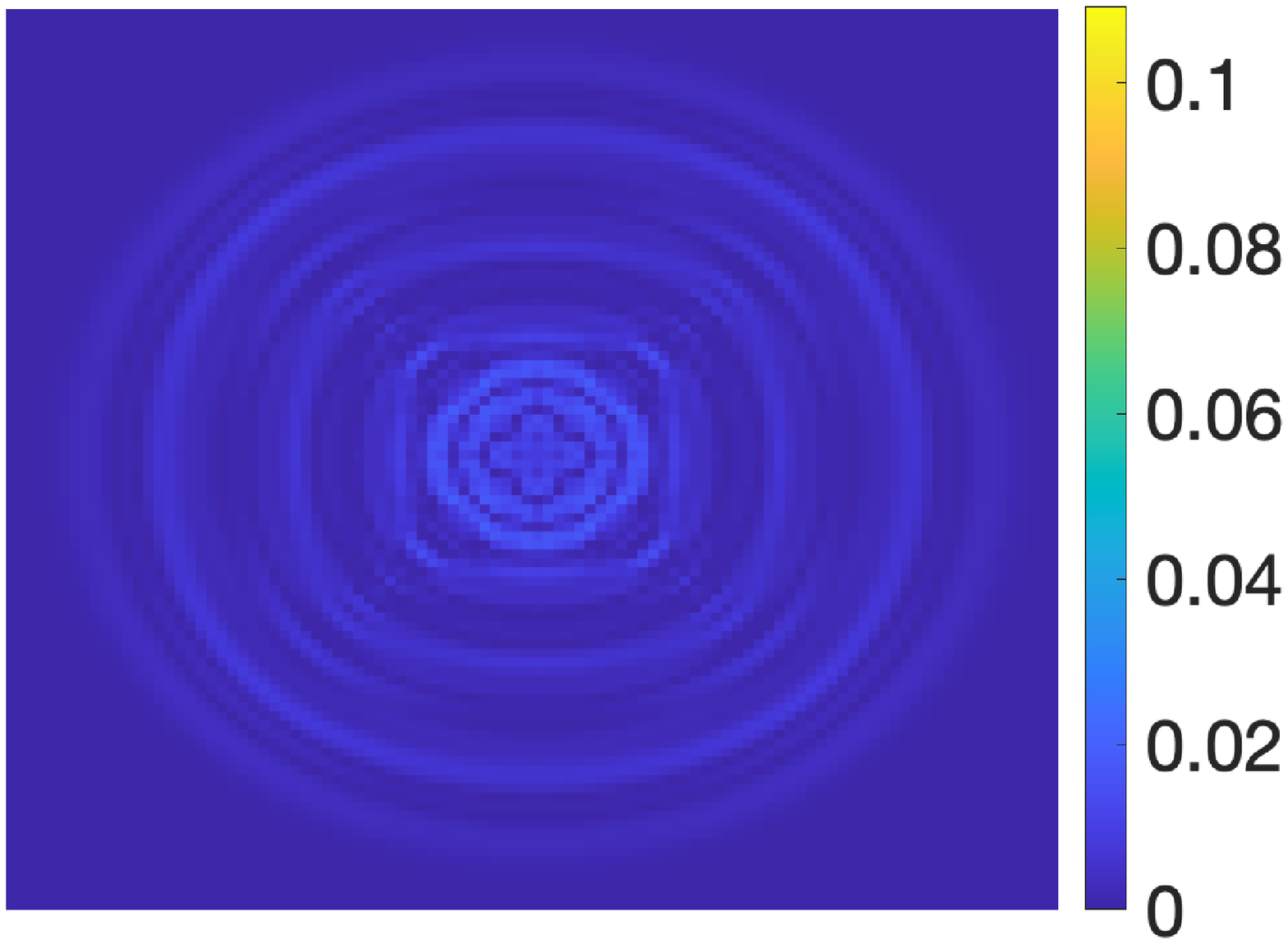}
 \end{subfigure}\\
 \tiny{Galerkin-BQ-ROM}
 \hspace{0.2cm}
 \begin{subfigure}{.27\textwidth}
\includegraphics[width=1.025\linewidth,height=0.8\linewidth]{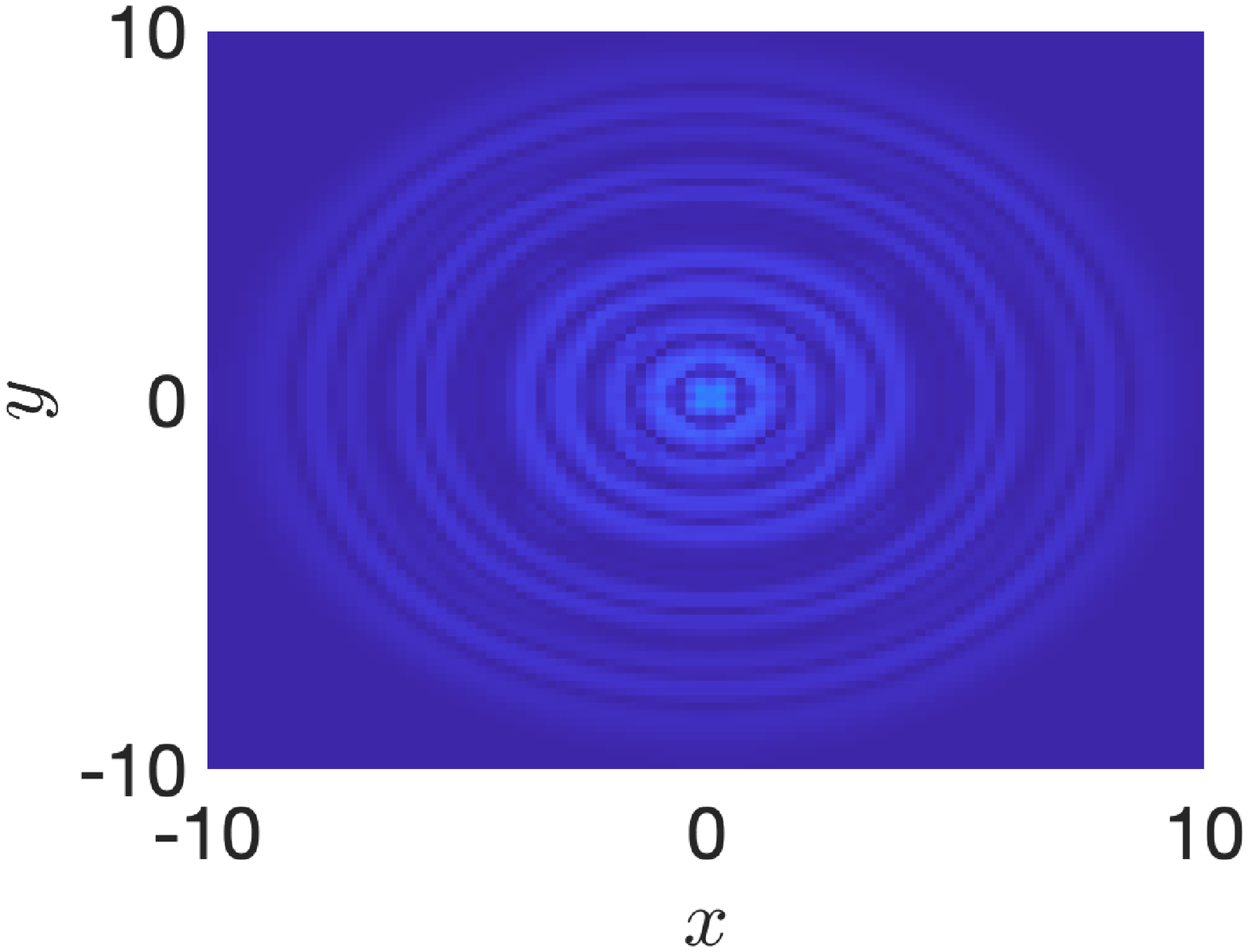}
\subcaption{$t=6$}
\end{subfigure}
\begin{subfigure}{.27\textwidth}
\includegraphics[width=1.025\linewidth,height=0.8\linewidth]{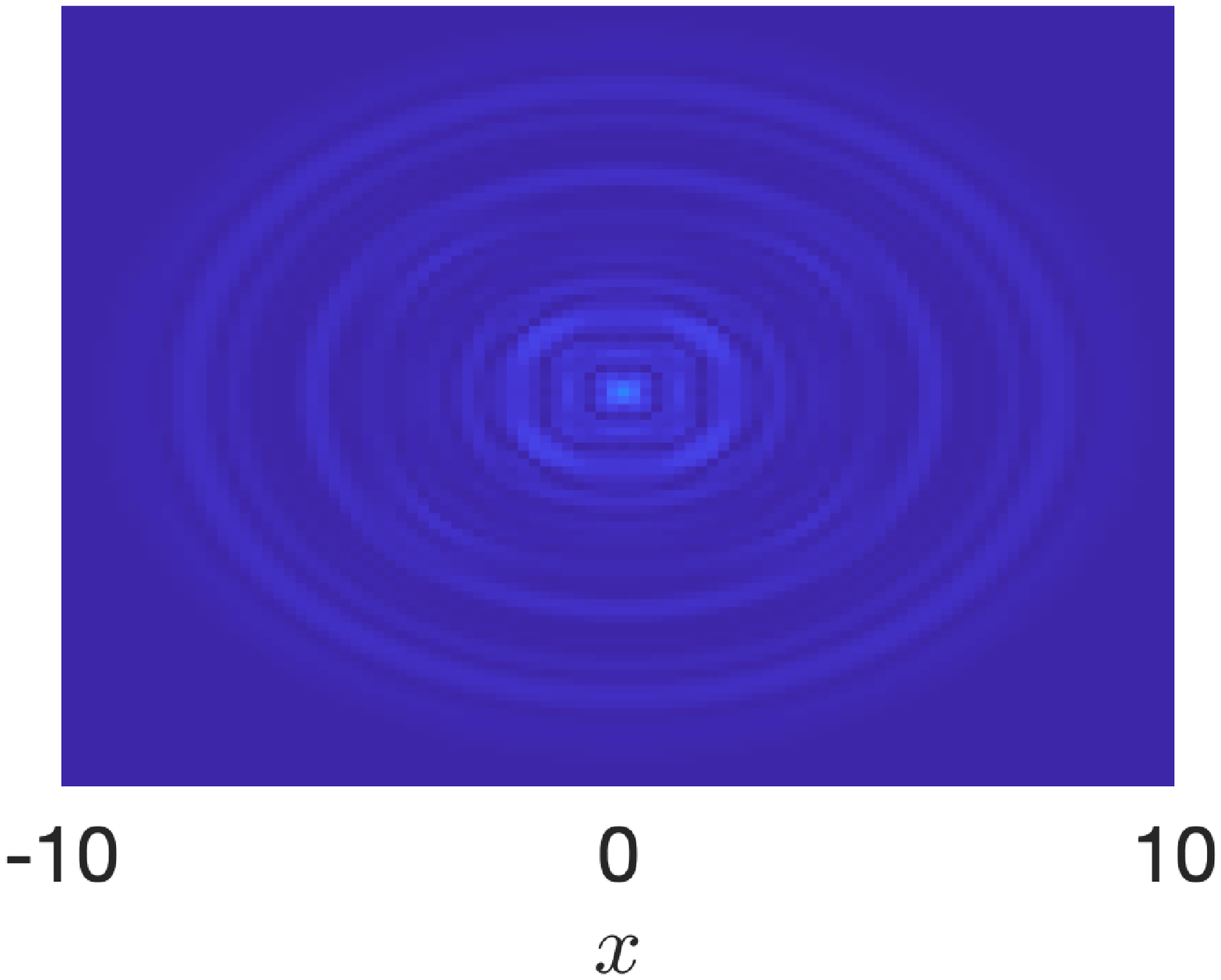}
\subcaption{$t=7$}
\end{subfigure} \  
 \begin{subfigure}{.27\textwidth}
\includegraphics[width=1.025\linewidth,height=0.8\linewidth]{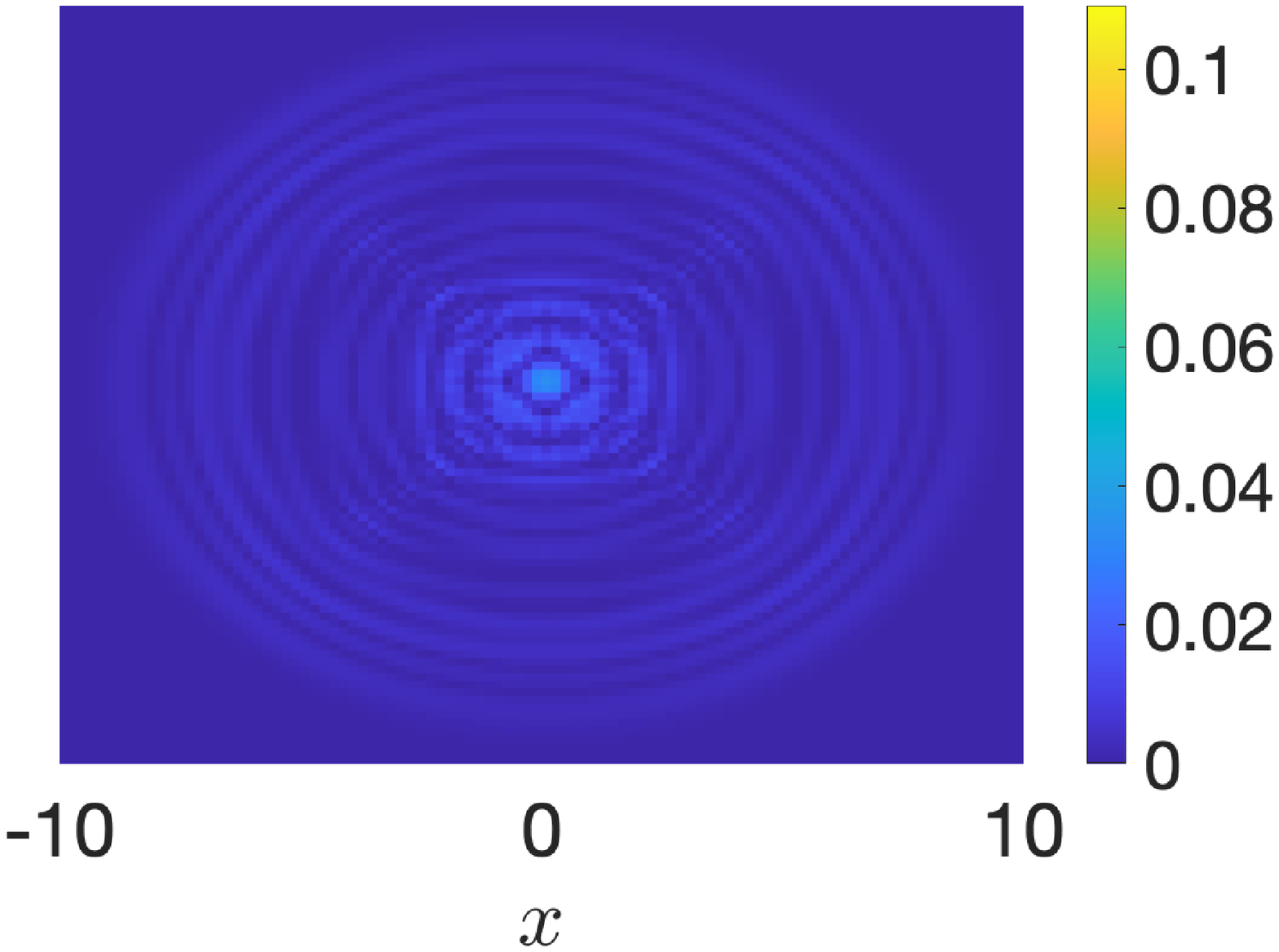}
\subcaption{$t=8$}
 \end{subfigure}
    \caption{{\Rb{Two-dimensional nonlinear wave equation. Plots compare the pointwise error in $\bq$ between the FOM solution and the reconstructed solution $\bGamma_{\bq}(\widetilde{\bQ}(\mu))$ of the ROM solution $\widetilde{\bQ}(\mu)$ using low-dimensional ($2r=48$) LSL-ROM, SMG-QMCL-ROM, and Galerkin-BQ-ROM at selected time instances $t\in\{ 6,7,8\}$ for $\mu_{\text{test,4}}=3$. The proposed approaches yield a lower pointwise error than the LSL-ROM at all three time instances.}} }
\label{fig:NW_gam_14}
\end{figure}

{\Rb{Let $\mu_1,\dots,\mu_{10} \in \mathcal P_{\text{train}}$ be $M=10$ parameters equidistantly distributed (including the endpoints) in $\mathcal P_{\text{train}}=[0.1,1]$. For this parameter extrapolation study, we build a training dataset by integrating the nonlinear Hamiltonian FOM~\eqref{eq:wave_2d_FOM} for each training parameter with the implicit midpoint method until final time $T=8$. We use a fixed time step of $\Delta t=0.1$. In this study, we do not shift the trajectory snapshot data, i.e., $\bq_{\text{ref}}=\bp_{\text{ref}}=\bzero$. From this training dataset, we first construct data-driven nonlinear approximation functions, i.e., $\bgammaQMCL$~\eqref{eq:quad_mnf_co_tan} and $\bgammaBQ$~\eqref{eq:gammqp}. We then derive the corresponding nonlinear SMG-QMCL-ROM~\eqref{eq:jac_based_nonlin_Ham} and nonlinear Galerkin-BQ-ROM~\eqref{pod_based_nonlinear_Ham} for different reduced dimensions. {\Rb{For this parameter extrapolation study, we found $\gamma_q=\gamma_p=10^{-1}$ to be a robust choice for both approaches.}} We also consider $M_{\text{test}}=4$ test parameters $\mu_{\text{test},1}=1.25$, $\mu_{\text{test},2}=1.5$, $\mu_{\text{test},3}=2$, and $\mu_{\text{test},4}=3$ to evaluate how the proposed nonlinear data-driven ROMs generalize for parameter values outside the training dataset.

In Figure~\ref{fig:NW_proj}, we compare the relative projection error in $\bq$ and $\bp$ for the training data consisting of ten trajectories. For both $\bq$ and $\bp$, the nonlinear approximations based on data-driven quadratic manifolds yield a substantially higher accuracy compared to the linear symplectic subspaces for all reduced dimensions. Similar to the linear wave equation example, we use the same quadratic approximations for $\bq$ in both the SMG-QMCL and the Galerkin-BQ approach, and therefore the relative projection error for $\bq$ in Figure~\ref{fig:NW_proj_q} is the same. In Figure~\ref{fig:NW_proj_p}, we observe that the BQ approximation mapping yields a lower projection error than the QMCL approximation for $\bp$.

In Figure~\ref{fig:NW_state}, we observe that the proposed SMG-QMCL and Galerkin-BQ approaches provide ROMs with lower relative state error than the linear symplectic ROMs for both training and test parameters. The comparison plots for the training parameters are shown in Figure~\ref{fig:NW_state_train} where we observe that the proposed approaches achieve substantially lower state error than the symplectic ROMs based on linear symplectic subspaces. The comparison for the test parameters in Figure~\ref{fig:NW_state_test} shows that the ROMs based on SMG-QMCL and Galerkin-BQ approaches achieve two to four times lower relative state error than the ROMs based on the linear symplectic subspace approach. 

The Hamiltonian error~\eqref{eq:err_hamiltonian_non} plots in Figure~\ref{fig:NW_energy} show that the approximate solutions obtained using SMG-QMCL and Galerkin-BQ ROMs of different reduced dimensions  accurately approximate the FOM Hamiltonian for both $\mu_{\text{test},1}=1.25$ and $\mu_{\text{test},4}=3$. For $\mu_{\text{test},1}=1.25$, the error in the Hamiltonian in Figure~\ref{fig:NW_energy_mu1} for the proposed approaches remain below $10^{-2}$. In Figure~\ref{fig:NW_energy_mu2}, we observe that the energy error for the Hamiltonian ROMs based on the SMG-QMCL approach level off at approximately $5 \times 10^{-3}$. The approximately Hamiltonian ROMs based on the Galerkin-BQ approach, on the other hand, exhibit oscillatory error behavior with the error in the Hamiltonian slowly growing with time. 

Finally, we study the accuracy of the proposed ROMs over the two-dimensional computational domain by comparing the pointwise error in $\bq$ between the FOM solution and the reconstructed solution $\bGamma_{\bq}(\widetilde{\bQ}(\mu))$ of the ROM solution $\widetilde{\bQ}(\mu)$ for two test parameter values. The time-evolution of the pointwise error in $\bq$ for the nonlinear Hamiltonian ROM of size $2r=48$ based on the linear symplectic subspace is compared with the corresponding errors for the SMG-QMCL-ROM and the Galerkin-BQ-ROM in Figure~\ref{fig:NW_gam_11} and Figure~\ref{fig:NW_gam_14} for $\mu_{\text{test},1}=1.25$ and $\mu_{\text{test},4}=3$, respectively. For $\mu_{\text{test},1}=1.25$, we observe that the nonlinear ROMs based on SMG-QMCL and Galerkin-BQ approaches yield reconstructed solution fields for $\bq$ with a maximal pointwise error of approximately $10^{-3}$ at all three time instances. In contrast, we observe that the maximal pointwise error in $\bq$ for the linear symplectic subspace approach increases from approximately $5 \times 10^{-3}$ to $2 \times 10^{-2}$ as we march forward in time. For $\mu_{\text{test},4}=3$, we observe that both the SMG-QMCL and the Galerkin-BQ ROM exhibit a maximal pointwise error of approximately $10^{-2}$ whereas the  linear symplectic subspace approach exhibits a maximal pointwise error of approximately $10^{-1}$ which is an order of magnitude higher than the error for proposed approaches. These results demonstrate the ability of both the SMG-QMCL and the Galerkin-BQ approach to provide accurate numerical solutions for multi-dimensional nonlinear problems, even for parameter values outside the range of training parameters. }}
\section{Conclusions \& future work} \label{sec:conclusions}
We have presented two projection-based model reduction approaches that use data-driven quadratic manifolds to derive accurate structure-preserving reduced-order models of canonical Hamiltonian systems. The SMG-QMCL approach derives Hamiltonian ROMs using the quadratic manifold cotangent lift mapping which is based on a quadratic manifold approximation for the generalized position vector. More generally, we proposed the manifold cotangent lift approximation mapping which is a symplectic map and can thus be used with the SMG projection to build a symplectic SMG-MCL-ROM. Note, that this construction does not require the approximation to be a quadratic mapping, as in the QMCL, so the SMG-MCL-ROM could be used to formulate symplectic MOR with other nonlinear approximation mappings like e.g., higher-order polynomials or autoencoders in future work. The Galerkin-BQ approach, on the other hand, augments the linear symplectic subspace with quadratic Kronecker product terms to derive approximately Hamiltonian ROMs. Both approaches are offline--online separable {\Rc{for linear Hamiltonian systems}} which means that the ROMs of both approaches are independent of the FOM dimension. The Galerkin-BQ-ROM approach is particularly attractive from a computational efficiency standpoint due to the ROMs admitting a convenient linear-quadratic model structure with computational complexity on the order of $\bigO(r^4)$ compared to $\bigO(r^5)$ for the SMG-QMCL-ROM. On the other hand, the SMG-QMCL-ROM is favorable in terms of long-term energy conservation. Both of these novel approaches together constitute a first step towards the model reduction of dynamical systems on nonlinear manifolds using interpretable (e.g., polynomial) manifold constructions that ensure that the approximate solution satisfies key physical properties as dictated by the original high-dimensional problem.
    
The numerical experiments with the parametrized linear {\Rc{and nonlinear}} wave equations demonstrate the advantages of using data-driven quadratic manifold approximations for structure-preserving model reduction of transport-dominated problems. The numerical results also show that the proposed approaches produce stable ROMs with higher accuracy and better predictive capabilities than their linear counterparts. Notably, the proposed methods yield accurate predictions of the high-dimensional state even for parameter values that were not included in the training data.
    
This work has opened a number of avenues for future work. {\Rb{Improving the computational efficiency of the proposed nonlinear ROMs based on data-driven quadratic manifolds is necessary. To reduce the online computational cost, the proposed approaches could be combined with the very recently developed structure-preserving DEIM~\cite{pagliantini2022gradient} to improve the computational efficiency of the nonlinear SMG-QMCL and Galerkin-BQ ROMs.}} {\Rc{To further reduce the computational cost}}, alternative formulations for learning manifold approximations in high-dimensional state spaces could be exploited to formulate more suitable symplectic mappings. {\Rb{In another direction, the proposed methodologies could be extended to noncanonical Hamiltonian systems with state-dependent and degenerate forms by building on the recent linear subspace approaches~\cite{hesthaven2021structure}.}} Finally, we seek to broaden the application of the presented framework by learning Hamiltonian dynamics directly from time-domain simulation data in a non-intrusive fashion using operator inference methods \cite{PEHERSTORFER2016196, sharma2022hamiltonian,gruber2023canonical}.
    
\section*{Acknowledgments}  H.S.\ and B.K.\ were in part financially supported by the U.S.\ Office of Naval Research under award number N00014-22-1-2624 and the Ministry of Trade, Industry and Energy (MOTIE) and the Korea Institute for Advancement of Technology (KIAT) through the International Cooperative R\&D program (No.~P0019804, Digital twin based intelligent unmanned facility inspection solutions). R.G.\ has been supported in part by the U.S.\ Department of Energy AEOLUS MMICC center under award DE-SC0019303, program manager W.\ Spotz, and by the AFOSR MURI on physics-based machine learning, award FA9550-21-1-0084, program manager F.\ Fahroo.
P.B.\ is funded by Deutsche Forschungsgemeinschaft (DFG, German Research Foundation) under Germany's Excellence Strategy - EXC 2075 - 390740016.
P.B.\ acknowledges the support by the Stuttgart Center for Simulation Science (SimTech).

\bibliography{references}
\bibliographystyle{abbrv}

%
\end{document}